\documentclass[11pt,letterpaper]{amsart}

\usepackage[OT2,T1]{fontenc}
\usepackage{indentfirst}
\usepackage{amstext}
\usepackage{amsthm}
\usepackage{amsopn}
\usepackage{amsfonts}
\usepackage{amsmath}
\usepackage{latexsym}
\usepackage[francais,english]{babel}
\usepackage{amscd}
\usepackage{amssymb}
\usepackage{amsmath}
\usepackage[all,cmtip]{xy}
\usepackage{leftidx}
\usepackage{graphicx}
\usepackage{etoolbox}
\patchcmd{\section}{\normalfont\scshape\centering}{\normalfont\bfseries}{}{}
\patchcmd{\subsection}{-.5em}{.5em}{}{}

\newtheorem{theo}{{Theorem}}[section]
\newtheorem{coro}[theo]{{Corollary}}
\newtheorem{lemma}[theo]{{Lemma}}
\newtheorem{prop}[theo]{Proposition}

\theoremstyle{definition}
\newtheorem{remark}[theo]{\textbf{Remark}}
\newtheorem{defn}[theo]{Definition}
\newtheorem{example}[theo]{Example}
\numberwithin{equation}{section}

\newtheorem{notation}[theo]{Notation}
\newtheorem{question}[theo]{Question}

\newcommand{\rad}{\mathrm{rad}}

\DeclareFontEncoding{OT2}{}{} 

\begin{document}
\tolerance 400 \pretolerance 200 \selectlanguage{english}

\title{Automorphisms of K3 surfaces, signatures, and isometries of lattices}
\author{Eva  Bayer-Fluckiger}
\date{\today}
\maketitle

\begin{abstract} Let $\alpha$ be a Salem number of degree $d$ with $4 \leqslant d \leqslant 18$. 
We show that if $d  \equiv 0, 4 \ {\rm or} \  6 \ {\rm (mod \ 8)}$, then $\alpha$ is the dynamical degree of
an automorphism of a complex (non-projective) $K3$ surface. We define a notion of {\it signature} of an automorphism, and
use it to give a criterion for Salem numbers of degree 10 and 18 to be realized as  the dynamical degree of such an automorphism.
The first part of the paper contains results on isometries of lattices.
\end{abstract}

\small{} \normalsize

\medskip

\selectlanguage{english}
\section{General introduction}

A {\it Salem polynomial} is a monic, irreducible polynomial $S \in {\bf Z}[X]$ such that 
$S(X) = X^{{\rm deg}(S)} S(X^{-1})$ and that $S$ has exactly two roots outside the unit circle, both positive
real numbers; hence ${\rm deg}(S)$ is even, and $\geqslant 2$; the unique real root $> 1$ of a Salem polynomial is called a {\it Salem number}.
The {\it degree} of the Salem number $\alpha$ is by definition the degree of the  Salem polynomial $S$, 

\medskip

Let $\mathcal X$ be a complex analytic $K3$ surface, and let $T : \mathcal X \to \mathcal X$ be an automorphism;
it induces an isomorphism $T^* : H^2(\mathcal X,{\bf Z}) \to H^2(\mathcal X,{\bf Z})$.  The characteristic polynomial of $T^*$ is a product of at most one Salem polynomial, and of a product of finitely many cyclotomic polynomials
(cf. \cite{Mc1}, Theorem 3.2. and Corollary 3.3).
The {\it dynamical degree} of
$T$ is by definition the spectral radius of $T^*$, i.e. the maximum of the absolute values of the eigenvalues of $T^*$; hence it is either 1 or a Salem number of degree $\leqslant 22$.
We denote by $\lambda(T)$ the dynamical degree of $T$.

\medskip McMullen raised the following question (see \cite {Mc1}, \cite{Mc2}, \cite{Mc3}~) 

\medskip
{\bf Question 1.} 
{\it Which Salem numbers occur as dynamical degrees of automorphisms of complex $K3$ surfaces ?}

\medskip

We say that a Salem number is {\it $K3$-realizable} (or realizable, for short)  if it is the dynamical
degree of an automorphism of a complex $K3$ surface. McMullen (see \cite {Mc1},  \cite{Mc2}, \cite{Mc3})
gave several examples of such Salem numbers (see also
Oguiso \cite {O}, Brandhorst-Gonzalez \cite {BG}, Iwasaki-Takada \cite {IT}), as well as an
{\it infinite} family of realizable degree $6$ Salem numbers; these are obtained by automorphisms of Kummer
surfaces (see \cite{Mc1}, \S 4). Other infinite families of realizable Salem numbers are given 
by Gross-McMullen (in degree $22$, see \cite{GM}, Theorems 1.7 and 1.6), by Reschke
(in degree $14$, see \cite {R}, Theorem 1.2), and Brandhorst proved that some {\it power} of every
Salem number of degree $\leqslant 20$ is realizable (cf. \cite{Br}, Theorem 1.1). The case of Salem numbers of degree $22$ was
treated in \cite {BT}.

\bigskip
We now consider a more precise question, as follows. Let $$H^2(\mathcal X,{\bf C}) = H^{2,0}(\mathcal X) \oplus H^{1,1}(\mathcal X) \oplus H^{0,2}(\mathcal X)$$ be the Hodge decomposition of $H^2(\mathcal X,{\bf C})$.
Since the subspace $H^{2,0}(\mathcal X)$ is one-dimensional, the isomorphism $T^*$ acts
on it by multiplication by a scalar, denoted by $\omega(T)$; we have $|\omega (T)| = 1$. 

\medskip
The value of $\omega(T)$ depends on the geometry of the surface : if 
$\mathcal X$ is projective, then $\omega(T)$ is a root of unity (see \cite {Mc1}, Theorem 3.5).

\medskip
\noindent
{\bf Definition 1.} If $\alpha$ is a Salem number, and if $\delta \in {\bf C}^{\times}$ is such that $|\delta| = 1$, we say that
$(\alpha,\delta)$ is {\it realizable} if there
there exists an automorphism $T$ of a $K3$ surface such that $\lambda(T) = \alpha$ and
$\omega(T) = \delta$.

\medskip Since the characteristic polynomial of $T$ is of the form $SC$, where $S$ is a Salem polynomial and $C$ a product of cyclotomic polynomials, 
$\omega(T)$ is either a root of $S$, or a root of unity. 

\medskip

Suppose first that
$\omega(T)$ is {\it not a root of unity}, i.e. $\lambda(T)$ and $\omega(T)$ are roots of the same Salem polynomial; this implies that the surface is not projective, 
it is of algebraic dimension $0$ (see Proposition \ref{algebraic dimension}). 
Note that this also implies that ${\rm deg}(S) \geqslant 4$.


\medskip
\noindent
{\bf Theorem 1.} {\it Let $S$ be a Salem polynomial of degree $d$ with $4 \leqslant d \leqslant 22$, and let $\alpha$ be the corresponding Salem number. The following
are equivalent :

\medskip
{\rm (1)} There exists a root $\delta$ of $S$ with $|\delta| = 1$ such that $(\alpha,\delta)$ is realizable.

\medskip
{\rm (2)}  $(\alpha,\delta)$ is realizable for all roots $\delta$ of $S$ with $|\delta| = 1$.}

\medskip This is proved in \ref{theorem 1}. For  Salem polynomials of degree $22$, we have the following result (see Theorem \ref{22}) :

\medskip
\noindent
{\bf Theorem 2.} {\it Let $\alpha$ be a Salem number of degree $22$ with 
Salem polynomial $S$, and let 
$\delta$  be a root of $S$ such that $|\delta| = 1$. Then 
$(\alpha,\delta)$  is realizable if and only if $|S(1)|$ and $|S(-1)|$ are squares.}

\medskip
Theorem 2 can also be deduced from the main result in \cite{BT} and Theorem 1. Indeed, it is proved in \cite{BT} that if
$|S(1)|$ and $|S(-1)|$ are squares, then 
$(\alpha,\delta)$  is realizable for {\it some} root $\delta$ of $S$ such that $|\delta| = 1$,
hence by Theorem 1 the pair $(\alpha,\delta)$ is realizable for all roots $\delta$ of $S$ with $|\delta| = 1$.

\medskip We now consider Salem polynomials $S$  of degree $d$ with $4 \leqslant d \leqslant 20$. 

\medskip
\noindent
{\bf Theorem 3.} {\it Let $\alpha$ be a Salem number of degree $d$ with Salem polynomial $S$; assume that 
$4 \leqslant d \leqslant 20$ and 
$d  \equiv 0, 4 \ {\rm or} \  6 \ {\rm (mod \ 8)}$. 
Let
$\delta$  be a root of $S$ such that $|\delta| = 1$.
Then $(\alpha,\delta)$  is realizable. }

\medskip
For $d = 20$, this is Theorem 1.3 in Takada's paper \cite{T}, for the other values of $d$, it  is a consequence of
Theorem \ref{congr 4}.




\medskip
If $d = 10$ or $18$, we give a condition for $(\alpha,\delta)$
to be realizable, in terms of signature maps (see Theorem \ref{equivalent conditions}). This condition
does not always hold when $d = 18$, leading to non-realizable pairs
 (see 
Example  \ref{second smallest}). On the other hand, the case $d = 10$ remains open; numerical 
evidence seems to indicate that all pairs $(\alpha,\delta)$ might be realizable.

 \medskip In a different direction, we show that if $\alpha$ is a Salem number of degree $d$ with 
$4 \leqslant d \leqslant 22$ and Salem polynomial $S$, and 
$\delta$  a root of $S$ such that $|\delta| = 1$, then there exist only finitely many $K3$ surfaces
having an automorphism $T$ with $\lambda(T) = \alpha$ and $\omega(T) = \delta$ (see 
Proposition \ref{finite}). This implies that, if $\alpha$ is a Salem number,  there
exists {\it only finitely many $K3$ surfaces of algebraic dimension $0$} having an
automorphism of dynamical degree $\alpha$. 

\medskip
The last sections of the paper contain some remarks on automorphisms of {\it projective} $K3$ surfaces.

\medskip
The paper is divided into three parts, each preceded by an introduction. The first part concerns {\it isometries of lattices} : these are needed for
the above mentioned results on automorphisms of $K3$ surfaces. Indeed, the intersection form $$H^2(\mathcal X,{\bf Z}) \times H^2(\mathcal X,{\bf Z}) \to {\bf Z}$$ of 
a $K3$ surface $\mathcal X$
is an even unimodular lattice of signature $(3,19)$.
 It is well-known that such a lattice is unique up
to isomorphism; we denote it by $\Lambda_{3,19}$.  Automorphisms of $\mathcal X$ induce
isometries of the intersection form, hence of the lattice $\Lambda_{3,19}$.

\medskip Part III concerns automorphisms of $K3$ surfaces. 
 A criterion of McMullen
(see \cite{Mc2}, Theorem 6.2) plays an important role in the proofs of Theorems 1-3;
he gives explicit conditions for isometries to be induced
by automorphisms of $K3$ surfaces.

\medskip
Part II is the bridge between the arithmetic results of Part I and the geometric applications of Part III. Starting with an isometry of $\Lambda_{3,19}$, one may have to modify it to be able to apply McMullen's criterion; in this process, the characteristic polynomial can change : we keep the same Salem factor, but the cyclotomic factor can be different.

\bigskip
I thank Marie Jos\'e Bertin, Serge Cantat, Curt McMullen, Chris Smyth and Yuta Takada for very useful comments and suggestions. 

\medskip

\bigskip
\centerline {\bf Table of contents}

\bigskip

\noindent
Part I : Isometries of even unimodular lattices \dotfill 3

\medskip
\noindent
Part II :  Salem signature maps and surgery on lattices \dotfill 34

\medskip
\noindent
Part III : Automorphisms of $K3$ surfaces \dotfill 43

\bigskip

\bigskip \centerline {\bf Part I : Isometries of even unimodular lattices}

\bigskip

A {\it lattice} is a pair $(L,q)$, where $L$ is a free ${\bf Z}$-module of finite rank, and $q : L \times L \to {\bf Z}$
is a symmetric bilinear form; it is {\it unimodular} if ${\rm det}(1) = \pm 1$, and {\it even} if $q(x,x)$ is an
even integer for all $x \in L$. 
Let $r,s \geqslant 0$ be integers such that $r  \equiv s \ {\rm (mod \ 8)}$; this congruence condition is equivalent to the
existence of an even unimodular lattice with signature $(r,s)$. When $r,s \geqslant 1$, such a lattice is unique up
to isomorphism (see for instance \cite{S}, chap. V); we denote it by $\Lambda_{r,s}$.
In \cite {GM}, Gross and McMullen raise
the following question (see \cite {GM}, Question 1.1) :

\bigskip
\noindent
{\bf Question 2.} {\it What are the possibilities for the characteristic polynomial $F(X) = {\rm det}(X - t)$ of an isometry
$t \in {\rm SO}(\Lambda_{r,s})$ ? }


\bigskip
The condition $t \in {\rm SO}(\Lambda_{r,s})$ implies that $F(X) = X^{{\rm deg}(F)}F(X^{-1})$,
hence $F$ is  a {\it symmetric} polynomial (cf. \S \ref{basic}). Let
$2n = {\rm deg}(F)$, and let $m(F)$ be the number of roots $z$ of $F$ such that $|z| > 1$. 
As shown in \cite{GM}, we have the further necessary conditions :

\bigskip
(C 1) {\it $|F(1)|$, $|F(-1)|$ and $(-1)^n F(1) F(-1)$ are squares.}

\bigskip

(C 2) {\it $r \geqslant m(F)$, $s \geqslant m(F)$, and if moreover $F(1)F(-1) \not = 0$, then $m(F) \equiv r \equiv s  \ {\rm (mod \ 2)}$.}

\bigskip

Gross and McMullen prove that if $F \in {\bf Z}[X]$ is an irreducible, symmetric and monic
polynomial satisfying condition (C 2) and such that  $|F(1)F(-1)| = 1$, then there exists  $t \in {\rm SO}(\Lambda_{r,s})$ with characteristic polynomial $F$ (see
\cite {GM}, Theorem 1.2). They
speculate that conditions (C 1) and (C 2) are sufficient for a monic {\it irreducible} polynomial to be realized as the characteristic
polynomial of an isometry of $\Lambda_{r,s}$; this is proved in \cite {BT}, Theorem A. More generally, Theorem A of \cite {BT}
implies that 
if a monic, irreducible 
and symmetric polynomial $F$ satisfies conditions (C 1) and (C 2), then {\it there exists an even unimodular lattice}
of signature $(r,s)$ having an isometry with characteristic polynomial $F$. This is also the point of view of the present
paper - we treat the definite and indefinite cases simultaneously. 

\medskip
The first result is that condition (C 1) is sufficient {\it locally} at the finite places.
 If $p$ is a prime number, we say
that a ${\bf Z}_p$-lattice $(L,q)$ is even if $q(x,x) \in 2{\bf Z}_p$ for all $x \in L$; note that if $p \not = 2$, then every lattice is even, since $2$
is a unit in ${\bf Z}_p$. The following is proved in Theorem \ref{local theorem}.

\medskip Let $F \in {\bf Z}[X]$ be a monic, symmetric polynomial of even degree.

\medskip
\noindent
{\bf Theorem 4.} {\it 
Condition {\rm (C 1)} holds if and only if for all prime numbers $p$, there
exists an even unimodular  ${\bf Z}_p$-lattice  having a semi-simple isometry
with characteristic polynomial $F$.}
 
 \medskip If $F$ is irreducible, then Theorem A of \cite {BT} implies that we have a {\it local-global
 principle}; however, this does not extend to the case where $F$ has distinct irreducible factors. 
 We
 define an obstruction to the Hasse principle in terms of an equivalence relation on the irreducible
 factors of $F$; this leads to an {\it obstruction group}, see \S \ref{combinatorial section}.
 
 \medskip Let us write $F(X) = F_1(X)(X-1)^{n^+}(X+1)^{n^-}$ with $F_1(1)F_1(-1) \not = 0$, and suppose for
 simplicity that $n^+ \not = 2$ and $n^- \not = 2$ (the case $n^+ = n^- = 2$ is also treated in the paper, but
 the results are more complicated to state). In this case, the obstruction group only depends on $F$, and
 we denote it by $G_F$. 
  In \S  \ref{sufficient section} we show the following :

 \medskip
 
 \noindent
 {\bf Theorem 5.} {\it Let  $r,s \geqslant 0$ be integers such that $r  \equiv s \ {\rm (mod \ 8)}$. Suppose that Conditions {\rm (C 1)} and {\rm (C 2)} hold and that $G_F = 0$. Let
 $t \in {\rm SO}_{r,s}({\bf R})$ be a semi-simple isometry with characteristic polynomial $F$. Then $t$
 preserves an
 even unimodular
 lattice.}

 \medskip
If $G_F \not = 0$, then we need additional conditions for a semi-simple isometry to stabilize an even unimodular lattice; finding such conditions
is the subject matter of sections  \S \ref{local data} -  \S \ref{sufficient section},  leading to Theorem 6 below :

\medskip Let $t \in {\rm SO}_{r,s}({\bf R})$ be a semi-simple isometry with characteristic polynomial $F$, and suppose that Condition (C 1) holds. In \S \ref{sufficient section}, we
define a homomorphism 
$$\epsilon_{t} : G_F \to {\bf Z}/2{\bf Z},$$
and we show (cf. Theorem \ref{Theo 8}) : 

\medskip
\noindent {\bf Theorem 6.} {\it The isometry $t$ preserves an even unimodular lattice 
if and only if $\epsilon_{t} = 0$. }

\medskip A more convenient way of expressing Theorem 6 is given in terms of {\it signature maps of isometries}, see \S \ref{sign}.

\section{Definitions, notation and basic facts}\label{basic}

We start by recalling some notions and results from \cite {M}, \cite {B 15}, \cite {BT} and \cite {B 21}.
Let $K$ be a field.

\medskip
{\it Isometries}

\medskip

Let $V$ be a finite dimensional
$K$-vector space, and let $q : V \times V \to K$ be a non-degenerate symmetric bilinear form. An {\it isometry} of $(V,q)$ is by definition
an isomorphism $t : V \to V$ such that $q(tx,ty) = q(x,y)$ for all $x, y \in V$. 


\medskip

If
$f \in K[X]$ is a monic polynomial such that $f(0) \not = 0$,  set $$f^*(X) = {f(0)}^{-1} X^{{\rm deg}(f)} f(X^{-1});$$ we say that $f$ is
{\it symmetric} if $f^* = f$.  
It is well-known that the characteristic polynomial of an isometry
is symmetric (see for instance \cite {B 15}, Proposition 1.1).

\medskip
{\it  Equivariant Witt groups and residue maps}

\medskip Let $\Gamma$ be the infinite cyclic group denoted multiplicatively, and
let $K[\Gamma]$ be the associated group ring. Sending $\gamma \in \Gamma$ to $\gamma^{-1}$
induces a $K$-linear involution of $K[\Gamma]$; let $W_{\Gamma}(K)$ be the  Witt
ring  of symmetric $K[\Gamma]$-bilinear forms (see  \cite{BT}, \S 3). If $M$ is a simple $K[\Gamma]$-module,
we denote by $W_{\Gamma}(K,M)$ the subgroup of $W_{\Gamma}(K)$  generated
by the classes of $K[\Gamma]$-bilinear forms $(M,q)$. 

\medskip
Let $V$ be a finite dimensional $K$-vector space, let $q : V \times V \to K$ be a non-degenerate
symmetric bilinear form, and let $t : V \to V$ be an isometry of $q$. Let $\gamma$ be a generator
of $\Gamma$. Let us define a $K[\Gamma]$-module structure on $V$ by setting 
$\gamma.x = t(x)$ for all $x \in V$. This gives rise to a $K[\Gamma]$-bilinear form $(V,q)$, and
we obtain an element $[V] = [V,q]$ of $W_{\Gamma}(K)$.

\medskip If $K$ is a local field with residue field $\kappa$, we have a residue map $W_{\Gamma}(K) \to
W_{\Gamma}(\kappa)$, see \cite{BT}, \S 4; a $K[\Gamma]$-bilinear form contains a unimodular
lattice if and only if it is defined on a bounded module, and its class in $W_{\Gamma}(K)$ is in the kernel of this map
(cf. \cite{BT}, Theorem 4.3 (iv)).

\medskip
{\it Symmetric polynomials over the integers}

\medskip The following definitions and notation are used throughout the paper.

\begin{defn} Let $f \in {\bf Z}[X]$ be a monic, symmetric polynomial. We say that $f$ is of 

\medskip
\noindent
$\bullet$ {\it type} 0 if $f$ is a product of powers of $X-1$ and of $X+1$;

\medskip
\noindent
$\bullet$ {\it type} 1  if $f$ is a product of powers of monic, symmetric, irreducible polynomials in
${\bf Z}[X]$ of even degree;

\medskip
\noindent
$\bullet$ {\it type} 2 if $f$ is a product of polynomials of the form $g g^*$, where $g \in {\bf Z}[X]$ is
monic, irreducible, and $g \not = \pm  g^*$. 

\end{defn}

It is well-known that every monic, symmetric polynomial  is a product of polynomials
of type {\rm 0, 1} and {\rm 2}.

\medskip If $F \in {\bf Z}[X]$ is a monic, symmetric polynomial, we write $F = F_0 F_1 F_2$, where $F_i$ is the product of the irreducible factors of type $i$ of $F$. We have
$F_0(X) = (X-1)^{n^+} (X+1)^{n^-}$ for some integers $n^+,n^- \geqslant 0$, and $F_ 1 = \underset{f \in I_1} \prod f^{n_f}$, where $I_1$ is the set of irreducible factors of type 1 of $F$.

\medskip

Finally, we need the following notation :

\begin{notation}\label{lambda}  Let $E_0$ be an \'etale $K$-algebra of finite rank, and let $E$ be an \'etale $E_0$-algebra 
which is free of rank 2 over $E_0$. Let $\sigma : E \to E$ be the involution fixing $E_0$.

\medskip
(a)  If $\lambda \in E_0^{\times}$, we denote by $b_{\lambda}$
the symmetric bilinear form $b_{\lambda} : E \times E \to K$ such that $b_{\lambda}(x,y) = {\rm Tr}_{E/K}(\lambda x \sigma(y))$. 

\medskip
(b) Set $T(E,\sigma) = E_0^{\times}/{\rm N}_{E/E_0}(E^{\times})$. 

\end{notation}

Note that the isomorphism class of $b_{\lambda}$ only depends on the class of $\lambda$ in
$T(E,\sigma)$.

\section{Local results}\label{local}

Let $F \in {\bf Z}[X]$ be a monic, symmetric polynomial of even degree such that $F(0) = 1$,  and let ${\rm deg}(F) = 2n$. Recall
from the introduction to Part I that the following condition is necessary for the existence of an even,
unimodular lattice having an isometry with characteristic polynomial $F$ :

\bigskip
(C 1) {\it $|F(1)|$, $|F(-1)|$ and $(-1)^n F(1) F(-1)$ are squares.}

\bigskip

The aim of this section is to show that this condition is necessary and sufficient for the existence
of such a lattice {\it at all the finite places} :

\begin{theo}\label{local theorem}
Condition {\rm (C 1)} holds if and only if for all prime numbers $p$, there
exists an even unimodular ${\bf Z}_p$-lattice  having a semi-simple isometry
with characteristic polynomial $F$.

\end{theo}

Let us write $F = F_0 F_1 F_2$ as in \S \ref{basic}, with $F_0(X) = (X-1)^{n^+} (X+1)^{n^-}$.
 Since the hyperbolic lattice of rank ${\rm deg}(F_2)$ has a
semi-simple isometry with characteristic polynomial $F_2$, we assume  that $F = F_0 F_1$.
Let ${\rm deg}(F_1) = 2n_1$.

\medskip
If $p$ is a prime number, we denote by $v_p$ the $p$-adic valuation.

\begin{notation}\label{w2} If $q$ is a non-degenerate quadratic form over ${\bf Q}_p$, we
denote by $w_2(q) \in {\rm Br}_2({\bf Q}_p)$ its Hasse-Witt invariant; we identify $ {\rm Br}_2({\bf Q}_p)$
with $\{0,1\}$. 

\end{notation}

\begin{notation} If $f \in {\bf Z}[X]$ is an irreducible, symmetric polynomial of even degree, set $E_f = {\bf Q}[X]/(f)$, let $\sigma_f : E_f \to E_f$ be the involution
induced by $X \mapsto X^{-1}$, and let $(E_f)_0$ be the fixed field of $\sigma$ in $E_f$. Let $\alpha \in E_f$ be the image of $X$.

\end{notation}

\begin{defn}\label{ramified} Let $f \in {\bf Z}[X]$ be an irreducible, symmetric polynomial of even degree, and let $p$ be a prime number. We say that
$f$ is {\it ramified at} $p$ if there exists a place $w$ of $(E_f)_0$ above $p$ that is ramified in 
$E_f$; otherwise, we say that $f$ is {\it unramified} at $p$. We denote by $\Pi^r_f$ the set of prime numbers $p$ such that $f$ is ramified at $p$.

\medskip Let $p \in \Pi^r_f$ be an odd prime number.  If $w$ is a place of $(E_f)_0$ above $p$ that is ramified in $E_f$, we denote by $\kappa_w$  the
residue field of $w$, and by $\overline \alpha$ be the image of $\alpha$ in $\kappa_w$; we denote by $S_+$ the set of places $w$ above $p$ such that $\overline \alpha = 1$,
and by $S_-$ the set of places $w$ above $p$ such that $\overline \alpha = - 1$. We denote by 
 $\Pi^{r,+}_f$ the set of prime numbers $p$ such that there exists a place $w$ above $p$ with $w \in S_+$, and by $\Pi^{r,-}_f$ the set of 
 prime numbers $p$ such that there exists a place $w$ above $p$ with $w \in S_-$.

\end{defn}

\begin{defn}\label{unramified}  Let $F \in {\bf Z}[X]$  be a monic, symmetric polynomial, and let $I_1$ be
the set of monic, symmetric, irreducible factors of $F$ of even degree.  If $p$ is a prime number, we say that $F$ is {\it ramified at} $p$ if there exists $f \in I_1$ such that $f$ is ramified 
at $p$; otherwise, we say that $F$ is {\it unramified} at $p$. We say that $F$ is {\it unramified} if $F$ is unramified at every prime number $p$; otherwise, we
say that $F$ is {\it ramified}. 

\end{defn}

If $|F(1)F(-1)| = 1$, then $F$ is unramified (see \cite{GM}, Proposition 3.1; note that in \cite{GM} an irreducible, symmetric polynomial $S$ is said to be
unramified if $|S(1)S(-1)| = 1$).

\begin{prop}\label{p} Let $p$ be a prime number. There exists an even,  unimodular ${\bf Z}_p$-lattice having a semi-simple isometry
with characteristic polynomial $F$ if and only if the following conditions hold :

\medskip
{\rm (a)} If $v_p(F_1(1))  \equiv 1 \ {\rm (mod \ 2)}$, then $n^+ \geqslant 1$;

\medskip
{\rm (b)} If $v_p(F_1(-1))  \equiv 1 \ {\rm (mod \ 2)}$, then $n^- \geqslant 1$;

\medskip
{\rm (c)} If moreover $p = 2$ and $n^+ = n^- = 0$, then the class of $(-1)^{n_1}F_1(1)F_1(-1)$ in ${\bf Q}_2^{\times}/{\bf Q_2}^{\times 2}$ lies in $\{1,-3 \}$.

\end{prop}

\noindent
{\bf Proof.} We start by showing that if  conditions (a), (b) and (c) hold, then there exists
an even,  unimodular ${\bf Z}_p$-lattice having a semi-simple isometry
with characteristic polynomial $F$. 

\medskip We use the notation of \S \ref{basic}; in particular, $F_1 = \underset{f \in I_1} \prod f^{n_f}$.
For all $f \in I_1$, we consider the \'etale $\bf Q$-algebra  $E_f = {\bf Q}[X]/(f)$; since $f$ is irreducible,
$E_f$ is a finite field extension of $\bf Q$. We denote by 
$\sigma : E_f \to E_f$ the involution induced by $X \mapsto X^{-1}$, and by $(E_f)_0$ the fixed field
of $\sigma_f$ in $E_f$. Let $L_f$ be a field extension of degree $n_f$ of $(E_f)_0$, linearly disjoint
from $E_f$ over $(E_f)_0$; set $\tilde E_f = E_f \otimes_{(E_f)_0} L_f$, and let
$\tilde \sigma_f : \tilde E_f \to  \tilde E_f$ be the involution induced by $\sigma_f$. Let $\alpha_f$ be a root of $f$
in $E_f$. 

\medskip
Set $E = \underset{f \in I_1} \prod \tilde {E_f}$,  $E_0 = \underset{f \in I_1} \prod L_f$, and let $\sigma = (\tilde \sigma_f)$ be the involution $\sigma : E \to E$ induced by the involutions $\tilde \sigma_f : \tilde E_f \to  \tilde E_f$; we have $E_0 = \{e \in E \ | \ \sigma(e) = e \}$. Let $\alpha = (\alpha_f)$; the characteristic polynomial of the multiplication by 
$\alpha$ on $E$ is equal to $F_1$. 

\medskip
Let $V^+$ and $V^-$ be ${\bf Q}_p$-vector spaces with ${\rm dim}(V^+) = n^+$,  ${\rm dim}(V^-) = n^-$,
and set $V = V^+ \oplus V^-$. Let $\epsilon = (\epsilon^+,\epsilon^-) : V \to V$ be the
isomorphism given by $\epsilon^+ : V^+ \to V^+$, $\epsilon^- : V^- \to V^-$, such that 
$\epsilon^+ = id$ and $\epsilon^- = -id$. 

\medskip
We denote by $\partial : W_{\Gamma}({\bf Q}_p) \to 
W_{\Gamma}({\bf F}_p)$ the residue map (see \S \ref{basic}). If $\lambda \in  E_0^{\times}$, we 
obtain a non-degenerate quadratic form $b_{\lambda} : E \times E \to {\bf Q}_p$, cf. notation  \ref{lambda}.

\medskip The rest of the proof is somewhat different for $p \not = 2$ and $p = 2$. We first consider
the case where $p \not = 2$ :

\medskip
\noindent {\it Claim} 1 : Assume that $p \not = 2$. Let $u_+,u_ -  \in {\bf Z}_p^{\times}$. There exists $\lambda \in E_0^{\times}$ and non-degenerate quadratic forms 
$q^+ : V^+ \times V^+ \to {\bf Q}_p$,  $q^- : V^- \times V^- \to {\bf Q}_p$ such that for $q = q^+ \oplus q^-$, we
have 
$$\partial[E \oplus V, b_{\lambda} \oplus q, \alpha \oplus \epsilon] = 0,$$
and, if moreover $n^+ \not = 0$, then
${\rm det}(q^+) = u_+ F_1(1)$, and if $n^- \not = 0$, then
 $ \ {\rm det}(q^-) = u_ - F_1(-1) .$

\medskip
\noindent {\it Proof of Claim} 1 :  If $v_p(F_1(1))  \equiv 0 \ {\rm (mod \ 2)}$ and
$v_p(F_1(-1))  \equiv 0 \ {\rm (mod \ 2)}$, then this is an immediate consequence of  \cite{BT}, Proposition 7.1.
The algebra $E_0$ decomposes as
a product of fields $E_0 = \underset {v \in S} \prod E_{0,v}$, where $S$ is the set of places above $p$. 
For all $v \in S$, 
set $E_v = E \otimes_{E_0} E_{0,v}$. 
The algebra $E_v$ is of one of
the following types :

\medskip
(sp) $E_v = E_{0,v} \times E_{0,v}$;
 
\medskip
(un) $E_v$ is an unramified extension of $E_{0,v}$;

\medskip
(+) $E_v$ is a ramified extension of $E_{0,v}$, and the image $\overline \alpha$ of $\alpha$ in the residue field $\kappa_v$ of $E_v$ is $1$;

\medskip
(--) $E_v$ is a ramified extension of $E_{0,v}$, and the image $\overline \alpha$ of $\alpha$ in the residue field $\kappa_v$ of $E_v$ is $-1$.

\medskip
This gives a partition $S = S_{sp} \cup S_{un} \cup S_+ \cup S_-$.

\medskip
Let $\gamma$ be a generator of $\Gamma$, and let $\chi_{\pm} : \Gamma \to \{ \pm \ 1\}$ be the character sending $\gamma$ to $\pm 1$. 

\medskip 
Let us choose $\lambda = (\lambda_v)_{v \in S}$ in $E_0^{\times} = \underset {v \in S} \prod E_{0,v}^{\times}$ such that 
for every $v \in S_{un}$, we have $\partial [E_v,b_{\lambda_v},\alpha] = 0$ in $W_{\Gamma}({\bf F}_p)$; this is possible by \cite {BT}, Proposition 6.4. 
For all $u \in {\bf Z}_p^{\times}$, let $\overline u$ be the image of $u$ in ${\bf F}_p$. 

\medskip
Assume first that  $v_p(F_1(1))  \equiv 1 \ {\rm (mod \ 2)}$ and
$v_p(F_1(-1))  \equiv 0 \ {\rm (mod \ 2)}$. For $v \in  S_-$  we choose $\lambda_v$ such that

\medskip

$\underset {v \in S_-} \sum \partial [E_v,b_{\lambda_v},\alpha] = 0$ in $W({\bf F}_p) = W({\bf F}_p,\chi_-) \subset W_{\Gamma}({\bf F}_p).$ This
is possible by \cite {BT}, Proposition 6.6; indeed,  by \cite {BT}, Lemma 6.8 (ii) we have  $$\underset {v \in S_-} \sum [\kappa_v : {\bf F}_p] \equiv \ v_p(F_1(-1)) \  ({\rm mod} \  2).$$ 

\medskip
We now come to the places in $S_+$. Recall that by \cite {BT}, Lemma 6.8 (i) we have $\underset {v \in S_+} \sum [\kappa_v : {\bf F}_p] \equiv \ v_p(F_1(1)) \  ({\rm mod} \  2) $. Since $v_p(F_1(1))  \equiv 1 \ {\rm (mod \ 2)}$ by hypothesis, this implies that $\underset {v \in S_+} \sum [\kappa_v : {\bf F}_p] \equiv \ 1 \  ({\rm mod} \  2) $. 
Therefore there exists $w \in S_+$ such that $[\kappa_w:{\bf F}_p]$ is odd. By \cite {BT}, Proposition 6.6, we can choose $\lambda_w$ such that
$\partial [E_w,b_{\lambda_w},\alpha]$ is either one of the two classes of $W({\bf F}_p) = W_{\Gamma}({\bf F}_p,\chi_+) \subset W_{\Gamma}({\bf F}_p)$ represented by a one-dimensional quadratic form over ${\bf F}_p$.
Let us choose the class of determinant $- \overline {u_+}$, and set $\partial [E_w,b_{\lambda_w},\alpha] = \delta$. 
Since $v_p(F_1(1))  \equiv 1 \ {\rm (mod \ 2)}$, by hypothesis we have $n^+ \geqslant 1$. Let $(V^+,q^+)$ be a non-degenerate 
quadratic form over ${\bf Q}_p$ such that ${\rm det}(q^+) = u_+ F_1(1)$, and that 
{$\partial[V^+,q^+,id] = - \delta$ in $W({\bf F}_p) = W({\bf F}_p,\chi_+) \subset W_{\Gamma}({\bf F}_p).$}

\medskip
Let $S_+' = S_+ - \{w \}$; we have $\underset {v \in S_+'} \sum [\kappa_v : {\bf F}_p] \equiv \ 0 \  ({\rm mod} \  2) $, hence by \cite {BT}, Proposition 6.6,
for all $v \in S_+'$ there exists $\lambda_v \in E_{0,v}^{\times}$ such that 
$\underset {v \in S_+'} \sum \partial [E_v,b_{\lambda_v},\alpha] = 0$ in $W({\bf F}_p) = W({\bf F}_p,\chi_+) \subset W_{\Gamma}({\bf F}_p).$ 
We have
$$\partial [E \oplus V^+, b_{\lambda} \oplus q^+, (\alpha,id)] = 0$$ in $W_{\Gamma}({\bf F}_p)$.
Taking for $(V^-,q^-)$ a  quadratic form over ${\bf Q}_p$ of determinant $u_- F_1(-1)$ and setting $q = q^+ \oplus q^-$, we get
$\partial [E \oplus V, b \oplus q, (\alpha,\epsilon)] = 0$ in $W_{\Gamma}({\bf F}_p).$
This completes the proof when $v_p(F_1(1))  \equiv 1 \ {\rm (mod \ 2)}$ and
$v_p(F_1(-1))  \equiv 0 \ {\rm (mod \ 2)}$. The proof is the same when
$v_p(F_1(1))  \equiv 0 \ {\rm (mod \ 2)}$ and
$v_p(F_1(-1))  \equiv 1 \ {\rm (mod \ 2)}$, exchanging the roles of $S_+$ and $S_-$. 

\medskip
Assume now that $v_p(F_1(1))  \equiv 1 \ {\rm (mod \ 2)}$ and
$v_p(F_1(-1))  \equiv 1 \ {\rm (mod \ 2)}$. By \cite {BT}, Lemma 6.8 (i) and (ii), we have 

\medskip
$\underset {v \in S_+} \sum [\kappa_v : {\bf F}_p] \equiv \ v_p(F_1(1)) \  ({\rm mod} \  2) $, and 
$\underset {v \in S_-} \sum [\kappa_v : {\bf F}_p] \equiv \ v_p(F_1(-1)) \  ({\rm mod} \  2) $.  
Therefore
$\underset {v \in S_+} \sum [\kappa_v : {\bf F}_p] \equiv \ \underset {v \in S_-} \sum [\kappa_v : {\bf F}_p] \equiv \ 1 \  ({\rm mod} \  2) $. Hence
there exist $w_{\pm} \in S_{\pm}$ such that $[\kappa_{w_+}: {\bf F}_p]$  and $[\kappa_{w_-}: {\bf F}_p]$ are odd. By \cite {BT}, Proposition 6.6, we
can choose $\lambda_{w_{\pm}}$ such that $\partial [E_{w_{\pm}},b_{\lambda_{w_{\pm}}},\alpha]$ is either one of the two classes of $\gamma \in W({\bf F}_p) = W(k,\chi_{\pm}) \subset W_{\Gamma}({\bf F}_p)$ 
with ${\rm dim}(\gamma) = 1$. Let us choose  $\lambda_{w_{\pm}}$ such that  $\partial [E_{w_{\pm}},b_{\lambda_{w_{\pm}}},\alpha]$ is
represented by a form of dimension 1 and determinant $\overline u_{\pm}$, and set $\delta_{\pm} = \partial [E_{w_{\pm}},b_{\lambda_{w_{\pm}}},\alpha].$
By hypothesis, we have $n^+ \geqslant 1$ and $n^- \geqslant 1$. Let $(V^{\pm},q^{\pm})$ be non-degenerate 
quadratic forms over ${\bf Q}_p$ such that ${\rm det}(q^{\pm}) = u_{\pm}f(\pm 1)$ and that

\medskip
\centerline {$\partial[V^{\pm},q^{\pm},\epsilon^{\pm}] = - \delta_{\pm}$ in $W({\bf F}_p) = W_{\Gamma}({\bf F}_p,\chi_{\pm}) \subset W_{\Gamma}({\bf F}_p)$.}

\medskip
\medskip
Let $S_+' = S_+ - \{w_+ \}$; we have $\underset {v \in S_+'} \sum [\kappa_v : {\bf F}_p] \equiv \ 0 \  ({\rm mod} \  2) $, hence by \cite {BT}, Proposition 6.6,
for all $v \in S_+'$ there exists $\lambda_v \in E_{0,v}^{\times}$ such that 
$\underset {v \in S_+'} \sum \partial [E_v,b_{\lambda_v},\alpha] = 0$ in $W({\bf F}_p) = W({\bf F}_p,\chi_+) \subset W_{\Gamma}({\bf F}_p).$ 
Similarly, set $S_-' = S_+ - \{w_- \}$; we have $\underset {v \in S_-'} \sum [\kappa_v : {\bf F}_p] \equiv \ 0 \  ({\rm mod} \  2) $, hence by \cite {BT}, Proposition 6.6,
for all $v \in S_-$ there exists $\lambda_v \in E_{0,v}^{\times}$ such that 
$\underset {v \in S_-'} \sum \partial [E_v,b_{\lambda_v},\alpha] = 0$ in $W({\bf F}_p) = W_{\Gamma}({\bf F}_p\chi_+)$ 
Set $q = q^+ \oplus q^-$,
and note that $\partial [E \oplus V, b_{\lambda} \oplus q, (\alpha,\epsilon)] = 0$ in $W_{\Gamma}({\bf F}_p)$.
This completes the proof of Claim 1. We now treat the case where $p = 2$. 

\medskip
\noindent {\it Claim} 2 : Let $u_+,u_ -  \in {\bf Z}_2^{\times}$ be such that  $u_+ u _ - = (-1)^{n}$.  Then there exists $\lambda \in E_0^{\times}$ and non-degenerate quadratic forms 
$q^+ : V^+ \times V^+ \to {\bf Q}_2$,  $q^- : V^- \times V^- \to {\bf Q}_2$ such that for $q = q^+ \oplus q^-$, we
have 
$$\partial[E \oplus V, b_{\lambda} \oplus q, \alpha \oplus \epsilon] = 0,$$
that $(E \oplus V, b_{\lambda} \oplus q)$ contains
an even unimodular ${\bf Z}_2$-lattice, 
and, if moreover $n^+ \not = 0$ then ${\rm det}(q^+) = u_+ F_1(1)$, and if $n^- \not = 0$ then
${\rm det}(q^-) = u_ - F_1(-1).$

\bigskip
\noindent {\it Proof of Claim} 2 :  
The algebra $E_0$ decomposes as
a product of fields $E_0 = \underset {v \in S} \prod E_{0,v}$. 
For all $v \in S$, 
set $E_v = E \otimes_{E_0} E_{0,v}$. 
 The algebra $E_v$ is of one of
the following types :

\medskip
(sp) $E_v = E_{0,v} \times E_{0,v}$;

\medskip
(un) $E_v$ is an unramified extension of $E_{0,v}$;

\medskip
(r) $E_v$ is a ramified extension of $E_{0,v}$.

\medskip
This gives a partition $S = S_{sp} \cup S_{un} \cup S_r$. Assume first that 
$v_2(F_1(1))  +
v_2(F_1(-1))  \equiv 0 \ {\rm (mod \ 2)}$. By \cite {BT}, Lemma 6.8 and Proposition 6.7, we can choose 
$\lambda_v$ such that 

\medskip

\centerline {$\underset {v \in S} \sum \partial [E_v,b_{\lambda_v},\alpha] = 0$ in $W({\bf F}_2) = W_{\Gamma}({\bf F}_2,1) \subset W_{\Gamma}({\bf F}_2).$}

\medskip Therefore
$\partial [E,b_{\lambda}, \alpha] = 0$ in $W_{\Gamma}({\bf F}_2)$.

\medskip

Suppose that $v_2(F_1(1))   \equiv  
v_2(F_1(-1))  \equiv 0 \ {\rm (mod \ 2)}$.
We have ${\rm det}(E,b_{\lambda}) = F_1(1)F_1(-1)$ in 
${\bf Q}_2^{\times}/{\bf Q_2}^{\times 2}$ (see for instance \cite {B 15}, Corollary 5.2). If 
$n^+ = n^- = 0$, then by 
condition (c)
the class of $(-1)^{n_1}F_1(1)F_1(-1)$ in ${\bf Q}_2^{\times}/{\bf Q_2}^{\times 2}$ lies in $\{1,-3 \}$,
and 
this implies that the discriminant of the form $(E,b_{\lambda})$ in ${\bf Q}_2^{\times}/{\bf Q_2}^{\times 2}$ belongs to  $\{1,-3 \}$.
Let us choose $\lambda$ such that the quadratic form
$(E,b_{\lambda})$ contains an even unimodular 
${\bf Z}_2$-lattice. If $S_r = \varnothing$, this is automatic; indeed, in that case every ${\bf Z}_2$-lattice in $(E,b_{\lambda})$  is even. If not, by  \cite {BT} Propositions  8.4 and 5.4 we can choose $\lambda$ having this additional property. 
If
$n^+ \not = 0$ and  $n^- = 0$, we choose $q^+$ such that ${\rm det}(q^+) = (-1)^{n_1} F_1(1) F_1(-1)$; since 
${\rm det}(b_{\lambda}) = F_1(1) F_1(-1)$, this implies that ${\rm det}(E \oplus V,b_{\lambda} \oplus q^+) = (-1)^n.$ Moreover,
let us choose the Hasse-Witt invariant of $b_{\lambda} \oplus q^+$  in such a way that the quadratic form 
$(E \oplus V,b_{\lambda} \oplus q)$ contains an even unimodular 
${\bf Z}_2$-lattice; this is possible by \cite {BT} Corollary 8.4. 
Note that since $v_p({\rm det}(q^+)) = 0$, we have $\partial [V,q^+] = 0$ in $W({\bf F}_2)$,
hence 
$\partial [E \oplus V, b_{\lambda} \oplus q, \alpha \oplus \epsilon] = 0$ in $W_{\Gamma}({\bf F}_2)$.
The same argument applies if  $n^+ = 0$ and  $n^- \not = 0$. 
Suppose that $n^+ \not = 0$ and  $n^- \not = 0$.  Let us choose $q^+$ such that ${\rm det}(q^+) = u_+ F_1(1)$ and
$q^-$ such that ${\rm det}(q^-) = u_-  F_1(-1)$. Since $u_+ u_- = (-1)^{n}$ and 
${\rm det}(b_{\lambda}) =  F_1(1) F_1(-1)$, this implies that ${\rm det}(E \oplus V,b_{\lambda} \oplus q^+ \oplus q^-) = (-1)^{n}.$
As in the previous cases, we can
choose $\lambda$, $q^+$ and $q^-$ such that $(E \oplus V,b_{\lambda} \oplus q)$ contains an even unimodular 
${\bf Z}_2$-lattice, and that $\partial [E \oplus V, b_{\lambda} \oplus q, \alpha \oplus \epsilon] = 0$ in $W_{\Gamma}({\bf F}_2)$.

\medskip
Suppose that $v_2(F_1(1))   \equiv \ 
v_2(F_1(-1))  \equiv 1 \ {\rm (mod \ 2)}$;
note that the hypothesis implies that $n^+, n^- \geqslant 1$.
With our previous choice of $\lambda$, we have $\partial [E,b_{\lambda}, \alpha] = 0$ in $W_{\Gamma}({\bf F}_2)$. Let us choose $q^+$ and $q^-$ such that ${\rm det}(q^{\pm}) = u_{\pm}F_1(\pm 1)$, and 
note that this implies that   ${\rm det}(E \oplus V,b_{\lambda} \oplus q) = (-1)^n$, and that
$v_2({\rm det}(q^-)) = v_2(F_1(-1))$.  Moreover, choose the Hasse-Witt invariants of $b_{\lambda}$, $q^+$ and $q^-$  so
that $(E \oplus V,b_{\lambda} \oplus q^+ \oplus q^-)$ contains an even unimodular 
${\bf Z}_2$-lattice; this is possible by \cite {BT} Corollary 8.4. We have
$\partial [V,q,\epsilon] = 0$ in $W_{\Gamma}({\bf F}_2)$, hence
$\partial [E \oplus V, b_{\lambda} \oplus q, \alpha \oplus \epsilon] = 0$ in $W_{\Gamma}({\bf F}_2)$.

\medskip
Assume now that $v_2(F_1(1))   \equiv  1 \  {\rm (mod \ 2)}$, and
$v_2(F_1(-1))  \equiv 0 \ {\rm (mod \ 2)}$. By hypothesis, this implies that $n^+ \geqslant 1$.
If $n^- \not=  0$, then choose $q^-$ 
such that ${\rm det}(q^{-}) = u_{-}F_1(- 1)$, and note that this implies that $v_2({\rm det}(q^-)) \equiv v_2(F_1(-1))   \equiv  0 \  {\rm (mod \ 2)}$;
choose $q^+$ such that 
${\rm det}(q^{+}) = u_{+}F_1( 1)$. Since $v_2(F_1(1))   \equiv \  1 \  {\rm (mod \ 2)}$, this implies that  $\partial [V^+,q^+,id]$ is the unique non-trivial element of $W({\bf F}_2) = W_{\Gamma}({\bf F}_2,1) \subset W_{\Gamma}({\bf F}_2).$ Note that ${\rm det}(E \oplus V,b_{\lambda} \oplus q^+ \oplus q^-) = (-1)^n$ in  ${\bf Q}_2^{\times}/{\bf Q}^{\times 2}$. Let us choose ${\lambda}$, $q^+$ and $q^-$ such that the quadratic form
$(E \oplus V,b_{\lambda} \oplus q^+ \oplus q^-)$ contains an even unimodular 
${\bf Z}_2$-lattice; this is possible by \cite {BT} Corollary 8.4. Note that
 $\partial [E,b_{\lambda}, \alpha]$ and $\partial [V^+,q^+,id]$ are both equal to the unique non-trivial element of $W({\bf F}_2) = W_{\Gamma}({\bf F}_2,1)$, which is a group of order $2$. Therefore 
$\partial [E \oplus V, b_{\lambda} \oplus q, \alpha \oplus \epsilon] = 0$ in $W_{\Gamma}({\bf F}_2)$.
If $v_2(F_1(1))   \equiv  0 \  {\rm (mod \ 2)}$ and
$v_2(F_1(-1))  \equiv 1 \ {\rm (mod \ 2)}$, the same argument gives the desired result.
This completes the proof of Claim 2.

\medskip
Let us show that conditions (a), (b) and (c) imply the existence of an even unimodular
lattice having a semi-simple isometry with characteristic polynomial $F$. If $p \not = 2$, then
this is an immediate consequence of Claim 1. If $p = 2$, we
apply Claim 2 and 
\cite {BT}, Theorem 8.1; note that $v_2({\rm det}(q^-)) = v_2(F_1(-1))$, hence Theorem 8.5 of
\cite {BT} implies that condition (iii) of this theorem is satisfied.

\medskip
Conversely, we now show that if there exists an even unimodular ${\bf Z}_p$-lattice having a 
semi-simple isometry with characteristic polynomial $F$, then conditions (a), (b) and (c) hold. 
Assume that such a lattice exists, and let $(V,q) = (V^0,q^0) \oplus (V^1,q^1) \oplus (V^2,q^2) = V^0 \oplus V^1 \oplus V^2$ the corresponding
orthogonal decomposition of ${\bf Q}_p[\Gamma]$-bilinear forms.
Since $\partial[V^2] = 0$, we have $\partial[V^0 \oplus V^1] = 0$. We have the further
orthogonal decomposition $V^1 =  \underset {f \in I_1} \oplus V_f$, and $V_f$ is the orthogonal
sum of quadratic forms of the type $b_{\lambda}$. 
By \cite{BT}, Lemma 6.8 and Proposition 6.6, the component of $\partial(V^1)$ in $W_{\Gamma}({\bf F}_p,\chi_{\pm})$ is represented by a quadratic form over ${\bf F}_p$ of dimension $v_p(F_1(\pm 1))$ if $p \not = 2$,
and by a form of dimension 
$v_2(F_1(1)) + v_2 (F(-1))$ if $p = 2$. 

\medskip Suppose that $p \not = 2$, and let us show that (a) and (b) hold. 
Assume that $n^+ = 0$. Then the component
of $\partial [V^0]$ in
$W_{\Gamma}({\bf F}_p,\chi_+)$ is trivial, and therefore $v_p(F_1(1))  \equiv 0 \ {\rm (mod \ 2)}$,
hence (a) holds. By the same argument, (b) holds as well.

\medskip 
Let $p = 2$. We have $v_2 (F_1(-1))  \equiv v_2({\rm det}(q^-))\ {\rm (mod \ 2)}$ by \cite {BT}, Proposition 8.6 and Theorem 8.5. 
If $n^+ = n^- = 0$, then $v_2 (F(-1))  \equiv 0 \ {\rm (mod \ 2)}$, and  the above argument shows that $v_2(F_1(1)) + v_2 (F_1(-1))  \equiv 0 \ {\rm (mod \ 2)}$.
This implies (a) and (b); property (c) also holds, since the discriminant of the lattice is $(-1)^{n_1}F_1(1)F_1(-1)$. If $n^+ \not = 0$ or $n^- \not = 0$, then $v_2(F_1(1)) + v_2 (F_1(-1))  \equiv 1 \ {\rm (mod \ 2)}$;
in both cases, we see that the congruence $v_2 (F_1(-1))  \equiv v_2({\rm det}(q^-))\ {\rm (mod \ 2)}$ implies
(a) and (b). 

\medskip
This completes the proof of the proposition. 

\medskip
\noindent
{\bf Proof of Theorem \ref{local theorem}.} Let $p$ be a prime number. By Proposition \ref{p}, the existence of an  even, unimodular lattice over ${\bf Z}_p$ having a semi-simple
isometry with characteristic polynomial $F$ implies that either $F(1) = 0$, or $v_p(F(1))$ is even; similarly,
either $F(-1) = 0$, or $v_p(F(-1))$ is even. This implies that $|F(1)|$ and $|F(-1)|$ are squares. If
$F(1)F(-1) = 0$, we are done. If not, we apply condition (c)~:  the class of
 $(-1)^{n}F_1(1)F_1(-1)$ in ${\bf Q}_2^{\times}/{\bf Q_2}^{\times 2}$ lies in $\{1,-3 \}$.
Since $|F(1)F(-1)|$ is a square of an integer, this implies that $(-1)^{n}F_1(1)F_1(-1)$ is a square, and hence
condition (C 1) holds. The converse is an immediate consequence of Proposition \ref{p}.

\medskip
The following results and notions will be useful in the next sections :

\begin{coro}\label{for next paper} If there exists an even unimodular lattice having an isometry with characteristic polynomial $F$, then 
the $F$ satisfies Condition {\rm (C 1)}.

\end{coro}

\noindent
{\bf Proof.} This is an immediate consequence of \ref{local theorem}. 

\begin{coro}\label{local coro}
Condition {\rm (C 1)} holds if and only if for all prime numbers $p$, there
exists an even ${\bf Z}_p$-lattice  of determinant $(-1)^n$ having a semi-simple isometry
with characteristic polynomial $F$.

\end{coro}

\noindent
{\bf Proof.} This follows from Theorem \ref{local theorem}, noting that if $(L,q)$ is an even,
unimodular  ${\bf Z}_p$-lattice  having a semi-simple isometry
with characteristic polynomial $F$, then there exists such a lattice of determinant $(-1)^n$.
If $F(1)F(-1) \not = 0$, then ${\rm det}(q) = F(1)F(-1)$ in ${\bf Q}_p^{\times}/{\bf Q}_p^{\times 2}$
(see for instance \cite {B 15}, Corollary 5.2), hence by Condition (C 1) we have
${\rm det}(q) = (-1)^n$ in ${\bf Q}_p^{\times}/{\bf Q}_p^{\times 2}$.
Since ${\rm det}(q)$ is a unit of ${\bf Z}_p$ for all $p$, this
implies that ${\rm det}(q) = (-1)^n$.
If $F(1)F(-1) = 0$, we apply Claim 1 or Claim 2 of the proof of Theorem \ref{local theorem}
to modify the determinant of $q^+$ or $q^-$, if necessary.

\begin{prop}\label{odd ramification} Let $f \in {\bf Z}[X]$  be an irreducible, monic, symmetric polynomial of even degree, and let $p$ be an odd prime number.
If $v_p(f(1))  \equiv 1 \ {\rm (mod \ 2)}$, then $p \in \Pi^{r,+}_f$; if 
$v_p(f(-1))  \equiv 1 \ {\rm (mod \ 2)}$, then $p \in \Pi^{r,-}_f$.

\end{prop}

\noindent
{\bf Proof.} With the notation of the proof of Proposition \ref{p}, Claim 1, we have

\medskip

$\underset {v \in S_+} \sum [\kappa_v : {\bf F}_p] \equiv \ v_p(f(1)) \  ({\rm mod} \  2) $, and 
$\underset {v \in S_-} \sum [\kappa_v : {\bf F}_p] \equiv \ v_p(f(-1)) \  ({\rm mod} \  2) $,

\medskip 
\noindent see \cite{BT}, Lemma 6.8.  This implies the proposition. 

\begin{prop}\label{ramification at 2} Let $f \in {\bf Z}[X]$  be an irreducible, monic, symmetric polynomial of even degree, and suppose that $f$ is unramified at the
prime number $2$. Then 
$v_2(f(1))$ and $v_2(f(-1))$ are even, and 
the class of $(-1)^n f(1) f(-1)$ in ${\bf Q}_2^{\times}/{\bf Q_2}^{\times 2}$ lies in $\{1,-3 \}$.

\end{prop}

\noindent
{\bf Proof.}
Set $E = E_f \otimes_{\bf Q}{\bf Q}_2$ and $E_0 = (E_f)_0 \otimes_{\bf Q}{\bf Q}_2$; let $\alpha$ be the image of $X$ in $E_f$. 
$O_E$ be the ring of integers of $E$, and $O_{E_0}$ the
ring of integers of $E_0$. By hypothesis, the polynomial $f$ is unramified at $2$, therefore all the places of $E_0$ are unramified in $E$. 
The algebra $E_0$ decomposes as a product of fields $\underset{v \in S} \prod E_{0,v}$, and $E = \underset{v \in S} \prod E_v$. 
By \cite{BT}, Proposition 6.4, for all $v \in S$ there exists $\lambda_v \in E_{0,v}^{\times}$ such that 
$\partial[E_v,b_{\lambda_v},\alpha] = 0$, and hence $\underset {v \in S} \sum \partial[E_v,b_{\lambda_v},\alpha] = 0$ in $W_{\Gamma}({\bf F}_2)$.  Setting
$\lambda = (\lambda_v)$, we obtain $\partial_v[E,b_{\lambda},\alpha] = 0$ in $W_{\Gamma}({\bf F}_2)$. By \cite{BT}, Theorem 4.3 (iv) this implies that $(E,b_{\lambda})$
contains a unimodular ${\bf Z_2}$-lattice stabilized by $\alpha$. 

\medskip Since $E/E_0$ is unramified, the trace map $O_E \to O_{E_0}$ is surjective, and therefore every lattice in $(E,b_{\lambda})$ is even (see for instance \cite{BT}, proof of Proposition 9.1). 
Therefore there exists an even unimodular 
${\bf Z}_2$-lattice having an isometry with characteristic polynomial $f$. By Proposition \ref{p}, this implies that $v_2(f(1))$ and $v_2(f(-1))$ are even, are squares. 
and that the class of $(-1)^n f(1) f(-1)$ in ${\bf Q}_2^{\times}/{\bf Q_2}^{\times 2}$ lies in $\{1,-3 \}$.

\begin{prop}\label{even ramification} Let $f \in {\bf Z}[X]$  be an irreducible, monic, symmetric polynomial of even degree, and suppose that $f$ is unramified. Then 
the polynomial $f$ satisfies condition {\rm (C 1)}.

\end{prop}

\noindent
{\bf Proof.} If $p$ is an odd prime number, then Proposition \ref{odd ramification} implies that $v_p(f(1))  \equiv 0 \ {\rm (mod \ 2)}$ and
$v_p(f(-1))  \equiv 0 \ {\rm (mod \ 2)}$. On the other hand, Proposition \ref{ramification at 2} implies that $v_2(f(1))  \equiv 0 \ {\rm (mod \ 2)}$ and
$v_2(f(-1))  \equiv 0 \ {\rm (mod \ 2)}$, and that the class of $(-1)^n f(1) f(-1)$ in ${\bf Q}_2^{\times}/{\bf Q_2}^{\times 2}$ lies in $\{1,-3 \}$.
But $f(1)$ and $f(-1)$ are integers, 
hence  this implies that $(-1)^n f(1) f(-1)$ is a square, and that
$|f(1)|$ and $|f(-1)|$ are squares. 
Therefore Condition (C 1) holds for $f$.




\section{The signature map of an isometry}\label{sign}

The aim of this section is to define a notion of signature for an isometry; 
we start with  reminders concerning signatures of quadratic forms.

\medskip
Let $V$ be a finite dimensional vector space over $\bf R$, and let $q : V \times V \to {\bf R}$ be
a quadratic form; it is the orthogonal sum of a non-degenerate quadratic form and the zero form;
its signature is by definition the signature of the non-degenerate form. The signature of $q$ is
an element of ${\bf N} \times {\bf N}$, denoted by ${\rm sign}(q)$. 
Let $${\rm sum} : {\bf N} \times {\bf N} \to {\bf N}$$ be defined by 
$(n,m) \mapsto n + m$, and 
$${\rm diff} : {\bf N} \times {\bf N} \to {\bf Z}$$  by $(n,m) \mapsto n - m$. Note that ${\rm sum}\circ{\rm sign}(q)$
is the dimension of the non-degenerate part of $q$, and that  ${\rm diff}\circ{\rm sign}(q)$
is the index of $q$, In particular,
the quadratic form $q$ is non-degenerate if
and only if ${\rm sum} \circ {\rm sign} (q) = {\rm dim}(V)$.

\begin{defn}\label{signature definition}
Let $t : V \to V$ be an isometry of $q$. If $f \in {\bf R}[X]$, set $V_f = {\rm Ker}(f(t))$ and let $q_f$ be the restriction of $q$ to $V_f$. 
Let $${\rm sign}_t : {\bf R}[X] \to {\bf N} \times {\bf N}$$
be the map sending $f \in {\bf R}[X]$ to the signature of $(V_f,q_f)$; it is called the {\it signature map
of the isometry} $t$. 

\end{defn}

\begin{lemma}\label{sum} Let $(V_1,q_1)$ and $(V_2,q_2)$ be two non-degenerate quadratic forms over $\bf R$; let $t_1$ be
an isometry of $q_1$, and $t_2$ an isometry of $q_2$. If $(V,q)$ is the orthogonal sum of  $(V_1,q_1)$ and $(V_2,q_2)$ 
and $t = t_1 \oplus t_2$, then ${\rm sign}_t = {\rm sign}_{t_1} + {\rm sign}_{t_2}$. 

\end{lemma}

\noindent
{\bf Proof.} For all $f \in {\bf R}[X]$, we have $V_f  = (V_1)_f \oplus (V_2)_f$ and $q_f = (q_1)_f \oplus (q_2)_f$, hence 
${\rm sign}_t (f) = {\rm sign}_{t_1}(f) + {\rm sign}_{t_2}(f)$. 

\bigskip

The aim of this section is to characterize signature maps of semi-simple isometries; we start with some notation and definitions :

\begin{defn}\label{Sym} Let ${\rm Sym}({\bf R}[X])$ be the set of symmetric polynomials in ${\bf R}[X]$. The 
{\it symmetric radical} of $f \in {\bf R}[X]$ is by definition the monic, symmetric divisor  of highest degree of $f$; we denote it by ${\rad}(f) \in {\bf R}[X]$.

\end{defn}

\begin{notation}\label{nf} If $F \in {\rm Sym}({\bf R}[X])$ and if $f \in {\rm Sym}({\bf R}[X])$ is irreducible, 
we denote by $n_f (F) $ be the largest integer $\geqslant 0$ such that $f^{n_f(F)}$ divides $F$. 

\end{notation}

\begin{notation}\label{proj} For $i=1,2$ let ${\rm proj}_i :  {\bf N} \times {\bf N} \to {\bf N},\; (n_1, n_2) \mapsto n_i$
   be the projection on the $i$-th component.

\medskip
Let ${\rm proj}_1 : {\bf N} \times {\bf N}$ be the projection of the first factor, ${\rm proj}_1 (n,m) = n$, and let 
 ${\rm proj}_2 : {\bf N} \times {\bf N}$ be the projection of the second factor, ${\rm proj}_1 (n,m) = m$.

\end{notation} 

\begin{defn}\label{bounded} We say that a map $\tau :  {\bf R}[X] \to {\bf N} \times {\bf N}$ is {\it bounded} if there
exists $(r,s) \in {\rm Im}(\tau)$ such that ${\rm proj}_1(\tau(f)) \leqslant r$ and ${\rm proj}_2(\tau(f)) \leqslant s$ for all $f \in  {\bf R}[X]$,
and $(r,s)$ is called the {\it maximum} of the bounded map $\tau$.

\end{defn}

\begin{prop}\label{signature characterization}
Let $\tau : {\bf R}[X] \to {\bf N} \times {\bf N}$ be a bounded map with maximum $(r,s)$
with the following properties :

\medskip
{\rm (a)} There exists a unique monic polynomial  $F \in {\rm Sym}({\bf R}[X])$ of degree $r+s$ such
that $\tau(F) = (r,s)$, and such that if $f \in {\bf R}[X]$ is irreducible with $\tau(f) \not = (0,0)$, then $f$ divides $F$.

\medskip
{\rm (b)}  If $f, g \in {\rm Sym}({\bf R}[X])$ are relatively prime polynomials, then 
$$\tau (fg) = \tau (f) + \tau (g).$$

{\rm (c)} 
$\tau(f) = \tau ({\rm rad}(f))$ for all $f \in {\bf R}[X]$ and $\tau(f) = 0$ if $f$ is constant.

\medskip
{\rm (d)} If $f \in {\rm Sym}({\bf R}[X])$ is irreducible, then 
$\tau (f) \in {\rm deg}(f) {\bf N} \times  {\rm deg}(f){\bf N},$

$${\rm sum} \circ \tau (f) = {\rm deg}(f)n_f(F)$$ and $$\tau(f^k) = \tau(f)$$  for all integers $k \geqslant 1$.

\medskip 
{\rm (e)} If $f(X) = (X-a)(X-a^{-1})$ with $a \not = \pm 1$, then 
$$ \tau (f) = (n_{X-a}(F),n_{X-a}(F)),$$
and $\tau(f^k) = \tau(f)$  for all integers $k \geqslant 1$.

\bigskip
Then there exists a non-degenerate quadratic form $q$ over $\bf R$ of signature $(r,s)$ and a semi-simple isometry $t$ of
$q$ such that ${\rm sign}_t = \tau$. 

\bigskip Conversely, if $t$ is a semi-simple isometry of a non-degenerate quadratic form  $q$ of signature $(r,s)$, then $\tau = {\rm sign}_t $
is bounded of maximum $(r,s)$ and 
satisfies conditions {\rm (a)-(e)}.

\end{prop}

\noindent 
{\bf Proof.} Let $F$ be the polynomial associated to $\tau$ as in (a), and let $I_1$ be the set of
irreducible factors of degree $2$ of $F$. For all $f \in I_1$, set $E_f = {\bf R}[X]/(f)$, and let
$\sigma_f : E_f \to E_f$ be the involution induced by $X \mapsto X^{-1}$. Condition (d) implies
that $\tau (f) \in 2{\bf N} \times 2{\bf N}$; let $r_f, s_f \geqslant 0$ be integers such that
$\tau(f) = (2r_f, 2s_f)$. Set $W_f = E_f^{n_f}$, let $h_f : W_f \times W_f \to E_f$ be
the hermitian form of signature $(r_f,s_f)$, and let $q_f : W_f \times W_f \to {\bf R}$ be 
defined by $q_f(x,y) = {\rm Tr}_{E_f/{\bf R}}(h_f(x,y))$; the
signature of this quadratic form is $(2r_f,2s_f)$, and it has a semi-simple isometry  with
characteristic polynomial $f^{n_f}$. 
Set $q_1 = \underset{f \in I_1} \oplus q_f$. Let $r_0,s_0 \geqslant 0$ be integers such
that $\tau(X-1) + \tau(X+1) = (r_0,s_0)$, and let $q_0$ be the non-degenerate quadratic
form over ${\bf R}$ with signature $(r_0,s_0)$. Finally, let $q_2$ be the hyperbolic
quadratic  form over ${\bf R}$ of dimension ${\rm deg}(F) - {\rm dim}(q_1) - {\rm dim}(q_0)$,
and set $q = q_0 \oplus q_1 \oplus q_2$; the quadratic form has a semi-simple isometry
with signature map $\tau$. 

\medskip The converse follows from Milnor's classification of isometries in \cite{M}, \S 3. 
Indeed, let $(V,q)$ be a non-degenerate quadratic form over ${\bf R}$ of signature $(r,s)$,
let $t : V \to V$ be a semi-simple isometry of $q$, 
and let $F \in {\bf R}[X]$ be the characteristic polynomial of $t$. 
Set $\tau = {\rm sign}_t$.

\medskip Let $I_1$ be the
set of irreducible factors of degree 2 of $F$, and let us write $F = F_0F_1F_2$, 
where $F_0(X) = (X-1)^{n_+}(X+1)^{n_-}$, $F_1 = \underset{f \in I_1} \prod f^{n_f(F)}$,
and $F_2$  is a product of polynomials of the form $(X-a)(X - a^{-1})$ with $a \in {\bf R}$
and $a \not = \pm 1$. By \cite {M}, Lemma 3.1, we have an orthogonal decomposition 
$V = V_0 \oplus V_1 \oplus V_2$; moreover, the restriction of $q$ to $V_2$ is hyperbolic
(see \cite {M}, page 94, Case 3). This implies that $\tau$ has properties (a), (b), (c) and (e). 
Property (d) follows from \cite{M}, Theorems 3.3 and 3.4, combined with \cite{M}, \S 1.

\begin{defn} A bounded map $\tau : {\bf R}[X] \to {\bf N} \times {\bf N}$ satisfying conditions
(a)-(e) is called a  (semi-simple) {\it signature map}. The polynomial of condition (a) is by definition the
{\it polynomial associated to} $\tau$. 

\end{defn}

Since we only consider semi-simple isometries in this paper, we use the terminology of {\it signature map}
without the adjective ``semi-simple". 

\medskip
Finally, we note that if $\tau$ is a signature map with polynomial $F$ and maximum $(r,s)$, then
condition (C 2) of Gross and McMullen (see the introduction to Part I) holds for $F$ and $(r,s)$.

\section{Signature maps of automorphisms of $K3$ surfaces}\label{signature automorphism}

Let $T : \mathcal X \to \mathcal X$ be an automorphism of a $K3$ surface. 
The {\it intersection form} of $\mathcal X$ induces a non-degenerate quadratic form
$$H^2(\mathcal X,{\bf R}) \times H^2(\mathcal X,{\bf R}) \to {\bf R}$$ of signature 
$(3,19)$.
 The isomorphism $T^*$ is an isometry of the intersection form, and
the {\it signature map of the automorphism} $T : \mathcal X \to \mathcal X$ is by definition the signature map
of the isometry $T^*$.

\medskip
The signature map of $T$ determines the characteristic polynomial of $T^*$, hence also the dynamical degree $\lambda(T)$. Note that
either $\lambda(T) = 1$, or $\lambda(T)$ is the unique real number $\alpha > 1$
such that $$\rm {sign}_{T*}((X-\alpha)(X-\alpha^{-1})) = (1,1).$$

\medskip In addition, the signature map of $T$ also determines $\omega(T)$, up to complex conjugation. 
Indeed, 
let $\delta$ be a complex number, let $\overline \delta$ be its complex conjugate, and let 
$f_{\delta}\ \in {\bf R}[X]$ be
the minimal polynomial of $\delta$ over $\bf R$. 

\medskip
Recall that ${\rm proj}_1 : {\bf N} \times {\bf N}$ denotes the projection of the first factor, ${\rm proj}_1 (n,m) = n$. 
We have

\medskip

\centerline {$\omega(T) = \delta$ or $\overline \delta$ $\iff$ $\rm {proj}_1 \circ \rm {sign}_{T*}(f_{\delta}) = 2.$}

\medskip
Let $F$ be the characteristic polynomial of $T^*$; we have $F = S C$, where $S$ is a Salem polynomial and $C$ is a product of cyclotomic polynomials; set
$d = {\rm deg}(S)$. We have either $\rm {sign}_{T*}(S) = (3,d-3)$ and, or 
$\rm {sign}_{T*}(S) =(1,d-1)$ and $\rm {sign}_{T*}(C) = (2,20-d)$.

\medskip

Let us  assume that $\lambda(T) > 1$ and that $\omega(T)$ is {\it not a root of unity}; this
implies that $\lambda(T)$ and $\omega(T)$ are roots of the Salem polynomial $S$. Then we are in the first case, $\rm {sign}_{T*}(S) = (3,d-3)$, and
the characterization
of $\lambda(T)$ and $\omega(T)$ in terms of the signature map takes a very simple form.
Let $\alpha \in {\bf R}$, $\delta \in {\bf C}^*$ be such that $\alpha > 1$ and $|\delta| = 1$. We have

\medskip

\centerline {$\lambda(T) = \alpha$  $\iff$ $\rm {sign}_{T*}((X-\alpha)(X-\alpha^{-1})) = (1,1)$.}

\medskip
and 

\centerline {$\omega(T) = \delta$ or $\overline \delta$ $\iff$ $\rm {sign}_{T*}((X-\delta)(X-\delta^{-1})) = (2,0).$}

\bigskip
Conversely, note that under this hypothesis,  the signature map $\rm {sign}_{T*}$ is determined by the characteristic polynomial $F = SC$, and by the choice of
$\delta \in {\bf C}^*$ such that $\rm {sign}_{T*}((X-\delta)(X-\delta^{-1})) = (2,0).$

\section{Local data}\label{local data}

To explain the results of this section, we start with a special case. As in the introduction, let
$S$ be a Salem polynomial and let $C$ be a product of cyclotomic polynomials.
Let $\tau$ be a signature map with associated polynomial $SC$. 

\medskip We consider the following questions - the first one being the central topic of Part I, the
second of Part II :

\medskip
(a) Does there exist an even unimodular lattice having an isometry with signature map $\tau$ ?

\medskip
 Let $t : L \to L$ be an isometry of an
even unimodular lattice $L$ with signature map $\tau$; set $L_S = {\rm Ker}(S(t))$ and $L_C = 
{\rm Ker}(C(t))$. The properties of these lattices 
determine whether or not $t : L \to L$ is induced by an automorphism of a $K3$ surface; this follows
from results of McMullen (see \cite{Mc2}, Theorem 6.2). 
This shows that it is not enough to answer question (a); we also need to understand the flexibility
we have in the choice of the isometry in case of a positive answer. 

\medskip
More precisely,  we ask :

\medskip
(b) Does there exist an even unimodular lattice $L$ having an isometry $t$ with signature map $\tau$,
such that the restriction of $t$ to $L_C$ satisfies McMullen's criterion~?

\medskip
Let $V = L \otimes_{\bf Z}{\bf Q}$, $V_S = L_S \otimes_{\bf Z}{\bf Q}$ and $V_C = L_C \otimes_{\bf Z}{\bf Q}$;
we have $V = V_S \oplus V_C$. As we will see, the ``building blocks"  $V_S$ and $V_C$ play an important
role in both questions (a) and (b).  
The first step is to describe, for all prime numbers $p$, the possibilities for the
quadratic spaces with isometry $V_S \otimes_{\bf Q}{\bf Q}_p$ and $V_C \otimes_{\bf Q}{\bf Q}_p$; 
 this will give rise to {\it combinatorial
data}, as indicated in the introduction to Part I.

\medskip
We now return to the general case. Let $\tau$ be a signature map with associated polynomial $F$,
and assume that $F \in {\bf Z}[X]$.
Let $p$ be a prime number, let $V$ be a finite dimensional ${\bf Q}_p$-vector space, let $q : V \times V
\to {\bf Q}_p$ a non-degenerate quadratic form, and let $t : V \to V$ a semi-simple isometry of the quadratic form $q$ with characteristic polynomial $F$. The aim of this section is to associate combinatorial data 
to this isometry.

\medskip This turns out to be especially simple in an important special case : when at least one of
$F(1)$ and $F(-1)$ is non-zero. Indeed, let 
$F = F_0 F_1$, where $F_i$ is the product of the irreducible factors $f \in {\bf Z}[X]$ of type $i$ of $F$,
with
$F_1 = \underset{f \in I_1} \prod f^{n_f}$ and 
$F_0(X) = (X-1)^{n^+} (X+1)^{n^-}$ for some integers $n^+,n^- \geqslant 0$. 

\medskip
 Set $V_{F_i} = {\rm Ker}(F_i(t))$; we have an
orthogonal decomposition of ${\bf Q}_p[\Gamma]$-quadratic forms $$(V,q)  = (V_{F_0},q_0) \oplus (V_{F_1},q_1).$$

Note that ${\rm det}(q) = {\rm det}(q_0) {\rm det}(q_1)$ in ${\bf Q}_p^{\times}/{\bf Q}_p^{\times 2}$
and that  $ {\rm det}(q_1) = F_1(1)F_1(-1)$. Set $V_f = {\rm Ker}(f(t))$; we  have
the further orthogonal decomposition 
of ${\bf Q}_p[\Gamma]$-quadratic forms $(V_{F_1},q_1) = \underset{f \in I_1}  \oplus (V_f,q_f)$, and 
${\rm det}(q_f) = [f(1)f(-1)]^{n_f}$. 
The dimensions and determinants of the quadratic forms $(V_f,q_f)$ are
determined by the polynomial $f$; these forms can however have different {\it Hasse-Witt invariants},
and this is what will determine the combinatorial data.

\medskip The dimension of $q_0$ is $n_0 = n^+ + n^-$, and its determinant is 

\medskip

\centerline {${\rm det}(q_0) = {\rm det}(q) {\rm det}(q_1) = {\rm det}(q) F_1(1) F_1(-1)$ in ${\bf Q}_p^{\times}/{\bf Q}_p^{\times 2}$,}

\medskip 
\noindent
hence ${\rm det}(q_0)$ only depends on the polynomial $F$ and the determinant of the quadratic form $q$.

\medskip
In the special case where at least one of $F(1)$ or $F(-1)$ is non-zero, then we have $n_0 = n^+$ or
$n_0 = n^-$, and we can define the combinatorial data as a function of $F$ and ${\rm det}(q)$; however,
if $F(1) = F(-1) = 0$, then we have a non-trivial orthogonal decomposition of ${\bf Q}_p[\Gamma]$-quadratic forms 
$(V_0,q_0) = (V_+,q_+) \oplus (V_-,q_-)$, and we only know ${\rm det}(q_0) = {\rm det}(q_+){\rm det}(q_-)$
as a function of $q$ and $F$, but not the factors  ${\rm det}(q_+)$ and ${\rm det}(q_-)$.

\medskip
We deal with this question in \S \ref{local 0 section}. For now, we assume that
at least one of $F(1)$ or $F(-1)$ is non-zero, and hence $n_0 = n^+$ or
$n_0 = n^-$. In the remainder of this section, we define the local data in this case : first, we
deal with the polynomials $f \in I_1$, the hypothesis is not needed for these, and then we pass
to $f \in I_0$, under the above assumption.

\medskip We start by introducing some notation. 

\begin{notation}\label{E}
 For all $f \in I_1$, 
let $Q_f : (E_f )^{n_f} \times (E_f )^{n_f}  \to {\bf Q}$ be the orthogonal sum of $n_f$ copies of the quadratic form $E_f \times E_f \to {\bf Q}$ defined by $(x,y) \mapsto {\rm Tr}_{E_f/{\bf Q}}(x \sigma_f(y))$.

\end{notation}

\medskip

Set $E_f^p = E_f \otimes_{\bf Q} {\bf Q}_p$ and
$(E_f)_0^p = (E_f)_0^p \otimes_{\bf Q}{\bf Q}_p$; there exists a unique hermitian form 
$(V_f,h^p_f)$ over
$(E_f^p,\sigma_f)$ such that $q_f(x,y) = {\rm Tr}_{E_f^p/{\bf Q}_p}(h^p_f(x,y))$ for all $x, y \in V_f$;
see for instance \cite {M}, Lemma 1.1 or \cite {B 15}, Proposition 3.6.

\medskip
Set 
$\lambda^p_f = {\rm det}(h^p_f) \in (E^p_f)_0^{\times}/{\rm N}_{E^p_f/(E^p_f)_0}({E_f^p}^{\times})$. Note that the hermitian
form $h^p_f$ is isomorphic to the $n_f$-dimensional diagonal hermitian form $\langle \lambda^p_f,1,\dots,1 \rangle$
over $E_f^p$; hence $q_f$ is determined by $\lambda^p_f$.

\begin{notation}\label{lambda bis} With the notation above, set 
$$\partial_p(\lambda^p_f) = \partial_p[q_f^p] \in W_{\Gamma}({\bf F}_p).$$
If $\lambda^p = (\lambda_f^p)_{f \in I_1}$, set $\partial_p(\lambda^p) =
\underset{f \in I_1} \oplus \partial_p (\lambda^p_f) =  \partial_p[\underset{f \in I_1} \oplus q_f^p] \in W_{\Gamma}({\bf F}_p).$

\end{notation}

\begin{notation}\label{w2} If $q$ is a non-degenerate quadratic form over ${\bf Q}_p$, we
denote by $w_2(q) \in {\rm Br}_2({\bf Q}_p)$ its Hasse-Witt invariant; we identify $ {\rm Br}_2({\bf Q}_p)$
with $\{0,1\}$. 

\medskip
For all $f \in I_1$, 
let $d_f \in (E_f)_0^{\times}$ be such that $E_f = (E_f)_0 (\sqrt d_f)$. 

\end{notation}

\begin{prop}\label{dim det w f} We have ${\rm dim}(q_f^p) = n_f{\rm deg}(f) = {\rm dim}(Q_f)$, ${\rm det}(q_f^p) = [f(1)f(-1)]^{n_f} = {\rm det}(Q_f)$ in ${\bf Q}_p^{\times}/{\bf Q}_p^{\times 2}$, and
$$w_2(q_f^p) + w_2(Q_f) = {\rm cor}_{(E_f)_0^p/{\bf Q}_p}(\lambda_f^p,d_f)$$ in ${\rm Br}_2({\bf Q}_p)$, where ${\rm cor}$ denotes
the corestriction map. 

\end{prop}

\noindent
{\bf Proof.} The assertion concerning the dimension is clear, the one on the determinant follows from 
 \cite{B 15}, Corollary 5.2, and the property of the Hasse-Witt invariants from \cite {B 21}, Proposition 12.8.

\begin{notation}\label{a 1} For $f \in I_1$, set $a^p(f) = {\rm cor}_{(E_f)_0^p/{\bf Q}_p}(\lambda_f^p,d_f)$
in ${\rm Br}_2({\bf Q}_p)$.

\end{notation}

By Proposition \ref {dim det w f}, we have $a^p(f)  = w_2(q_f^p) + w_2(Q_f)$.

\bigskip We now {\it assume that} $F(1) \not = 0$ {\it or} $F(-1) \not = 0$.

\medskip
If $f \in I_0$, then 
either $f(X) = X-1$  or $f(X) = X + 1$. Set $n_0 = n_f$ and 
note that $F_0(X) = (X-1)^{n_0}$ in the first case, and $F_0(X) = (X+1)^{n_0}$ in the second one. 

\medskip We have ${\rm dim}(V_{F_0}^p) = n_0$ and ${\rm det}(q_0^p) =  {\rm det}(q) F_1(1) F_1(-1)$ in ${\bf Q}_p^{\times}/{\bf Q}_p^{\times 2}$.
 
\medskip Set
 $D_0 =  {\rm det}(q) F_1(1) F_1(-1)$, and let
$Q_0$ be the diagonal quadratic form  of dimension $n_0$ over ${\bf Q}$ defined by
$Q_0 = \langle D_0,1,\dots,1 \rangle$. 

\begin{notation}\label{0} With $f \in I_0$ as above, set $a^p(f) = w_2(q_0^p) + w_2(Q_0)$ in ${\rm Br}_2({\bf Q}_p)$.

\end{notation} 

Set $I = I_0 \cup I_1$. If $f \in I_0$, $a^p(f)$ is defined as in Notation \ref{0}, and if $f \in I_1$, then
$a^p(f)$ is as in Notation \ref{a 1}.
Let us denote by $C(I)$ the set of maps $I \to {\bf Z}/2{\bf Z}$; we obtain 
$a^p \in C(I)$. In summary, we have

\begin{defn}\label{a} Assume that $F(1) \not = 0$ or $F(-1) \not = 0$, and let
 $\lambda^p = (\lambda_f^p)_{f \in I_1}$ and $q_0^p$ be as above. We define 
$$a^p = a^p[\lambda^p,q_0] \in C(I)$$
by setting $a^p(f) = w_2(q_0^p) + w_2(Q_0)$ if $f \in I_0$, and
$a^p(f)  = w_2(q_f^p) + w_2(Q_f) =  {\rm cor}_{(E_f)_0^p/{\bf Q}_p}(\lambda_f^p,d_f)$ if $f \in I_1$. 

\end{defn}

\begin{notation}\label {curly C}

Let  $\delta \in W_{\Gamma}({\bf F}_p)$, and let ${\mathcal C}_{\delta}^p$ be the set of 
$a^p  =  a^p[\lambda^p,q_0] \in C(I)$ as in definition \ref{a} such that $\partial_p(\lambda^p) \oplus \partial_p[q_0]  = \delta$. 

\end{notation}

\section{Local data - continued}\label{local 0 section}

We keep the notation of \S \ref{local data}; in particular, $F \in {\bf Z}[X]$ is
a monic, symmetric polynomial, $p$ is a prime number, $(V,q)$ is
a non-degenerate quadratic form over ${\bf Q}_p$,  and $t : V \to V$ is a semi-simple isometry of $q$ with characteristic
polynomial $F$. Recall that $F = F_0 F_1$, with 

\medskip
\centerline {$F_1 = \underset{f \in I_1} \prod f^{n_f}$ and
$F_0(X) = (X-1)^{n^+} (X+1)^{n^-}$} 
\noindent for some integers $n^+,n^- \geqslant 0$.

\medskip In \S \ref{local data}, we defined a set $\mathcal C^p$ under the hypothesis that $n^+ = 0$
or $n^- = 0$; we now do this in general. Recall that $D_0 = {\rm det}(q)F_1(1)F_1(-1)$, and that
we have the orthogonal decompositions of ${\bf Q}_p[\Gamma]$-quadratic forms $$(V,q)  = (V_{F_0},q_0) \oplus (V_{F_1},q_1) 
$$ and $(V_0,q_0) = (V_+,q_+) \oplus (V_-,q_-)$. We have ${\rm det}(q_0) = D_0$,
but ${\rm det}(q_+)$ and ${\rm det}(q_-)$  are not determined by $F$ and $q$.

\medskip 
Let $D_+, D_- \in {\bf Q}_p^{\times}/{\bf Q}_p^{\times 2}$ be such that $D_0 = D_+ D_-$; the set
$\mathcal C^p$ we now define {\it also depends on the choice of} $D_+$ and $D_-$. 

\medskip Let
$Q_+$ be the diagonal quadratic form  of dimension $n_+$ over ${\bf Q}$ defined by
$Q_+ = \langle D_+,1,\dots,1 \rangle$; similarly, let
$Q_-$ be the diagonal quadratic form  of dimension $n_-$ over ${\bf Q}$ defined by
$Q_- = \langle D_-,1,\dots,1 \rangle$.

\begin{notation}\label{+ -} Set $a^p(X - 1) = w_2(q_+) + w_2(Q_+)$ in ${\rm Br}_2({\bf Q}_p)$, and 
$a^p(X +1) = w_2(q_-) + w_2(Q_-)$ in ${\rm Br}_2({\bf Q}_p)$.

\end{notation} 

This defines $a^p(f)$ for $f \in I_0$, and 
if $f \in I_1$, then
$a^p(f)$ is as in Notation \ref{a 1}. We identify ${\rm Br}_2({\bf Q}_p)$
with $\{0,1 \}$, and obtain $a^p \in C(I)$. In summary, we have

\begin{defn}\label{aa} Let 
let $\lambda^p = (\lambda_f^p)_{f \in I_1}$ be as in \S \ref{local data}. We define 
$$a^p = a^p[\lambda^p,q_{\pm}] \in C(I)$$
by setting $a^p(X -1) = w_2(q_+) + w_2(Q_+)$, $a^p(X +1) = w_2(q_-) + w_2(Q_-)$, and 

\medskip
\centerline {$a^p(f)  = w_2(q_f^p) + w_2(Q_f) =  {\rm cor}_{(E_f)_0^p/{\bf Q}_p}(\lambda_f^p,d_f)$ if $f \in I_1$.}

\end{defn}

\begin{notation}\label {curly C bis}

Let  $\delta \in W_{\Gamma}({\bf F}_p)$, and let ${\mathcal C}_{\delta}^p$ be the set of 
$a^p  =  a^p[\lambda^p,q_{\pm}] \in C(I)$ as in definition \ref{aa} such that $\partial_p(\lambda^p) \oplus \partial_p[q_+]  \oplus \partial_p[q_-] = \delta$. 

\end{notation}




\section{Combinatorial data - obstruction group and equivalence relations}\label{combinatorial section}

The aim of this section is to define some combinatorial
objects that we need in the next sections. We keep the notation of \S \ref{local 0 section} : we fix monic, symmetric polynomial $F \in {\bf Z}[X]$ of even degree.
Let $I$ be the set of irreducible, symmetric factors of $F$, and we write $F$ as $F = F_0 F_1$, with 

\medskip
\centerline {$F_1 = \underset{f \in I_1} \prod f^{n_f}$ and
$F_0(X) = (X-1)^{n^+} (X+1)^{n^-}$} 
\noindent for some integers $n^+,n^- \geqslant 0$.

\medskip If $f \in {\bf Z}[X]$ is a symmetric, irreducible polynomial of even degree, recall that we denote by 
$\Pi^r_f$ the set of prime numbers $p$ such
that $f$ is ramified at $p$ (see Definition \ref{ramified}).

\begin{notation} If $f, g \in {\bf Z}[X]$ are monic, irreducible, symmetric polynomials of even degree, we denote by  $\Pi_{f,g}$  the set of prime numbers $p$ such that  the
following conditions holds :

\medskip
The polynomial $f$ has a symmetric, irreducible factor $f' \in {\bf Z}_p[X]$, the polynomial $g$ has a symmetric, irreducible factor $g' \in {\bf Z}_p[X]$, such that
 $f'  \ {\rm (mod \ {\it p})}$ and $g'  \ {\rm (mod \ {\it p})}$
 have a common irreducible, 
symmetric factor in ${\bf F}_p[X]$.

\end{notation}

\medskip

For all prime numbers $p$, let $D^p_+, D^p_- \in {\bf Q}_p^{\times}/{\bf Q}_p^{\times 2}$.

\begin{notation} If  $f \in I_1$, let $\Pi_{f,X-1}$ be the set of prime numbers $p$
such that $p \in \Pi^r_f$,  that $f  \ {\rm (mod \ {\it p})}$
is divisible by $X-1$ 
 in ${\bf F}_p[X]$, and that if $n^+ = 2$, then $D^p_+ \not = -1$.
 
 \medskip
 Let $\Pi_{f,X+1}$ be the set of prime numbers $p$
such that $p \in \Pi^r_f$, and that $f  \ {\rm (mod \ {\it p})}$
is divisible by $X+1$ 
 in ${\bf F}_p[X]$, and that if $n^- = 2$, then $D^p_- \not = -1$.
 
 \medskip
 Let $\Pi_{X-1,X+1} = \{2\}$ if the following conditions hold : $n^+ \not = 0$, $n^- \not = 0$, and  if $n^+ = 2$, then  $D^2_+ \not = -1$; if $n^- = 2$, then $D^2_- \not = -1$.
 Otherwise, set $\Pi_{X-1,X+1} = \varnothing$.

\end{notation}

\medskip
We denote by $C(I)$ the set of maps $I \to {\bf Z}/2{\bf Z}$.

\begin{notation}\label{f,g}
 
 If $f, g \in I$, let $c_{f,g} \in C(I)$ be such that  

\medskip
\centerline {$c_{f,g}(f) = c_{f,g}(g) = 1$, and $c_{f,g}(h) = 0$
if $h \not = f,g$. }

\medskip
Let $(f,g): C(I) \to C(I)$
be the map  sending $c$ to $c + c_{f,g}$. 

\end{notation}

\begin{notation}

Let $C_0 (I)$ be the set of $c \in C(I)$ such that

\medskip 
\centerline { $c(f) = c(g)$ \  if  \ $\Pi_{f,g} \not =  \varnothing$,} 

\medskip 
\noindent
and we denote by $G_F(D_+,D_-)$  the quotient of the group
 $C_0(I)$ by the subgroup of constant maps. 
 
 \medskip 
 In general,  the group  depends on $D_+ = (D_+^p)$ and $D_- = (D_-^p)$. If $n^+ \not = 2$ and $n^- \not = 2$, then $G_F(D_+,D_-)$ only depends on $F$, and
 we denote it by $G_F$; similarly, we use the notation $G_F(D_-)$ if $n^+ \not = 2$,   and $G_F(D_+)$ if $n^- \not = 2$.

 \end{notation}
 
 The following proposition is useful in some of the examples.
 
 \begin{prop}\label{for examples} Let $f, g \in {\bf Z}[X]$ be monic, irreducible, symmetric polynomials of even degree, and $p$ be a prime number not dividing $(fg)(1)(fg)(-1)$. If
 $f  \ {\rm (mod \ {\it p})}$ and $g \ {\rm (mod \ {\it p})}$
 have a common irreducible, 
symmetric factor of even degree in ${\bf F}_p[X]$, then $p \in \Pi_{f,g}$.

 \end{prop}

The proposition is a consequence of the following lemma :

\begin{lemma}\label{lemma for examples} Let $p$ be a prime number, and let $f_0 \in {\bf Z}_p[X]$ be a monic, irreducible polynomial of degree $n$. Set $f(X) = X^{n}f_0(X + X^{-1})$, and
assume that $p$ does not divide $f(1)f(-1)$. Let $\overline f$ be the image of $f$ in ${\bf F}_p[X]$. If $\overline f$ has an irreducible, symmetric factor in ${\bf F}_p[X]$, then
$f  \in {\bf Z}_p[X]$ is irreducible.

\end{lemma}

\noindent
{\bf Proof.} Set $E = {\bf Q}_p[X]/(f)$, and let $\alpha$ be the image of $X$ in $E$; set $E_0 = {\bf Q}_p[X]/(f_0)$. If $f$ is irreducible, then $E/E_0$ is a quadratic
field extension, otherwise $E = {\bf Q}_p[X]/(h) \times {\bf Q}_p[X]/(h^*)$ for some irreducible polynomial $h \in {\bf Q}_p[X]$ such that $h \not = h^*$, and $f = h h^*$.
Set $\delta = (\alpha - \alpha^{-1})^2$. Since $p$ does not divide $f(1)f(-1)$ by hypothesis, the elements $\alpha + 1$ and $\alpha -1$ are units, and hence
$\delta$ is a unit of $E_0$. Let $O_0$ be the ring of integers of $E_0$, let $m_0$ be its maximal ideal, and set $\kappa_0 = O_0/m_0$. 
Suppose that $f = hh^*$ with $h \in {\bf Q}_p[X]$ such that $h \not = h^*$; then $\delta$ is a square in $O_0$, hence its image $\overline \delta$ is a square 
in $\kappa_0$. This implies that the irreducible factors of $\overline f$ in ${\bf F}_p[X]$ are not symmetric, and this concludes the proof of the lemma.

\medskip
\noindent
{\bf Proof of Proposition \ref{for examples}.} Let $h \in  {\bf F}_p[X]$ be a common irreducible, symmetric factor of $f  \ {\rm (mod \ {\it p})}$ and $g \ {\rm (mod \ {\it p})}$.
By Lemma \ref{lemma for examples} there exist symmetric, irreducible polynomials $f',g' \in {\bf Z}_p[X]$ such that $f'$ divides $f$, $g'$ divides $g$, and that
$h$ is a common factor of $f' \ {\rm (mod \ {\it p})}$ and $g' \ {\rm (mod \ {\it p})}$. This implies that $p \in \Pi_{f,g}$.



 



\medskip
In the case where $F$ has no linear factors, the group $G_F$ was defined in  \cite{B 21}, and
several examples are given in  \cite{B 21}, \S 25 and \S 31. Here are some more examples :

 \begin{example}\label{Sa}  Let $a$ be an integer $\geqslant 0$, set
$$S_a(X) = X^6 - aX^5 - X^4 + (2a-1) X^3 - X^2 - aX +1,$$ and $R_a(X)  = (X^2 - 4)(X-a) - 1$; we
have $S_a(X) = X^3 R_a(X+X^{-1})$. 

\medskip If $m$ is an integer $\geqslant 3$, recall that  $\Phi_m$ is the $m$-th cyclotomic polynomial, and 
that ${\rm Res}(S_a,\Phi_m)$ is  the resultant of $S_a$ and $\Phi_m$.

\medskip
(i) We have $\Pi_{S_a,\Phi_3} \not  = \varnothing$. Indeed, $${\rm Res}(S_a,\Phi_3) = {\rm N}_{{\bf Q}(\zeta_3)/{\bf Q}}(S(\zeta_3)) = {\rm N}_{{\bf Q}(\zeta_3)/{\bf Q}}(R(-1)) = 3(a+1) -1.$$ Since  $a \geqslant 0$, we have
$3(a+1) -1 \equiv 2 \ {\rm (mod \ 3)}$, therefore 
${\rm Res}(S_a,\Phi_3)$ is divisible by a prime number $p$ such that $p \equiv 2 \ {\rm (mod \ 3)}$.
The polynomial $\Phi_3$ is irreducible modulo $p$, hence $p \in \Pi_{S_a,\Phi_3}$.

\medskip Set $F_a(X) = S_a (X) \Phi_3^2(X) (X-1)^{12}$; note that with the above notation we have $n^+ = 12$ and $n^- = 0$, hence the obstruction group
only depends on $F_a$. We have $\Phi_3(1) = 3$ and $3 \in \Pi^r_{\Phi_3}$, hence $3 \in \Pi_{\Phi_3(X),X-1}$, and therefore $G_{F_a} = 0$.

\medskip
(ii) Set $P_a(X) = S_a(X) \Phi_4^2 (X) (X-1)^{12}$. We have $2 \in \Pi_{\Phi_4(X),X-1}$, and
$${\rm Res}(S_a,\Phi_4) =
{\rm N}_{{\bf Q}(i)/{\bf Q}}(R(0)) = 4a -1.$$
If $a \geqslant 1$, then $\Pi_{S_a,\Phi_4} \not  = \varnothing$, and $G_{P_a} = 0$. 
On the other hand, $\Pi_{S_0,\Phi_4} = \varnothing$, and $G_{P_a}$ is of order $2$.

\end{example}

\begin{example}\label{a=3} With the notation of Example \ref{Sa}, set $a = 3$, 
$$S(X) =  S_3(X) = X^6 - 3X^5 - X^4 + 5X^3 - X^2 - 3X + 1.$$ 

\medskip {\rm (i)} Set $F = S \Phi_{10}^4$. Then
$G_F$ is of order $2$. Indeed, the resultant of $S$ and $\Phi_{10}$ is $121$, and 
these two polynomials have the irreducible common factors $X+3$ and $X + 4$ in ${\bf F}_{11}[X]$. These polynomials
are not symmetric, therefore
$\Pi_{S,\Phi_{10}} = \varnothing$; hence $G_F \simeq {\bf Z}/2{\bf Z}$. 

\medskip {\rm (ii)} Set $F(X) = S(X)  \Phi_{10}(X)^2 (X-1)^8$; in this case, we have $G_F \simeq ({\bf Z}/2{\bf Z})^2$. 
Indeed, we have seen in (i) that $\Pi_{S,\Phi_{10}} = \varnothing$, and it is clear that $\Pi_{S(X),X-1} = \varnothing$, $\Pi_{\Phi_{10}(X),X-1} = \varnothing$. 

\end{example}

\begin{defn}\label{equivrel}

 If $p$ is a prime number, we  consider the equivalence relation on $C(I)$ generated by the elementary equivalence
$$a \sim b \iff  b = (f,g) \circ a \ {\rm with} \ p \in \Pi_{f,g},$$ and we denote by $ \sim_p$ the equivalence
relation on $C(I)$ generated by this elementary equivalence. 

\end{defn}

 \begin{lemma}\label{c and a} Let $c \in C_0(I)$, let $p$ be prime number, and let $a, b \in C(I)$
 such that $a \sim_p b $. Then $$ \underset{h \in I} \sum c(h) a(h) =  \underset{h \in I} \sum c(h) b(h).$$
 
 \end{lemma}

\noindent
{\bf Proof.} We can assume that $b = (f,g) \circ a$ with 
$p \in \Pi_{f,g}.$ By definition, we have $b(h) = a(h)$ if $h \not = f,g$, $b(f) = a(f) + 1$ and 
$b(g) = a(g) + 1$, hence $$\underset{h \in I} \sum c(h) b(h) =  \underset{h \in I} \sum c(h) a(h) + c(f) + c(g).$$
Since $c \in C_0(I)$ and $\Pi_{f,g} \not = \varnothing$ we have $c(f) = c(g)$, and this proves the claim. 

\bigskip

\section{Local data and equivalence relation}

We now combine the results and notions of \S \ref{local data}, \S \ref{local 0 section} and \S \ref{combinatorial section}.  We keep the notation
of these sections : $F \in {\bf Z}[X]$ is a monic, symmetric polynomial, $p$ is a prime number, $(V,q)$ is
a non-degenerate quadratic form over ${\bf Q}_p$,  and $t : V \to V$ is a semi-simple isometry of $q$ with characteristic
polynomial $F$. Recall that $F = F_0 F_1$, with 

\medskip
\centerline {$F_1 = \underset{f \in I_1} \prod f^{n_f}$ and
$F_0(X) = (X-1)^{n^+} (X+1)^{n^-}$} 
\noindent for some integers $n^+,n^- \geqslant 0$.

\medskip 

Recall that $D_0 = {\rm det}(q)F_1(1)F_1(-1)$, and let $D_+, D_- \in {\bf Q}_p^{\times}/{\bf Q}_p^{\times 2}$ be such that $D_0 = D_+ D_-$.

\bigskip We fix  $\delta \in W_{\Gamma}({\bf F}_p)$, and set  ${\mathcal C}^p = {\mathcal C}_{\delta}^p$; see Notation \ref {curly C bis} for the notation ${\mathcal C}^p $,
and Definition \ref{equivrel} for the definition of the equivalence relation $\sim_p$.

 \begin{prop}\label{equivalence class} The set  $\mathcal C^p$  is a $\sim_p$-equivalence class of $C(I)$. 
 
 \end{prop}

 \noindent
 {\bf Proof.} Set $A^p = w_2(q) + w_2(Q)$
in ${\rm Br}_2({\bf Q}_p) = {\bf Z}/2{\bf Z}$, and note that for all $a^p \in \mathcal C^p$, we have 
 $\underset{f \in J}\sum a^p(f) = A^p$.
 
 We start by proving that the set $\mathcal C^p$ is stabilized by the maps $(f,g)$ for $p \in \Pi_{f,g}$.
  Let $a^p[\lambda^p,q_{\pm}^p] \in \mathcal C^p$, let $f,g \in I$ be such that $p \in \Pi_{f,g}$, and let
 us show that $(f,g)(a^p[\lambda^p,q_{\pm}^p])  \in \mathcal C^p$. Note that if $f \in I_1$, then $p \in \Pi_{f,g}$ implies that  
 $(E^p_f)_0^{\times}/{\rm N}_{E^p_f/(E^p_f)_0}(E^p_f) \not = 0$. 
 Assume first that $f,g \in I_1$. There exist $\mu_f, \mu_g \in (E^p_f)_0^{\times}/{\rm N}_{E^p_f/(E^p_f)_0}(E^p_f)$
 such that ${\rm cor}_{(E_f)_0^p/{\bf Q}_p}(\mu_f,d_f) \not = {\rm cor}_{(E_f)_0^p/{\bf Q}_p}(\lambda_f,d_f)$
 and ${\rm cor}_{(E_f)_0^p/{\bf Q}_p}(\mu_g,d_g) \not = {\rm cor}_{(E_f)_0^p/{\bf Q}_p}(\lambda_g,d_g)$.
 Let $\mu^p \in E_0^p$ be obtained by replacing  $\lambda^p_f$ by $\mu^p_f$, $\lambda^p_f$ by $\mu^p_f$, and
 leaving the other components unchanged. 
 We have $a^p[\mu^p,q_{\pm}^p] = (f,g)(a^p[\lambda^p,q^p_{\pm}])$. Using the arguments of  \cite {B 21}, Propositions 16.5 and 22.1
 we see that  $a^p[\mu^p,q_{\pm}^p] \in \mathcal C^p$. Assume now that $f \in I_1$ and $g(X) = X -1$. In this case, 
 the hypothesis  $p \in \Pi_{f,g}$ implies that there exists a place $w$ of $(E_f)_0$ above $p$ that
 ramifies in $E_f$, and such that the $w$-component $\lambda^w$ of $\lambda^p$ is such that with the notation
 of \cite {B 21}, \S 22, $\partial_p(\lambda^w)$ is in $W_{\Gamma}({\bf F}_p,N_+)$. 
 We modify the $w$-component of $\lambda^p$ to obtain $\mu_f  \in (E^p_f)_0^{\times}/{\rm N}_{E^p_f/(E^p_f)_0}(E^p_f)$
 such that 
 ${\rm cor}_{(E_f)_0^p/{\bf Q}_p}(\mu_f,d_f) \not = {\rm cor}_{(E_f)_0^p/{\bf Q}_p}(\lambda_f,d_f)$, and let
 $b^p$ be a quadratic form over $\bf Q_p$ with ${\rm dim}(b^p) = {\dim}(q^p_+)$, ${\rm det}(b^p) = {\det}(q^p_+)$,  and
 $w_2(b^p) = w_2(q^p_+) + 1$. We have $a^p[\mu^p,b^p,q^p_-] = (f,g)(a^p[\lambda^p,q^p_{\pm}])$.
 The arguments of  \cite {B 21}, Propositions 16.5 and 22.1 show that $a^p[\mu^p,b^p,q^p_-]  \in \mathcal C^p$. 
 The proof is the same when $f \in I_1$ and $g(X) = X + 1$.
 
 \medskip 
 Conversely, let us show that if $a^p[\lambda^p,q_{\pm}^p]$ and $a^p[\mu^p,b_{\pm}^p]$ are in $\mathcal C^p$, then 
 $a^p[\lambda^p,q_{\pm}^p] \sim_p a^p[\mu^p,b_{\pm}^p]$.
 Let $I'$ be the set of $f \in I$ such that   $a^p[\lambda^p,q_{\pm}^p](f) \not = a^p[\mu^p,b_{\pm}^p](f)$. Since
 $\underset{h \in I}\sum a^p(h) = A^p$ for
 all $a^p \in \mathcal C^p$, the set $I'$ has an even number of elements. 
 
 \medskip
 Assume first that $p \not = 2$. This implies that for all $f \in I'$, we have $\partial_p(\lambda^p_f) \not = \partial_p(\mu^p_f)$
 and that if $f(X) = X \pm 1$, then $\partial_p(q_{\pm}^p) \not = \partial_p(b_{\pm}^p)$. Hence there exist $f, g \in I'$
 with $f \not = g$ such that $\partial_p(W_{\Gamma}({\bf Q}_p,M_f^p))$ and $\partial_p(W_{\Gamma}({\bf Q}_p,M_g^p))$ have a
  non-zero intersection. This implies that $p \in \Pi_{f,g}$. The element $(f,g)(a^p[\lambda^p,q^p_{\pm}])$ differs
  from $a^p[\mu^p,b_{\pm}^p]$ in less elements than 
$a^p[\lambda^p,q_{\pm}^p]$. Since $I'$ is a finite set, continuing this way we see that $a^p[\lambda^p,q_{\pm}^p] \sim_p a^p[\mu^p,b_{\pm}^p]$.

\medskip
Suppose now that $p = 2$. Let $I''$ be the set of $f \in I'$ such that $\partial_2(\lambda^2_f) \not = \partial_2(\mu^2_f)$,
and note that $I'$ has an even number of elements. The same argument as in the case $p \not = 2$ shows
that applying maps $(f,g)$, we can assume that $I'' = \varnothing$. If $f \in I'$ and $f \not \in I''$, then
$\partial_2(\lambda^2_f)$ belongs to $W_{\Gamma}({\bf F}_2,1) \subset W_{\Gamma}({\bf F}_2)$.  Therefore $f,g  \in I'$ and $f,g  \not \in I''$, then 
$2 \in \Pi_{f,g}$. The number of these elements is also even, hence after a finite number of elementary
equivalences we see that $a^p[\lambda^p,q_{\pm}^p] \sim_p a^p[\mu^p,b_{\pm}^p]$. This completes
the proof of the proposition.

\medskip
 \begin{notation}\label {epsilon p} Let $a^p \in \mathcal C^p$, and let $c \in C(I)$. Set $$\epsilon_{a^p}(c) = 
 \underset{f \in I} \sum c(f) a^p(f).$$
 
 \end{notation}

 \begin{lemma} Let $a^p$, $b^p$ be two elements of $\mathcal C^p$, and let $c \in C_0(I)$. Then
  $$\epsilon_{a^p}(c) = \epsilon_{b^p}(c).$$
 
 \end{lemma}

\noindent
{\bf Proof.} By Proposition \ref{equivalence class}, we have $a^p \sim_p b^p $; therefore Lemma
\ref{c and a}
implies the desired result.

\bigskip
Since $\epsilon_{a^p}(c)$ does not depend on the choice of  $a^p \in \mathcal C^p$, we set
$\epsilon^p (c) = \epsilon_{a^p}(c)$ for some $a^p \in \mathcal C^p$, and obtain a map $$\epsilon^p : C_0(I) \to {\bf Z}/2{\bf Z}.$$

\section{A Hasse principle}\label{Hasse section}

The aim of this section is to decide whether the ``building blocks" of \S \ref{local data} and \S \ref{local 0 section}
give rise to a global isometry. We do this in a more general context than needed in Part I; this is done to prepare
Part II, by including some results that will be useful there. 

\medskip
Let $\tau$ be a signature map with associated
polynomial $F$, and assume that $F \in {\bf Z}[X]$; we write $F = F_0 F_1$ as in the previous sections. Let $V$ be a finite dimensional ${\bf Q}$-vector
space,  let $q : V \times V \to {\bf Q}$ be a non-degenerate quadratic form, and let
$D_+, D_- \in {\bf Q}^{\times}/{\bf Q}^{\times 2}$ be such that $D_+ D_- = {\rm det}(q) F_1(1) F_1(-1)$. 
Let $\delta = (\delta_p)$ with $\delta_p \in W_{\Gamma}({\bf F}_p)$ be such that $\delta_p = 0$ for
almost all prime numbers $p$. 

\medskip
Suppose that for all prime numbers $p$, there exists a semi-simple isometry 
 $t : V \otimes_{\bf Q} {\bf Q}_p \to V \otimes_{\bf Q} {\bf Q}_p$ of
$(V,q) \otimes_{\bf Q} {\bf Q}_p$ with characteristic polynomial $F$ such that
the associated ${\bf Q}_p[\Gamma]$-bilinear form $(V,q) \otimes_{\bf Q} {\bf Q}_p$  satisfies
$$\partial_p [(V,q) \otimes_{\bf Q} {\bf Q}_p] = \delta_p$$ in $W_{\Gamma}({\bf F}_p)$. Set $V_{F_i}^p = {\rm Ker}(F_i(t))$; 
recall that we have 
orthogonal decompositions of ${\bf Q}_p[\Gamma]$-quadratic forms $$(V,q) \otimes_{\bf Q} {\bf Q}_p = (V_{F_0},q_0^p) \oplus (V_{F_1},q_1^p) 
\ {\rm and} \ (V_{F_0},q_0^p) \simeq (V_+,q_+^p) \oplus (V_-,q_-^p);$$
we assume that ${\rm det}(q_+^p) = D_+$ and ${\rm det}(q_+^p) = D_+$ in ${\bf Q}_p^{\times}/{\bf Q}_p^{\times 2}$. 

\medskip To this data, we associate $a^p \in \mathcal{C}^p = \mathcal{C}_{\delta}^p$  as in \S \ref{local 0 section}.

 \begin{prop}\label{almost zero} For almost all prime numbers $p$, the zero map belongs to the set $\mathcal C^p$. 
 
 \end{prop}

 \noindent
 {\bf Proof.} Let $S$ be the set of prime numbers such that $p$ is ramified in the extension $E_f/{\bf Q}$ for some
 $f \in I_1$, or $w_2(q) \not = w_2(Q)$ in ${\rm Br}_2({\bf Q}_p)$; this is a finite set. We claim that if $p \not \in S$, then 
 the zero map belongs to $\mathcal C^p$. Indeed, set $q_f^p = Q_f^p$ for all $f \in I_1$, and 
 $q_{\pm}^p = Q_{\pm}^p$. We have ${\rm det}(q) = {\rm det}(Q)$ in $\bf Q^{\times}/{\bf Q}^{\times 2}$, and if $p \not \in S$ we have $w_2(q)  = w_2(Q)$ in ${\rm Br}_2({\bf Q}_p)$, therefore, for $p \not \in S$, we have
 
 $$(V,q) \otimes_{\bf Q} {\bf Q}_p = \underset {f \in I_1} \oplus (V^p_f,q_f^p) \oplus (V^p_+,q_+^p) \oplus (V^p_-,q_-^p).$$
 
 If $p$ is unramified in all the extensions  $E_f/{\bf Q}$ for $f \in I_1$, by \cite{B 21}, Lemma 11.2 we
 have $\partial_p[V^p_f,q_f^p] = 0$ in $W_{\Gamma}({\bf F}_p)$; moreover,  $v_p(D_{\pm}) = 0$, hence 
 $\partial [V^p_{\pm},q_{\pm}^p] = 0$ in $W_{\Gamma}({\bf F}_p)$. 
 
 \medskip The above arguments show that if $p \not \in S$, then the choice of $q_f^p = Q_f^p$ for all $f \in I_1$ and 
 $q_{\pm}^p = Q_{\pm}^p$ gives rise to the element $a^p = 0$ of  $\mathcal C^p$; therefore the zero map is in $\mathcal C^p$,
 as claimed. This completes the proof of the proposition.
 
 \bigskip
 Let us denote by $\mathcal V$ the set of places of ${\bf Q}$, by $\mathcal V'$ the set of finite places,
 and by $\infty$ the unique infinite place.
We have an orthogonal decomposition  $$(V,q) \otimes_{\bf Q} {\bf R} = (V_{F_0},q_0^{\infty}) \oplus (V_{F_1},q_1^{\infty}) 
$$ 
of ${\bf R}[\Gamma]$-quadratic forms determined by the signature map $\tau$, with
$(V_{F_1}, q_1^{\infty}) = \underset {f \in I_1} \oplus (V^{\infty}_f,q_f^{\infty})$. If moreover
$n_+ > 0$ and $n_- > 0$, 
we then have the further orthogonal decomposition of ${\bf R}[\Gamma]$-quadratic forms  $(V_{F_0},q_0^{\infty}) \simeq (V_+,q_+^{\infty}) \oplus (V_-,q_-^{\infty})$. 
The  ${\bf R}[\Gamma]$-quadratic forms $(V^{\infty}_f,q_f^{\infty})$ and $(V^{\infty}_{\pm},q_{\pm}^{\infty})$ are 
also determined by the signature map $\tau$. 

\begin{prop}\label{real Hasse-Witt} We have ${\rm dim}(q_f^{\infty}) = {\rm deg}(f)n_f$, ${\rm det}(q_f^{\infty}) = [f(1)f(-1)]^{n_f}$, and
the Hasse-Witt invariant of $q_f^{\infty}$ satisfies
$$w_2(q_f^{\infty}) + w_2(Q_f) = {\rm cor}_{(E_f)_0^{\infty}/{\bf R}}({\rm det}(h^{\infty}_f),d_f)$$ in ${\rm Br}_2({\bf R})$. 

\end{prop}

\noindent
{\bf Proof.} The assertion concerning the dimension is clear, the one on the determinant follows from 
\cite{B 15}, Corollary 5.2,  and the property of the Hasse-Witt invariants from \cite {B 21}, Proposition 12.8.

\medskip
 For $f \in I_1$, set
$E_f^{\infty} = E_f \otimes_{\bf Q} {\bf R}$ and $(E_f)_0^{\infty} = (E_f)_0 \otimes_{\bf Q} {\bf R}$. There exists a unique non-degenerate hermitian form $(M_f^{\infty},h_f)$ over
$(E_f^{\infty},\sigma_f)$ such that $$q_f^{\infty}(x,y) = {\rm Tr}_{E_f^{\infty}/{\bf R}}(h^{\infty}_f(x,y)),$$
see for instance \cite {M}, Lemma 1.1 or \cite {B 15}, Proposition 3.6.
Set 
$$\lambda^{\infty}_f = {\rm det}(h^{\infty}_f) \in (E^{\infty}_f)_0^{\times}/{\rm N}_{E^{\infty}_f/(E^{\infty}_f)_0}.$$

Set $a^{\infty}_{\tau}(f) = w_2(q_f^{\infty}) + w_2(Q_f)$ if $f \in I_1$, and note that by Proposition \ref{real Hasse-Witt}
we have $a_{\tau}^{\infty}(f) =  {\rm cor}_{(E_f)_0^{\infty}/{\bf R}}(\lambda^{\infty}_f,d_f)$.

\medskip
 Set $$a_{\tau}^{\infty}(X \pm 1) = w_2(q^{\infty}_{\pm}) + w_2(Q_{\pm}).$$ 
 

\medskip
We record the following observation for later use : 

\begin{prop}\label{dependance on tau} If $\tau$ and $\tau'$ are two signature maps with the same
maximum and associated polynomial, and if $f \in I$ is such that 
${\rm proj}_2(\tau(f)) \equiv  {\rm proj}_2(\tau'(f))  \ {\rm (mod \ 4)}$, then $a^{\infty}_{\tau}(f) =
a^{\infty}_{\tau'}(f)$.

\end{prop}

\noindent
{\bf Proof.} Let $f \in I$, let $\tau(f) = (r_f,s_f)$. Note that ${\rm proj}_2(\tau(f)) = s_f$, and
that $q_f^{\infty}$ is the diagonal quadratic form over $\bf R$ with $r_f$ entries $1$ and $s_f$
entries $-1$. The Hasse-Witt invariant $w_2(q_f^{\infty})$ is determined by the class of $s_f$
$ {\rm (mod \ 4)}$, and this implies the proposition. 

\bigskip
We obtain a map

$$\epsilon^{\infty}_{\tau} :  C(I) \to {\bf Z}/2{\bf Z}$$ by setting $$\epsilon^{\infty}_{\tau}(c) =  \underset{f \in I} \sum c(f) a_{\tau}^{\infty}(f).$$

For $v \in \mathcal V$, set $\epsilon^v = \epsilon^p$ if $v = v_p$, and $\epsilon^v = \epsilon^{\infty}_{\tau}$ if $v = v_{\infty}$. Set

$$\epsilon_{\tau}(c) = \underset{v \in \mathcal V} \sum \epsilon^v(c).$$

Since $\epsilon^v = 0$ for almost all $v \in \mathcal V$ (cf. Proposition \ref{almost zero}), this is a finite sum. 

\bigskip Set $\epsilon_{\tau}^{\rm finite} = \underset {v \in {\mathcal V'}} \sum \epsilon^v$, and note that 
if $n^+ = 0$ or $n^- = 0$, then $\epsilon_{\tau}^{\rm finite}$ does not depend on $D_+$ and $D_-$. 
By definition, we have
$\epsilon_{\tau} = \epsilon_{\tau}^{\rm finite} + \epsilon^{\infty}_{\tau}$. We obtain
a homomorphism

$$\epsilon _{\tau}: C_0(I) \to {\bf Z}/2{\bf Z}.$$

\begin{prop}\label{epsilon} The homomorphism $\epsilon_{\tau} : C_0(I) \to {\bf Z}/2{\bf Z}$ induces a homomorphism
 $$\epsilon_{\tau} : G_F(D_+,D_-)  \to {\bf Z}/2{\bf Z}.$$

\end{prop}

\noindent
{\bf Proof.} It suffices to show that  if $c(f) = 1$ for all $f \in J$, then $\epsilon_{\tau}(c) = 0$. 
For all $v \in \mathcal V$, set $A^v = w_2(q) + w_2(Q)$
in ${\rm Br}_2({\bf Q}_v) = {\bf Z}/2{\bf Z}$, where 
${\bf Q}_v$ is either $\bf R$ or ${\bf Q}_p$, for a prime number $p$. Note that 
$A^v = 0$ for almost all $v \in \mathcal V$, and that
$\underset {v \in \mathcal V} \sum A^v = 0$. Moreover, for all $a^v \in \mathcal C^v$, we have by 
definition $\underset{f \in J}\sum a^v(f) = A^v$.

\medskip Let $c \in C(I)$ be such that $c(f) = 1$ for all $f \in J$. We have 
$$\epsilon_{\tau}(c) = \underset{v \in \mathcal V} \sum \  \underset {f \in J} \sum c(f) \ a^v(f) = \underset{v \in \mathcal V} \sum \  \underset {f \in J} \sum  \ a^v(f) = \underset {v \in \mathcal V} \sum A^v = 0.$$

\begin{notation}
We denote by $\mathcal C$ the set of $(a^v)$ with $a^v \in \mathcal C^v$ and $a^v = 0$ for almost all $v$.

\end{notation}

By Proposition \ref{almost zero}, the set  $\mathcal C$  is not empty.

\begin{theo}\label{realization}  Let $(a^v) \in \mathcal C$ be such that $\underset{v \in \mathcal V} \sum a^v(f) = 0$ for
all $f \in I$. Then there exists a semi-simple isometry $t : V \to V$ of $(V,q)$ with signature map $\tau$
and local data $(a^v)$.

\end{theo}

\noindent
{\bf Proof.} 
 \bigskip
 If $f \in I_1$, we have 
 $a^{v_p}(f) = a[\lambda^p,q_+^p,q_-^p](f) = {\rm cor}_{(E_f)_0^p/{\bf Q}_p}(\lambda^p_f,d_f),$ and 
 $a^{v_{\infty}}(f) = a[\lambda^{\infty},q_+^{\infty},q_-^{\infty}](f) =  {\rm cor}_{(E_f)_0^{\infty}/{\bf R}}(\lambda^{\infty}_f),d_f)$ (cf. Propositions \ref{dim det w f} and  \ref{real Hasse-Witt}).
  Since $\underset {v \in \mathcal V} \sum a^v(f) = 0$, this implies that
  $$\underset {v \in \mathcal V} \sum {\rm cor}_{(E_f)_0^v/{\bf Q}_v}(\lambda^v_f,d_f) = 0,$$
  where ${\bf Q}_v = {\bf Q}_p$ if $v = v_p$ and ${\bf Q}_v = {\bf R}$ if $v = v_{\infty}$. Therefore we have
   $$\underset {w \in {\mathcal W}} \sum (\lambda^w_f,d_f) = 0,$$ where $\mathcal W$ is the set of primes of $(E_f)_0$. 
   This implies that there exists $\lambda_f \in (E_f)_0^{\times}/N_{E_f/(E_f)_0}(E_f^{\times})$ mapping to $\lambda_f^w$ for all
   $w \in {\mathcal W}$ (see for instance \cite{B 21}, Theorem 10.1). In particular, we have $(\lambda_f,d_f) = (\lambda^w_f,d_f)$
   in ${\rm Br}_2((E_f)_0^w)$ for all $w \in {\mathcal W}$. 
   
   \medskip
   Note that ${\rm diff} \circ \tau(f)$ is an even integer. Let $h_f : V_f \times V_f \to E_f$ be a hermitian form of 
determinant $ \lambda_f$ and index  ${1 \over 2 }{\rm diff} \circ \tau(f)$; such a hermitian
   form exists (see for instance  \cite{Sch}, 10.6.9). Let us define   
   $q_f : V_f \times V_f \to {\bf Q}$ by $q_f(x,y) = {\rm Tr}_{E_f/{\bf Q}}(h_f(x,y)).$

   \bigskip
   Let $f = X \pm 1$. We have  $\underset {v \in \mathcal V} \sum a^v(f) = 0$,  hence by the Brauer-Hasse-Noether theorem
   there exists $a(\pm) \in {\rm Br}_2({\bf Q})$ mapping to $a^v(f)$ in ${\rm Br}_2({\bf Q}_v)$ for all $v \in \mathcal V$.
   Let $q_{\pm}$ be a quadratic form over $\bf Q$ of dimension $n_{\pm}$, determinant $D_{\pm}$, Hasse-Witt
   invariant $w_2(q_{\pm}) = a(\pm) + w_2(Q_{\pm})$ and signature $\tau(X\pm 1) = (r_{\pm},s_{\pm})$. 
   Such a quadratic form exists; see for instance \cite{S}, Proposition 7.

\medskip
   Let $q' : V \times V \to {\bf Q}$ be the quadratic form given by $$(V,q')= \underset{f \in I_1} \oplus (V_f,q_f) \oplus 
   \underset {f \in I_0}\oplus (V_f,q_f).$$ By construction,
   $(V,q')$ has the same dimension, determinant, Hasse-Witt invariant and signature as $(V,q)$, hence
   the quadratic forms $(V,q')$ and $(V,q)$ are isomorphic. 
   
   \medskip Let $t : V \to V$ be defined by $t(x) = \gamma x$, where $\gamma$ is a generator of $\Gamma$. 
  By construction, $t$ is an isometry of $(V,q')$ and it is semi-simple with signature map $\tau$ and
  local data $(a^v)$.

\begin{theo}\label{tau} Suppose that $\epsilon_{\tau} = 0$. Then there exists $(a^v) \in \mathcal C$ such
that $\underset{v \in \mathcal V} \sum a^v(f) = 0$ for all $f \in I$.

\end{theo} 

\noindent
{\bf Proof.} This follows from \cite{B 21}, Theorem 13.5. 

\begin{coro}\label{HP} The quadratic form $(V,q)$ has a semi-simple isometry with signature map $\tau$ such that
the associated ${\bf Q}\Gamma$ quadratic form $(V,q)$ satisfies $\partial_p[V,q] = \delta_p$ for all prime numbers $p$ if and only if 
$\epsilon_{\tau} = 0$.

\end{coro} 

\noindent {\bf Proof.} Assume that $(V,q)$ has a semi-simple isometry as above, and let $(a^v)$ be
the associated local date. Recall that the
element  $a^p \in \mathcal C^p$ given by $a^p(f) = w_2(q^p_f) + w_2(Q_f)$ if $f \in I_1$, by $a^p(X \pm 1) = w_2(q^p_{\pm}) + w_2(Q_{\pm}),$ and  $a^p(f) = 0$ if $f \in J$ with $f \not \in I$, $f \not = X \pm 1$ (see \S \ref{local data}). 
Similarly, we have the element $a^{\infty} \in C(I)$ given by
 $a^{\infty}(f) = w_2(q^{\infty}_f) + w_2(Q_f)$ if $f \in I_1$,  $a^{\infty}(X \pm 1) = w_2(q^{\infty}_{\pm}) + w_2(Q_{\pm}),$ and 
$a^{\infty}(f) = 0$ if $f \in J$ with $f \not \in I$, $f \not = X \pm 1$. 
Since $q^p_f = q_f \otimes_{\bf Q}{\bf Q}_p$ and $q^{\infty}_f = q_f \otimes_{\bf Q}{\bf R}$, we have

\medskip

\centerline { $ \underset{v \in \mathcal V} \sum \  a^v(f) = 0$ for all $f \in J$.}

\medskip \noindent This implies that $\epsilon_{\tau} = 0$. The converse follows from Theorems \ref {realization}
and \ref{tau}.

\section{Necessary conditions}\label{necessary section}

Our aim is to determine the signature maps of semi-simple
isometries with determinant $1$ of  even unimodular lattices; in this section, we collect
some necessary conditions for this. 



\medskip We start by recalling a well-known result :

\begin{lemma} \label{mod 8} Let $(L,q)$ be an even unimodular lattice. Then

\medskip

{\rm (a)} For all prime numbers $p$, we have $w_2(q) = w_2(h)$ in ${\rm Br}_2({\bf Q}_p)$, where $h$ is the hyperbolic form over ${\bf Q}_p$
of dimension equal to the rank of $L$. 

\medskip {\rm (b)} Let $(r,s)$ be the signature of $(L,q)$. Then $r  \equiv s \ {\rm (mod \ 8)}$.

\end{lemma}

\noindent
{\bf Proof.} (a) The reduction mod 2 of $(L,q)$ is a non-degenerate alternating symmetric bilinear form over ${\bf F}_2$, hence
it is even dimensional; its dimension is equal to the rank of $L$, therefore this rank is even. Set
${\rm rank}(L) = 2n$. The determinant of $q$ is $\pm 1$; by \cite{bg}, Proposition 5.2  this implies
that 
$q \otimes_{\bf Z} {\bf Z}_2$ is isomorphic to $h \otimes_{\bf Z} {\bf Z}_2$, and hence
 $w_2(q) = w_2(h)$ in ${\rm Br}_2({\bf Q}_2)$. If $p$ is a prime number with $p \not = 2$, then
 $w_2(q) = w_2(h) = 0$; this proves (a). 
 
 \medskip (b) This is proved  for instance  in \cite{S}, Chap. V, Corollaire 1. Here is another proof : 
 since $q$ is a global form and $w_2(q) = 0$ in ${\rm Br}_2({\bf Q}_p)$ for $p \not = 2$, we have
 $w_2(q) = 0 $ in ${\rm Br}_2({\bf Q}_2)$ if and only if $w_2(q) = 0 $ in ${\rm Br}_2({\bf R})$. Moreover,
 $w_2(q) = 0$ in ${\rm Br}_2({\bf R})$ if and only if $s \equiv  \ 0  \ {\rm or} \ 1\ {\rm (mod \ 4)}$.
 On the other hand, by (a) we know that $w_2(q) = w_2(h)$ in ${\rm Br}_2({\bf Q}_2)$, therefore
 $w_2(q) = 0$ in ${\rm Br}_2({\bf Q}_2)$ if and only if $n \equiv  \ 0  \ {\rm or} \ 1\ {\rm (mod \ 4)}$.
 This implies that $n  \equiv s \ {\rm (mod \ 4)}$. We have $2n = r + s$, therefore 
 $r - s = 2(n-s)$, hence $r  \equiv s \ {\rm (mod \ 8)}$, as claimed. 

\bigskip

This motivates the following definition. Let  $F \in {\bf R}[X]$ and let $r,s \geqslant 0$ be two integers.

\begin{defn}\label{C 0} We say that condition (C 0) holds for  
$F$ and $(r,s)$ if

\medskip
(C 0) {\it $F \in {\bf Z}[X]$ with $F(0) = 1$, ${\rm deg}(F) = r + s$, and $r  \equiv s \ {\rm (mod \ 8)}$}.

\bigskip
\noindent
We say that condition (C 0) holds for a signature map $\tau$ if it holds for the maximum $(r,s)$ and the associated polynomial $F$
of  $\tau$.

\end{defn}

\begin{lemma}\label{conditionc0} Suppose that $\tau$ is the signature map of a semi-simple isometry
with determinant 1 of an even unimodular lattice. Then condition {\rm (C 0)} holds.

\end{lemma}

\noindent
{\bf Proof.}
It is clear that $F \in {\bf Z}[X]$ and that $F(0) = 1$, we have 
${\rm deg}(F) = r + s$ by definition, and $r  \equiv s \ {\rm (mod \ 8)}$ follows from Lemma \ref{mod 8}.

\bigskip

Recall from 
\S \ref{local} that condition (C 1) is said to hold for a polynomial $F \in {\bf Z}[X]$ with  $2n = {\rm deg}(F)$ if

\medskip
(C 1) {\it $|F(1)|$, $|F(-1)|$ and $(-1)^n F(1) F(-1)$ are squares.}

\begin{defn} 
We say that condition (C 1) holds for a signature map $\tau$ if it holds for the polynomial $F \in {\bf Z}[X]$ associated
to $\tau$. 

\end{defn}

\begin{lemma}\label{conditionc0} Suppose that $\tau$ is the signature map of a semi-simple isometry
with determinant 1 of an even unimodular lattice. Then condition {\rm (C 1)} holds.

\end{lemma}

\noindent
{\bf Proof.} 
This follows from 
Theorem \ref{local theorem} (see also \cite{GM}, Theorem 6.1). 

\begin{notation}\label {q r,s} If $r,s \geqslant 0$ are two integers, we denote by $(V_{r,s},q_{r,s})$
the diagonal quadratic form over ${\bf Q}$ with $r$ diagonal entries $1$ and $s$ diagonal entries $-1$.

\end{notation}

\begin{lemma}\label{next} Suppose that $r  \equiv s \ {\rm (mod \ 8)}$. Then for all prime numbers $p$, we have 
$w_2(q_{r,s}) = w_2(h)$ in ${\rm Br}_2({\bf Q}_p)$, where $h$ is the hyperbolic form over ${\bf Q}_p$
of dimension $r + s$. 

\end{lemma} 

\noindent
{\bf Proof.} This follows from an easy computation. 

\begin{lemma} Suppose that $\tau$ is the signature map of a semi-simple isometry
with determinant 1 of an even unimodular lattice $(L,q)$, and let $(r,s)$ be the maximum of $\tau$. 
 Then $(L,q)\otimes_{\bf Z}{\bf Q} \simeq (V_{r,s},q_{r,s})$. 

\end{lemma}

\noindent
{\bf Proof.} It is clear that $q$ and $q_{r,s}$ have the same dimension, determinant and signature. 
They have the same Hasse-Witt invariant by Lemma \ref{mod 8} (a) and Lemma \ref{next}, hence
they are isomorphic. 

\bigskip Let $\tau$ be a signature map, let $F \in {\bf Z}[X]$ be the associated polynomial, and
assume that $F = F_0F_1$ as in \S \ref{local data}, with 
$$F_0(X) = (X-1)^{n^+} (X+1)^{n^-}$$ for some integers $n^+,n^- \geqslant 0$. Let $(V,q) = (V_{r,s},q_{r,s})$.

\begin{notation}\label{decomposition}
Let $p$ be a prime number, and let $t : V \otimes_{\bf Q}{\bf Q}_p \to V \otimes_{\bf Q}{\bf Q}_p$ be
a semi-simple isometry with characteristic polynomial $F$. Let
$$(V,q) \otimes_{\bf Q} {\bf Q}_p = (V_{F_0},q_0^p) \oplus (V_{F_1},q_1^p) 
$$ be 
the orthogonal decomposition of ${\bf Q}_p[\Gamma]$-quadratic forms of \S \ref{local data}. If moreover
$n_+ > 0$ and $n_- > 0$, 
we then have the further orthogonal decomposition of ${\bf Q}_p[\Gamma]$-quadratic forms  $(V_{F_0},q_0^p) \simeq (V_+,q_+^p) \oplus (V_-,q_-^p)$, as in \S \ref{local data}.

\end{notation} 

\begin{prop}\label{tau and det} Assume that $n^+ > 0$ and $n^- > 0$, and let
$p$ be a prime number. With the above notation, we have

\medskip

\centerline {${\rm det}(q^p_+) = u_+ F_1(1)$, ${\rm det}(q^p_-) = u_- F_1(-1)$ in ${\bf Q}_p^{\times}/{\bf Q}_p^{\times 2}$}

\medskip
\noindent for some $u_+, u_- \in {\bf Q}_p^{\times}$ such that $u_+ u_- = (-1)^n$.

\end{prop}

\noindent
{\bf Proof.} This follows from \cite{BT}, Lemma 6.8 and proposition 6.6 if
$p \not = 2$, from \cite {BT}, Lemma 6.8, Proposition 6.7 and the fact that $v_2({\rm det}(q_-))  \equiv v_2(F_1(-1)) \ {\rm (mod \ 2)}$ since the lattice is even if $p = 2$.

\bigskip
We also need the following observation : 

\begin{lemma}\label{R} Let $\tau(X-1) = (r_+,s_+)$ and
$\tau(X+1) = (r_-,s_-)$. Let 
$$(V,q) = (V_{F_0},q_0^{\infty}) \oplus (V_{F_1},q_1^{\infty}) 
$$ be 
the orthogonal decomposition of ${\bf R}[\Gamma]$-quadratic forms of corresponding to
the signature map $\tau$.  If moreover
$n^+ > 0$ and $n^- > 0$, 
we then have the further orthogonal decomposition of ${\bf R}[\Gamma]$-quadratic forms  $(V_{F_0},q_0^{\infty}) \simeq (V_+,q_+^{\infty}) \oplus (V_-,q_-^{\infty})$.
Then we have

\medskip
\centerline {${\rm det}(q^{\infty}_+) = (-1)^{s_+}$, ${\rm det}(q^{\infty}_-) = (-1)^{s_-}$ in ${\bf R}^{\times}/{\bf R}^{\times 2}$.}

\end{lemma}

\noindent
{\bf Proof.} This is clear. 

\section{Local data - isometries of even unimodular lattices}\label{local uni}

Let $\tau$ be a signature map with maximum $(r,s)$ and associated polynomial $F$. We assume
that conditions (C 0) and (C 1) hold for $\tau$; let $n \geqslant 1$ be an integer such that
${\rm deg}(F) = 2n$.  Suppose that $F = F_0 F_1$, with the notation of  the previous section. 

\medskip
In the next
section, we give a necessary and sufficient condition for $\tau$ to be the signature map of
a semi-simple isometry with determinant $1$ of an even unimodular lattice. In order to apply the results
of  \S \ref{Hasse section} we need to fix a non-degenerate quadratic form $q : V \times V \to {\bf Q}$,
and $D_+, D_- \in {\bf Q}^{\times}/{\bf Q}^{\times 2}$ such that $D_+ D_- = {\rm det}(q) F_1(1) F_1(-1)$.

\medskip 

The results of the previous section imply that the only possible choice is $(V,q) = (V_{r,s},q_{r,s})$ 
(see Lemma \ref{next}). We have ${\rm det}(q) = (-1)^n$, hence the condition for $D_+$ and $D_-$
becomes $D_+ D_- = (-1)^n F_1(1) F_1(-1)$. 

\medskip We now determine the values of $D_+$ and $D_-$, also using the results of the previous
section. Let $\tau(X-1) = (r_+,s_+)$ and
$\tau(X+1) = (r_-,s_-)$. 
With the notation of Lemma \ref{R}, we have ${\rm det}(q) = (-1)^n = F_1(1) F_1(-1)  (-1)^{s_+}  (-1)^{s_-}$
in ${\bf R}^{\times}/{\bf R}^{\times 2}$. Set $\epsilon_+ = {F_1(1) \over {|F_1(1)|}}$ and 
$\epsilon_- = {F_1(-1) \over {|F_1(-1)|}}$; then we have $$(-1)^n = \epsilon_+ (-1)^{s_+} \epsilon_- (-1)^{s_-}.$$

 \medskip Set $u_+ =  \epsilon_+ (-1)^{s_+} $ and $u_- = \epsilon_- (-1)^{s_-}$. Then Lemma \ref{R}
 and Lemma \ref{tau and det} show that the only possible choice for $D_+$ and $D_-$ is
 $$D_+ = u_+ F_1(1) = (-1)^{s_+} |F_1(1)|  \ ,  \ D_- = u_- F_1(-1) = (-1)^{s_-} |F_1(-1)|.$$ 
 
 Note that $D_+$ and $D_-$ are determined by the signature map $\tau$. 
 We now show that this choice gives rise to non-empty
local data. 

\begin{prop}\label{necessary proposition} Let $u_+, u_- \in \{\pm 1 \}$ be such that $u_+ u_- = (-1)^n$,
and let $(V,q) = (V_{r,s},q_{r,s})$.
Then for every
prime number $p$, the quadratic form $(V,q)\otimes_{\bf Q} {\bf Q}_p$ has a semi-simple isometry 
with characteristic polynomial $F$
that stabilizes an even unimodular ${\bf Z}_p$-lattice such that if $n_+ > 0$ and $n_- > 0$, then
${\rm det}(q^p_+) = u_+ F_1(1)$, ${\rm det}(q^p_-) = u_- F_1(-1)$ in ${\bf Q}_p^{\times}/{\bf Q}_p^{\times 2}$.

\end{prop}

\noindent
{\bf Proof.} 
Condition (C 1) implies that for every prime number $p$, there exists an even ${\bf Z}_p$-lattice
of determinant $(-1)^n$ having a semi-simple isometry with characteristic polynomial $F$ (cf. 
Corollary \ref{local coro}); the statement concerning ${\rm det}(q^p_+)$ and ${\rm det}(q^p_+)$ follows
from Claim 1 and Claim 2 (see \S \ref{local}, proof of Theorem \ref{local theorem}).  It remains to
show that  for every prime number $p$, the extension to ${\bf Q}_p$ of such a  lattice is isomorphic to $(V,q) \otimes_{\bf Q} {\bf Q}_p$. Note that the conditions $2n =  r + s$ and $r  \equiv s \ {\rm (mod \ 8)}$
imply that $n  \equiv s \ {\rm (mod \ 4)}$, therefore $(-1)^s = (-1)^n$; this implies that ${\rm det}(q) = (-1)^n$.
Moreover, for all prime numbers $p$ we have  $w_2(q) = w_2(h)$ in ${\rm Br}_p({\bf Q}_p)$, where
$h$ is the hyperbolic lattice over ${\bf Q}_p$
of rank $2n$. This coincides with the Hasse-Witt invariants of an even unimodular 
${\bf Z}_p$-lattice of rank $2n$ (this is clear  for $p \not = 2$, and follows from \cite {bg}, Proposition 5.2 for $p = 2$).
Hence  the extension to ${\bf Q}_p$ of an even ${\bf Z}_p$-lattice
of rank $2n$ and determinant $(-1)^n$ is isomorphic to $(V,q) \otimes_{\bf Q} {\bf Q}_p$.

\section{A necessary and sufficient condition}\label{sufficient section}

Let $\tau$ be a signature map. The aim of this section is to apply the results of \S \ref{Hasse section} to 
decide whether $\tau$ is the signature map of a semi-simple isometry of an even unimodular lattice.
Let $F$ be the polynomial associated to $\tau$, let $(r,s)$ be the maximum of $\tau$, 
and set $(V,q) = (V_{r,s},q_{r,s})$. 

\medskip

Recall that condition (C 0) holds for  
$F$ and $(r,s)$ if

\medskip
(C 0) {\it $F \in {\bf Z}[X]$ with $F(0) = 1$, ${\rm deg}(F) = r + s$, and $r  \equiv s \ {\rm (mod \ 8)}$}.

\medskip
\noindent
and that condition (C 0) holds for  $\tau$ if it holds for  $(r,s)$ and $F$.

\medskip

Recall that condition (C 1) is said to hold for a polynomial $F \in {\bf Z}[X]$ with  $2n = {\rm deg}(F)$ if

\medskip
(C 1) {\it $|F(1)|$, $|F(-1)|$ and $(-1)^n F(1) F(-1)$ are squares.}

\medskip
\noindent
and that condition (C 1) holds for  $\tau$ if it holds for $F$.

\medskip
In the following, we {\it assume that conditions} (C 0) and (C 1)  {\it hold for $\tau$}.

\medskip
Let us write $F = F_0 F_1$ as
in the previous sections. With the notation of \S \ref{local uni}, set $$D_+ = u_+ F_1(1) = (-1)^{s_+} |F_1(1)| \ , \ D_- = u_- F_1(-1) = (-1)^{s_-} |F_1(-1)| .$$

Note that the group $G_F(D_+,D_-) = G_F ((-1)^{s_+}|F_1(1)|,(-1)^{s_-}|F_1(-1)|)$ only depends on the polynomial $F$ and the integers $s_+, s_1$.

\medskip

Let $\delta_p = 0$ for all prime numbers $p$.

\medskip For every
prime number $p$, the quadratic form $(V,q)\otimes_{\bf Q} {\bf Q}_p$ has a semi-simple isometry 
with characteristic polynomial $F$
that stabilizes an even unimodular ${\bf Z}_p$ lattice (cf. Proposition \ref{necessary proposition}). To
this isometry, we associate a homomorphism $\epsilon^p : C_0(I) \to {\bf Z}/2{\bf Z}$ as in \S \ref{local 0 section}
with $\epsilon^p = 0$ for almost all $p$, and we obtain a homomorphism (cf. \S \ref{Hasse section})
$$\epsilon_{\tau}^{\rm finite} : C_0(I) \to {\bf Z}/2{\bf Z}.$$

Note that if $F(1) \not = 0$ or $F(-1) \not =0$, then $\epsilon_{\tau}^{\rm finite}$ only depends only on the polynomial $F$, and in this case,  is denoted by
$$\epsilon_F^{\rm finite} : C_0(I) \to {\bf Z}/2{\bf Z}.$$

We associate to $\tau$ a homomorphism $\epsilon_{\tau}^{\infty} : C_0(I) \to {\bf Z}/2{\bf Z}$,
and set $\epsilon_{\tau} = \epsilon_{\tau} ^{\rm finite} + \epsilon_{\tau}^{\infty}$, obtaining a homomorphism 

$$\epsilon_{\tau} : G_{F}(D_+,D_-)  \to {\bf Z}/2{\bf Z}$$ (see Proposition \ref{epsilon}). 

\begin{theo}\label {theo special} There exists an even unimodular lattice having a semi-simple isometry with signature map $\tau$ if
and only if $\epsilon_{\tau} = 0$. 

\end{theo}

\noindent
{\bf Proof.} Assume that $\epsilon_{\tau} = 0$. By Corollary \ref{HP}, the quadratic form $(V,q)$ has a
semi-simple  isometry with signature map $\tau$ such that for all prime numbers $p$, the
quadratic form  $(V,q)\otimes_{\bf Q} {\bf Q}_p$ contains a unimodular ${\bf Z}_p$-lattice $L_p$  
that is stabilized by $t$; for $p = 2$, the ${\bf Z}_2$-lattice $L_2$ can be chosen to be even (this follows from \cite{BT}, Theorem 8.1, in which conditions (i) and (ii) are clearly satisfied,  and condition (iii) follows from  
\cite{BT}, Theorem 8.4 and Claim 2). We can also suppose that for almost all 
prime numbers $p$, the lattice $L_p$ is equal to $L' \otimes_{\bf Z}{\bf Z}_p$ for some ${\bf Z}$-lattice $L'$
of $(V,q)$ stabilized by $t$. Set $L = \underset{p} \cap L_p$, where the intersection is taken over all
prime numbers $p$; the lattice $L$ is even unimodular and stabilized by $t$. 
The converse is an immediate consequence of Corollary \ref{HP}. 

\bigskip
Note that if $n^+ \not = 2$ and $n^- \not = 2$, then  the group $G_{\tau}$ only depends on the polynomial $F$; in this case, we
denote it by $G_F$. Similarly, if $n^+ \not = 2$, we use the notation $G_F(D_-)$, and if $n^- \not = 2$, the notation $G_F(D_+)$.

\bigskip
We now derive some consequences of Theorem \ref{theo special} : in particular, we prove the results announced in the introduction to Part I.

\medskip
Let $t \in {\rm SO}_{r,s}({\bf R})$ be an semi-simple isometry with characteristic polynomial $F$, and suppose that $n^+ \not = 2$ and $n^- \not = 2$. 

\medskip

The following result is  Theorem 6 of the introduction :

\begin{theo}\label{Theo 8}
The isometry $t$ preserves an even unimodular lattice 
if and only if $\epsilon_{{\rm sign}_t} = 0$. 
\end{theo}

Theorem \ref{theo special}  also implies the following corollaries : 

\begin{coro}\label{G=0} If $G_F = 0$, then there exists an even unimodular lattice having an
isometry with signature map $\tau$. 

\end{coro}

\begin{coro}\label{all signatures} If $G_F = 0$, then there exists an even unimodular lattice of signature $(r,s)$  having an
isometry with characteristic polynomial $F$.

\end{coro}

If $S$ is a Salem polynomial of degree $22$, then Condition (C 1) holds for $S$ if and only if $|S(1)|$ and $|S(-1)|$ are squares; since $S$ is irreducible, we have $G_S = 0$. Therefore
the following result is a consequence of Corollary \ref{G=0} :

\begin{coro}\label{BT corollary} Let $S$ be a Salem polynomial of degree $22$ such that $|S(1)|$ and $|S(-1)|$ are squares, and let
$\tau$ be a signature map with maximum $(3,19)$ and polynomial $S$. Then $\Lambda_{3,19}$ has an isometry
with signature map $\tau$. 

\end{coro} 

In particular, the lattice $\Lambda_{3,19}$ has an isometry
with characteristic polynomial $S$; this was proved in \cite{BT}. 

\bigskip

The following sections contain some special cases and examples that are of interest for the applications to $K3$ surfaces; the first of these
consists of reminders on Salem numbers and Salem polynomials. 

\section{Salem polynomials}\label{Salem section}

The aim of this section is to record some basic properties of Salem polynomials, and more generally, of symmetric polynomials. 

\begin{lemma}\label{mod 4}

Let $n \geqslant 1$ be an integer, and let $F \in {\bf Z}[X]$  be a symmetric polynomial of degree $2n$.
Then $F(1) \equiv (-1)^n F(-1) \ {\rm (mod \ 4)}$.



\end{lemma}

\noindent{{\bf Proof.} Let $I = \{ 1,\dots,n-1 \}$ and let $a_1,\dots, a_n$ be integers such that  $$F(X) = X^{2n} + \underset{i \in I} \sum a_i (X^i + X^{2n -i}) + a_n X^n + 1.$$

Let $I_{\rm {even}}$ be the set of $i \in I$ such that $i$ is even, and let $I_{\rm {odd}}$ be the set of $i \in I$ such that $i$ is odd.
If $n$ is even, then we have $F(1) - F(-1) = 4 \underset{i \in I_{\rm {odd}}} \sum a_i$; if 
$n$ is odd, then  $F(1) + F(-1) = 4 + 4 \underset{i \in I_{\rm {even}}} \sum a_i.$

\begin{lemma}\label{1 and -1}
 
Let $S$ be a Salem polynomial. Then $S(1) < 0$ and $S(-1) > 0$.

 \end{lemma} 
 
 \noindent{Proof.} We have $S(X) = (X-\alpha)(X - \alpha^{-1}) \underset{\rho \in P} \prod (X - \rho) (X - \overline {\rho})$,
 where $\alpha > 1$ is a real number. Note that 
for all $\rho \in P$ we have $(1-\rho)(1- \overline {\rho}) > 0$ and $(1 + \rho)(1 + \overline {\rho}) > 0$. On the other hand,
we have $(1-\alpha)(1- \alpha^{-1}) < 0$ and $(1+ \alpha)(1+  \alpha^{-1}) < 0$. Therefore we have 
$S(1) < 0$ and $S(-1) > 0$, as claimed.

\bigskip
Recall that a monic polynomial $F \in {\bf Z}[X]$ of degree $2n$ satisfies condition (C 1) if $|F(1)|$, 
$|F(-1)|$ and $(-1)^n F(1)F(-1)$ are squares.

\begin{prop}\label{GM} Let $n \geqslant 1$ be an integer, and let $S$ be a Salem polynomial of degree $2n$. 
Then we have

\medskip
{\rm (a)} If $S$ satisfies condition {\rm (C 1)}, then $n$ is odd.

\medskip
{\rm (b)}
If $|S(1)S(-1)| = 1$, then $n$ is odd.

\end{prop}

\noindent
{\bf Proof.} (a) The condition that $(-1)^n S(1)S(-1)$ is a square implies that $n$ is odd; indeed, by
Lemma \ref{1 and -1} we have $S(1) < 0$ and $S(-1) > 0$, therefore $S(1)S(-1) < 0$.

\medskip
(b) Suppose that $n$ is even; then by Lemma \ref{mod 4} we have $S(1) \equiv S(-1) \ {\rm (mod \ 4)}$. On the
other hand, Lemma \ref{1 and -1} implies that $S(1) < 0$ and $S(-1) > 0$; therefore 
$|S(1)S(-1)| > 1$.

\bigskip Note that part (b) of the lemma was proved by Gross and McMullen, see \cite{GM}, Proposition 3.3.

\begin{prop}\label{powers} If $\alpha$ is a Salem number of degree $d$, then so is $\alpha^n$, for all
integers $n$ with $n \geqslant 1$. 

\end{prop}

\noindent
{\bf Proof.} For $d \geqslant 4$, this is a result of Salem, see \cite{Sa}, page 169; see also \cite{Sm}, Lemma 2,
and it is well-known if $d = 2$. 

\begin{prop}\label{single Salem} If $K = {\bf Q}(\alpha)$ for some Salem number $\alpha$, then there 
exists a Salem number $\beta$ such that the set of Salem numbers in $K$ consists of the powers
of $\beta$. 

\end{prop} 

\noindent
{\bf Proof.} This is well-known if $d = 2$; for $d \geqslant 4$, see \cite{Sm}, Proposition 3, (iii) (and \cite{Sa}, page 169, for a slightly less general
statement).

\medskip We refer to the survey articles \cite{GH} and \cite{Sm} for more information about Salem
polynomials and Salem numbers.

\section{Examples}\label{example section}

The aim of this section and the next one  is to illustrate the results of \S \ref{sufficient section} by some examples; in view
of applications to automorphisms of $K3$ surfaces, we look for isometries of the lattice $\Lambda_{3,19}$
with characteristic polynomial $SC$, where $S$ is a Salem polynomial, and $C$ a product of cyclotomic
polynomials. We start with the case where $C$ is a {\it power} of a cyclotomic polynomial.

\begin{example}\label{Salem example} Let $F = SC$, where $S$ is a Salem polynomial
of degree $d$ with $4 \leqslant d \leqslant 20$ and $C$ is a power of a cyclotomic polynomial; assume
that ${\rm deg}(F) = 22$, and that $C(X) \not = (X-1)^2$, $C(X) \not = (X+1)^2$. If $z$ is a root of $F$ with $|z| = 1$, let $f_z \in {\bf R}[X]$ be the minimal
polynomial of $z$ over ${\bf R}$, and let $\tau_z$ be the signature map with maximum $(3,19)$
and polynomial $F$ such that ${\rm proj}_1 \circ \tau (f_z) = 2$, where ${\rm proj}_1 : {\bf N} \times {\bf N}
\to {\bf N}$ is the projection on the first factor. 

\medskip
We say that $\tau_z$ is {\it realizable} if it is the signature map of a semi-simple isometry of
the lattice $\Lambda_{3,19}$. 

\medskip
(a) Suppose that $F$ satisfies condition (C 1), and that $G_F = 0$; then $\tau_z$ is realizable
for all $z$.

\medskip  
(b) Suppose that $F$ satisfies condition (C 1), and that $G_F \not = 0$. 
We have
$I = \{S,\Phi \}$ where $\Phi$ is cyclotomic and $C$ is a power of $\Phi$. 
Let $c : I \to {\bf Z}/2{\bf Z}$ be such that $c(S) = 1$ and $c(\Phi) = 0$, and let
$c' : I \to {\bf Z}/2{\bf Z}$ be such that $c'(S) = 0$ and $c'(\Phi) = 1$. 

\medskip

The map $\epsilon^{\rm finite}_F : C(I) \to {\bf Z}/2{\bf Z}$ only depends on the polynomial $F$. 
On the other hand, the map $\epsilon^{\infty} : C(I) \to {\bf Z}/2{\bf Z}$ also depends on the choice of $z$, as
follows :

\medskip
$\bullet$ Assume that {\it $z$ is a root of $S$}. Then we have :

\medskip $d  \equiv 0 \ {\rm (mod \ 4)}$ $\implies$ $a^{\infty}(S) = 0$ and $a^{\infty}(\Phi) = 1$.

\medskip $d  \equiv 2 \ {\rm (mod \ 4)}$ $\implies$ $a^{\infty}(S) = 1$ and $a^{\infty}(\Phi) = 0$.

\medskip With $c$ and $c'$ as above, we have :

\medskip $d  \equiv 0 \ {\rm (mod \ 4)}$ $\implies$ $\epsilon_{\tau_z}^{\infty}(c) = 0$ and 
$\epsilon_{\tau_z}^{\infty}(c') = 1$.

\medskip $d  \equiv 2 \ {\rm (mod \ 4)}$ $\implies$ $\epsilon_{\tau_z}^{\infty}(c) = 1$ and 
$\epsilon_{\tau_z}^{\infty}(c') = 0$.

\bigskip
This implies that $\tau_z$ is realizable if and only if either

\medskip $d  \equiv 0 \ {\rm (mod \ 4)}$ and $\epsilon_F^{\rm finite}(c) = 0$, 
$\epsilon_F^{\rm finite}(c') = 1$, or 

\medskip $d  \equiv 2 \ {\rm (mod \ 4)}$ and $\epsilon_F^{\rm finite}(c) = 1$,
$\epsilon_F^{\rm finite}(c') = 0$.

\bigskip
$\bullet$ Assume that {\it $z$ is a root of $\Phi$}. Then we have :

\medskip $d  \equiv 0 \ {\rm (mod \ 4)}$ $\implies$ $a^{\infty}(S) = 1$ and $a^{\infty}(\Phi) = 0$.

\medskip $d  \equiv 2 \ {\rm (mod \ 4)}$ $\implies$ $a^{\infty}(S) = 0$ and $a^{\infty}(\Phi) = 1$.

\medskip With $c$ and $c'$ as above, we have :

\medskip $d  \equiv 0 \ {\rm (mod \ 4)}$ $\implies$ $\epsilon_{\tau_z}^{\infty}(c) = 1$ and 
$\epsilon_{\tau_z}^{\infty}(c') = 0$.

\medskip $d  \equiv 2 \ {\rm (mod \ 4)}$ $\implies$ $\epsilon_{\tau_z}^{\infty}(c) = 0$ and 
$\epsilon_{\tau_z}^{\infty}(c') = 1$.

\bigskip
This implies that $\tau_z$ is realizable if and only if either

\medskip $d  \equiv 0 \ {\rm (mod \ 4)}$ and $\epsilon_F^{\rm finite}(c) = 1$, 
$\epsilon_F^{\rm finite}(c') = 0$, or 

\medskip $d  \equiv 2 \ {\rm (mod \ 4)}$ and $\epsilon_F^{\rm finite}(c) = 0$,
$\epsilon_F^{\rm finite}(c') = 1$.

\bigskip
In summary, this implies that whether or not $\tau_z$ is realizable depends on $\epsilon^{\rm finite}_F$,
and whether $z$ is a root of $S$ or of $\Phi$. In particular, 

\medskip
$\bullet$ if $z$ and $z'$ are roots of the same irreducible
polynomial, then 

\medskip
\centerline {$\tau_z$ is realizable $\iff$ $\tau_{z'}$ is realizable.}

\medskip
On the other hand, 

\medskip
$\bullet$ if $z$ and $z'$ are roots of different irreducible polynomials, then 

\medskip
\centerline {$\tau_z$ realizable $\implies$ $\tau_{z'}$ is {\it not} realizable.}

\bigskip
Moreover, by Proposition \ref{epsilon}, we have $\epsilon_{\tau_z} (c) = \epsilon_{\tau_z} (c')$, and 
this implies the stronger statement

\medskip
$\bullet$ if $z$ and $z'$ are roots of different irreducible polynomials, then 

\medskip
\centerline {$\tau_z$ realizable $\iff$ $\tau_{z'}$ is not realizable.}

\bigskip
Therefore, either $\tau_z$ is realizable for all roots $z$ of $S$ with $|z| = 1$, or $\tau_z$ is realizable for all roots $z$ of $\Phi$ --  but not both. In particular, the lattice $\Lambda_{3,19}$ has a semi-simple isometry
with characteristic polynomial $S C$. 

\medskip
(c) Assume that both $S$ and $C$ satisfy condition (C 1). In this case, one gets a more precise result.
Suppose that $G_F \not = 0$ -- otherwise, every signature map is realizable by (a). We keep the notation of (b). 

\medskip
Since $S$ is irreducible, and $C$ is a power of an irreducible polynomial, we have $G_S = G_C = 0$. 
By hypothesis, $S$ satisfies condition (C 1), and this implies that $d  \equiv \pm 2 \ {\rm (mod \ 8)}$. 

\medskip
Assume first that $z$ is a root of $S$ with $|z| = 1$; then $\tau_z(S) = (3,d-3)$. Condition (C 0) holds
for $S$ and $(3,d-3)$ if and only if $d  \equiv - 2 \ {\rm (mod \ 8)}$. If this is the case, then 
Condition (C 0) also holds for $C$ and $(0,22-d)$, and
by Theorem \ref{theo special} the signature map $\tau_z$ is realizable. 

\medskip
Note that by (b), this implies that if $d  \equiv - 2 \ {\rm (mod \ 8)}$ and $z'$ is a root of $\Phi$, then $\tau_{z'}$ is not realizable.

\medskip
Assume that $z$ is a root of $\Phi$. In this case, we have $\tau_z(S) = (1,d-1)$ and $\tau_z(C) = (2,20-d)$.
Condition (C 0) holds for $S$ and $(1,d-1)$ if and only if it holds for $C$ and $(2,20-d)$, and this
happens if and only if 
$d  \equiv 2 \ {\rm (mod \ 8)}$. In this case, by Theorem \ref{theo special} the signature map $\tau_z$ is realizable.

\medskip
By (b), this implies that if $d  \equiv  2 \ {\rm (mod \ 8)}$ and $z'$ is a root of $S$, then $\tau_{z'}$ is not realizable.

\end{example}

The following examples are numerical illustrations of Example \ref{Salem example}. 

\begin{example} Let $S(X)  = X^6 - 3X^5 - X^4 + 5X^3 - X^2 - 3X +1$; this is a Salem polynomial (see for
instance  \cite{Mc1}, \S 4, or \cite{GM}, \S 7, Example 1  with $a = 3$), let $C = \Phi_{10}^4$ and $F = SC$. We keep the 
notation of Example \ref{Salem example}. We are in case (c) : the polynomials $S$ and $C$ both 
satisfy condition (C~1), and Example
\ref{a=3} shows that 
$G_F \simeq {\bf Z}/2{\bf Z}$.

\medskip If $z$ is a root of $S$ with $|z| = 1$; then Example \ref{Salem example} (c) shows that 
$\tau_z$ is realizable. On the other hand, if $z'$ is a root of $\Phi_{10}$, then $\tau_{z'}$ is not realizable.

\end{example}

\begin{example}\label{example 3} Let $S(X) = X^{14} - X^{11} - X^{10} + X^7 - X^4 - X^3 + 1$, let $C = \Phi_{20}$, and set $F = SC$. Note that $S$ is a Salem polynomial, associated to the third smallest known Salem number (see for instance \cite{Bo}). We have $G_F \simeq {\bf Z}/2{\bf Z}$. Indeed, the resultant of $S$ and $\Phi_{20}$ is equal
to $41^2$, and the common factors of $S$ and $C$ modulo $41$ are $X + 5$ and $X + 33$; these polynomials are not symmetric, hence $\Pi_{S,C} = \varnothing$. 

\medskip
We are in case (c) : the polynomals $S$ and $C$ both satisfy Condition (C 1), and  $G_F \simeq {\bf Z}/2{\bf Z}$.  Example \ref{Salem example}
shows that $\tau_z$ is realizable if $z$ is a root of $S$ with $|z| = 1$, and is not realizable if $z$ is a root of $C$. 

\end{example}

The next examples are slight variations of the same principle.

\begin{example}\label{Mc example} Let $S(X) = X^{10} + X^ 9 - X^7 - X^6 - X^5 - X^ 4 - X^ 3 + X + 1$,
the minimal polynomial of the smallest known Salem number (see for instance \cite{Bo}),
and let $C(X) = \Phi_{22}(X) (X+1)^2$; set $F = S C$, and note that $G_F \simeq {\bf Z}/2{\bf Z}$.
Indeed, ${\rm Res}(S,\Phi_{22}) = 23^2$, and the common irreducible  factors of $S$ and $\Phi_{22}$ in
${\bf F}_{23}$ are $X+4$ and $X+6$; these polynomials are not symmetric, hence $\Pi_{S,\Phi_{22}} =
\varnothing$. We also have $\Pi_{S(X),X+1} =
\varnothing$, since $S(-1) = 1$. On the other hand,  $\Phi_{22}(-1) = 11$, hence $\Pi_{\Phi_{22}(X),X+1} = 
\{11\}$,
and $G_F \simeq {\bf Z}/2{\bf Z}$.

\medskip The polynomial $F$ satisfies condition (C 1); 
with the notation of Example \ref{Salem example}
and reasoning as in part (b) of this example, we see that the signature map $\tau_z$ is realizable
either for all roots $z$ of $S$ with $|z| = 1$, or for all roots of $C$, but not both. 
It is easy to see that $\tau_z$ is realizable for all roots of $C$. Indeed, both polynomials $S$ and $C$ satisfy
condition (C 1), and 
we have $G_S = G_C = 0$. If $z$ is a root of $C$, then condition (C 0) also holds,
and hence $\tau_z$ is realizable. On the other hand, $\tau_z$ is not realizable if $z$ is a root of $S$. 

\end{example}

\begin{example}\label{intro} Let $a$ be an integer $\geqslant 0$, and set 
$$S_a(X) = X^6 - aX^5 - X^4 + (2a-1) X^3 - X^2 - aX +1.$$ The polynomial $S_a$ is a Salem polynomial
(see for instance  \cite{Mc1}, \S 4, or \cite{GM}, \S 7, Example 1 ). Set 
$$F_a(X) = S_a(X) (X^2+1)^2(X-1)^{12}.$$ By Example \ref{Sa}, we have 
$G_{F_a} = 0$ if $a \geqslant 1$, and $G_{F_0}$ is of order $2$. Therefore  if $a \geqslant 1$, then every
signature map of maximum $(3,19)$ and polynomial $F_a$ is realizable. On the other
hand, if $a = 0$, then the method of the previous examples shows that 
$\tau$ is realizable if and only if $\tau = \tau_z$ where $z$ is a root of $S$ with
$|z| = 1$; equivalently,  if we have $\tau(S) = (3,3)$.

\end{example}

\medskip
{\bf Part II : Salem signature maps and surgery on lattices}

\medskip
This part  is the bridge between the arithmetic results of Part I, and the applications to $K3$ surfaces
of Part III. The strategy is the following :

\medskip
Let $t$ be a semi-simple isometry of an even unimodular lattice $(L,q)$ of signature $(3,19)$
with characteristic polynomial $SC$, where $S$ is a Salem polynomial of degree $d$ with
$4 \leqslant d \leqslant 20$, and $C$ is a product of cyclotomic polynomial. Let $\alpha$ be
the Salem number corresponding to $S$. The criterion of
McMullen  \cite{Mc2}, Theorem 6.2 (see also \cite{Mc3}, Theorem 6.1) determines whether
$t$ is induced by an automorphism of a $K3$ surface.

\medskip
Set $L_S = {\rm Ker}(S(t))$ and $L_C = {\rm Ker}(C(t))$;  then $(L_C,q)$ is either negative definite,
or of signature $(2,{\rm deg}(C)-2)$. As in the introduction, we 
{{\it assume that $(L_C,q)$ is
negative definite}. 
In this case, McMullen's criterion concerns the restriction of $t$ to $L_C$ : namely,
it is satisfied if and only if $t$ stabilizes a chamber of the root system of $L_C$ (see \cite{Mc2}, \S 5, \S 6 
and \cite{Mc3}, \S 2, for details; see also Theorem \ref{Mc6.2}). If this criterion is fulfilled, then we are done~:
$\alpha$ is realizable. If not, we replace $t$ by an isometry $t'$ satisfying McMullen's 
criterion, of characteristic polynomial $SC'$, where $C'$ is a product of cyclotomic polynomials. 

 

\section{Salem signature maps and realizable isometries}\label{Salem signature section}

We start with some definitions and notation.

\begin{defn}\label{Salem signature definition} A signature map $\tau$ is called a {\it Salem signature map}
if its associated polynomial is of the form $F = SC$, where
$S$ is a Salem polynomial of degree $d$ with $4 \leqslant d \leqslant 22$ and $C$ is a product of 
cyclotomic polynomial, and moreover

\medskip
$\bullet$ $\tau(S) = (3,d-3)$ and $\tau(C) = (0,{\rm deg}(C))$. 

\medskip

$\bullet$ $\tau$ satisfies conditions (C 0) and (C 1).

\medskip
The Salem polynomial $S$ is called the {\it Salem polynomial of $\tau$}, and the Salem number of the polynomial $S$ is called the {\it Salem number of 
$\tau$.}

\end{defn}




\begin{defn}\label{Salem equivalence definition} Let $S$ be a Salem polynomial. Two signature
maps $\tau$ and $\tau'$ are said to be {\it equivalent} at $S$ if 

\medskip
\centerline {$\tau(f) = \tau'(f)$ for all factors $f \in {\bf R}[X]$ of $S$.}

\end{defn}

\begin{notation}\label{z signature map} Let $S$ be a Salem polynomial of degree $d$ with $4 \leqslant d \leqslant 22$, and 
let $z \in {\bf C}$ be a root of $S$ with $|z| = 1$. Set $f_z(X) = (X-z)(X-z^{-1})$. We denote by $\tau_{S,z}$ the signature map of
maximum $(3,d-3)$ and polynomial $S$ characterized by $\tau_{S,z}(f_z) = (2,0)$. 

\end{notation}

\medskip

We next define the notion of {\it realizable isometry}.

\begin{defn}\label{type P} If $(N,q)$  is an even, negative definite lattice, let
$$\Delta(N) = \{x \in N \ | \ q(x,x) = -2 \}$$ be the set of {\it roots} of $N$; it is a root system
(see for instance \cite{E}, Theorem 1.2). 
Set $N_{\bf R} = N \otimes_{\bf Z} {\bf R}$. For all $x \in \Delta(N)$, let  $H_x$ be the hyperplane of
$N_{\bf R}$ orthogonal to $x$, and let $H$ be the union of the hyperplanes $H_x$ for $x \in \Delta(L)$. The
complement of $H$ in $N_{\bf R}$ is a  disjoint union of connected open subsets, called the {\it chambers} of $N_{\bf R}$ with
respect to $q$. 
An isometry $t : N \to N$ is said to be of {\it type} (P) if $t$ stabilizes a chamber of
$N_{\bf R}$ with
respect to $q$.
\end{defn}

\begin{defn}\label{realizable isometry} Let $t$ be a semi-simple isometry of an even unimodular lattice $(L,q)$ 
with characteristic polynomial $SC$, where 
$S$ is a Salem polynomial of degree $d$ with
$4 \leqslant d \leqslant 20$, and $C$ is a product of cyclotomic polynomials. 
Set $L_S = {\rm Ker}(S(t))$ and $L_C = {\rm Ker}(C(t))$; suppose that $(L_C,q)$ is
negative definite. We say that $t$ is {\it realizable} if the restriction of $t$ to $L_C$ is of type (P).

\end{defn}

Let  $S$  be a Salem polynomial of degree $d$ with $4 \leqslant d \leqslant 22$, and let $z \in {\bf C}$ be a root of $S$ with $|z| = 1$.

\begin{theo}\label{realizable theorem} Suppose that  $d \leqslant 16$ and $d  \equiv 0, 4 \ {\rm or} \ 6  \ {\rm (mod \ 8)}$.
Then $\Lambda_{3,19}$ has
a realizable isometry with signature map equivalent to $\tau_{S,z}$.

\end{theo}


We start with some preliminary results.








\begin{lemma}\label{-id} Let $(L,q)$ be an even negative definite lattice of rank $2$, and set $V = L \otimes_{\bf Z} {\bf Q}$. The lattice $(L,q)$ has an isometry $t : L \to L$  of type {\rm (P)} such that $$\partial_v[V,q,t] =
\partial_v[V,q,-id]$$ for all $v \in \mathcal V'$.

\end{lemma}

\noindent
{\bf Proof.} 
There exist integers $D \geqslant 1$ and $f \geqslant 1$  such that 
${\rm det}(L,q) = f^2D$, where $-D$ is the discriminant of an imaginary quadratic field. 
If $(L,q)$ contains no roots, then we are done : we can take $t = -id$. Otherwise, the lattice $(L,q)$  is isomorphic to $-q'$, where $q'$ is a quadratic form associated to an 
order $O$ of the imaginary quadratic field ${\bf Q}(\sqrt {-D})$ (see for instance \cite {Co}, Theorem 7.7 (ii)). Complex
conjugation induces an isometry of the quadratic form $(O,q')$ with characteristic polynomial $X^2 - 1$. If $D = 3$ and $f = 1$, then
$(O,q')$  is isomorphic to the root lattice $A_2$, and complex conjugation is an isometry of type (P) of $(O,q')$ (see \cite {Mc2}, \S 5, Example); if $D = 4$ and $f = 1$,
then $(O,q')$ is isomorphic to $A_1 \oplus A_1$, and exchanging the two copies of $A_1$ gives rise to an isometry of type (P) of $(O,q')$; 
otherwise, $(O,q')$ contains only two roots,  fixed by complex conjugation, hence we obtain an isometry $t$ of 
type (P) of $(O, q')$ in
this case as well, and it is clear that $\partial_v[V,q,t] =
\partial_v[V,q,-id]$ for all $v \in \mathcal V'$.

\begin{prop}\label{not both squares} Suppose that  $d \leqslant 18$ and that one of the following holds :

\medskip {\rm (i)} $|S(1)|$ and $S(-1)$ are not both squares. 

\medskip {\rm (ii)} $S$ is ramified at the prime $2$.

\medskip
Then $\Lambda_{3,19}$ has
a realizable isometry with signature map equivalent to $\tau_{S,z}$.

\end{prop}

\noindent
{\bf Proof.} (a) Suppose first  that  there exists a prime number $p$ such that $v_p(S(-1)) \equiv 1 \ {\rm (mod \ 2)}$. If $p$ is odd, then  by Lemma \ref{odd ramification} we
have $p \in \Pi^{r,-}_S$; if $p = 2$, then  $2 \in \Pi^{r}_S$ by Proposition \ref{ramification at 2}. Set $F(X) = S(X)(X-1)^{20-d}(X+1)^2$. We have $D_- = (-1)^{s_-}S(-1) = S(-1)$; this implies that $D_- \not = -1$ in 
${\bf Q}_p^{\times}/{\bf Q}_p^{\times 2}$, since $v_p(D_-) = v_p(S(-1)) \equiv 1 \ {\rm (mod \ 2)}$. Therefore $p \in \Pi_{S,X+1}$, and $G_F(D_-) = 0$.

\medskip (b) Suppose now that no such prime number exists; this implies that $S(-1)$ is a square. If we are in case (i), then this implies that $|S(1)|$ is  not a square, therefore 
there exists a prime number $p$ such that $v_p(S(-1)) \equiv 1 \ {\rm (mod \ 2)}$. By the same argument as before, we have $p \in \Pi^{r,+}_S$ if $p$ is odd,
and $2 \in \Pi^{r}_S$ if $p = 2$. If we are in case (ii), then by hypothesis we have $2 \in \Pi^{r}_S$. Set $F(X) = S(X)(X-1)^{22-d}$. We have $d \leqslant 18$, hence $22-d \not = 2$; this implies that $p \in \Pi_{S,X-1}$ and hence
$G_F = 0$. 

\medskip

Let $\tau$
be the  signature map with polynomial $F$ and maximum $(3,19)$ such that
$\tau(f) = \tau_{S,z}(f)$ for all factors $f \in {\bf R}[X]$ of $S$, and note that Conditions (C 0) and (C~1) hold for $\tau$. 
Let $(L,q)$ be an even unimodular lattice of signature $(3,19)$, and let $t : L \to L$ be a semi-simple
isometry with signature map $\tau$; such an isometry exists by Theorem \ref{theo special}, since 
$G_F(D_-) = 0$ in case (a) and $G_F = 0$ in
case (b).

\medskip
Set $L_S = {\rm Ker}(S(t))$, and let $S_C$ be the orthogonal of $L_S$ in $L$. If we are in case (b),
then the restriction of $t$ to $L_C$ is the identity, hence it is of type (P); this implies that
$t : L \to L$ is realizable.

\medskip
Suppose now that we are in case (b). Let $L_1 = L_S = {\rm Ker}(S(t))$,
$L_2 = {\rm Ker}(t+1)$ and $L_3 = {\rm Ker}(t-1)$; let $t_i : L_i \to L_i$ be the restriction of $t$ to $L_i$. 
Since $(L_2,q)$ is of rank $2$, even and negative
definite, we can apply Lemma \ref{-id}; 
let $t_2' : L_2 \to L_2$ be the isometry of $(L_2,q)$ of type (P) provided by this lemma.
Then $(L,q)$ has a semi-simple isometry $t' : L \to L$ inducing  the 
isometry $t_2$ or $t_2'$ on
$L_2$ and $t_i$ on $L_i$ for $i = 1,3$; hence $t'$ is realizable. This completes the proof of the proposition.

\medskip

\noindent{\bf Proof of Theorem \ref{realizable theorem}.}  If $|S(1)|$ and $S(-1)$ are not
both squares, then the result follows from Proposition \ref{not both squares} (i).

\medskip Assume that $|S(1)|$ and $S(-1)$ are both squares, and set $n = d/2$. 
Since $S$ is a Salem polynomial, we have $S(1) < 0$
and $S(-1) > 0$.

\medskip
Suppose first that $d  \equiv 0 \ {\rm (mod \ 4)}$; then $n$ is even, and hence $(-1)^{n}S(1)S(-1) <0$. If $S$ is unramified at $2$, then by
Proposition \ref{ramification at 2} 
the class of $(-1)^n S(1) S(-1)$ in ${\bf Q}_2^{\times}/{\bf Q_2}^{\times 2}$ lies in $\{1,-3 \}$, and this leads to a contradiction. Therefore $S$
is ramified at $2$, and 
Proposition  \ref{not both squares} (ii) gives
the desired result.

\medskip Assume now that $d  \equiv 6 \ {\rm (mod \ 8)}$. Then
$n$ is odd, hence $(-1)^{n}S(1)S(-1) > 0$. By hypothesis,
$|S(1)|$ and $S(-1)$ are squares, therefore condition (C 1) holds for $S$.

\medskip Set $(r,s) = (3,d-3)$; note that the hypothesis $d  \equiv 6 \ {\rm (mod \ 8)}$ implies that
$r  \equiv s \ {\rm (mod \ 8)}$. Since $S$ is irreducible, we have $G_S = 0$; therefore by Corollary
 \ref{G=0}
the lattice $\Lambda_{3,d-3}$ has an isometry with signature map $\tau_{S,z}$. 
The identity is a semi-simple isometry of the lattice $-E_8$; since $\Lambda_{3,19} =
\Lambda_{3,3} \oplus (-E_8) \oplus (-E_8) = \Lambda_{3,11} \oplus (-E_8)$, we obtain a realizable
isometry, as claimed. 

\begin{prop}\label{22 theo} Suppose that $d = 22$ and that $|S(1)|$ and $S(-1)$ are both squares. 
Then  $\Lambda_{3,19}$ has
a realizable isometry with signature map equivalent to $\tau_{S,z}$.

\end{prop}

\noindent
{\bf Proof.} This is an immediate consequence of Corollary \ref{BT corollary}.

\medskip If $d = 10$ or $18$ and if $|S(1)|$ and $S(-1)$ are not  both squares, then 
$\Lambda_{3,19}$ has
a realizable isometry with signature map equivalent to $\tau_{S,z}$ by Proposition \ref{not both squares}. The aim of 
the next sections is to give a necessary and sufficient condition in the case where $|S(1)S(-1)| = 1$, and
more generally, when $S$ is {\it unramified}. We start with a definition. 

\medskip Recall that if $\tau$ is a signature map, we use the notation $\tau(X-1) = (r_+,s_+)$ and $\tau(X+1) = (r_-,s_-)$.
Note that if $\tau$ is a Salem signature map, then $s_+$ and $s_-$ are both even, and hence $D_+ = |F_1(1)|$, $D_- = |F_1(-1)|$. 

\begin{defn}\label{signature map with trivial obstruction}
Let $\tau$ be a Salem signature map with polynomial $F$. We say that $\tau$ is a {\it Salem signature map with trivial obstruction} if
 $G_F(D_+,D_-) = 0$.

\end{defn}

\begin{theo}\label{condition 10} Suppose that $d = 10$ or $18$ and that $S$ is unramified. Suppose that there exists
a Salem signature map with trivial obstruction of maximum $(3,19)$ or $(3,11)$ and Salem factor $S$. Then the lattice 
$\Lambda_{3,19}$ has
a realizable isometry with signature map equivalent to $\tau_{S,z}$.

\end{theo}

This result will be proved in the next sections.

\section {\bf Preparing for surgery}\label{surgery section}

Let $S$ be a Salem polynomial of degree $d$ with $4 \leqslant d \leqslant 18$ such that
$d  \equiv 2 \ {\rm (mod \ 4)}$. 
We assume that $S$ is {\it unramified}; for instance, this is the case if $|S(1)S(-1)| = 1$. The aim of
this section is to prove a preliminary result that will be useful in the proof of Theorem \ref{condition 10} in
the next sections.

\medskip Let $\Phi$ be a cyclotomic polynomial of even degree, let $k \geqslant 1$ be an integer, and
set $C = \Phi^k$. Assume that ${\rm deg}(C) \equiv 4 \ {\rm (mod \ 8)}$, and that 
condition (C 1) holds for $C$. 
Recall that if $t : L \to L$ is a semi-simple isometry
with characteristic polynomial $SC$, we set $L_S = {\rm Ker}(S(t))$ and $L_C = {\rm Ker}(\Phi(t))$,
$V = L \otimes_{\bf Z}{\bf Q}$, 
$V_S = L_S \otimes_{\bf Z}{\bf Q}$ and $V_C = L_C \otimes_{\bf Z}{\bf Q}$.
Recall that $\mathcal V$ is the set of places of ${\bf Q}$ and that $\mathcal V'$ is
the set of finite places.

\begin{prop}\label{unique p} Let $\tau$ be a Salem signature map with trivial obstruction and polynomial $SC$, and let $p \in \Pi_{S,\Phi}$.
 Then there exists an even unimodular lattice $L$
having a semi-simple isometry with characteristic polynomial $F$ and signature map $\tau$ such
that $\partial_v([V_S]) = \partial_v([V_C]) = 0$ for all $v \in \mathcal V'$ with $v \not = v_p$.

\end{prop}

\noindent
{\bf Proof.} Let $L$ be an even unimodular lattice and let 
$t : L \to L$ be a semi-simple isometry with signature map $\tau$; 
these exist by Theorem \ref{theo special}.
Let $V = L \otimes_{\bf Z}{\bf Q}$,
and note that $\partial_v([V]) = 0$ for all $v \in \mathcal V'$. We have $V = V_S \oplus V_C$, therefore
for all $v \in \mathcal V'$, we have $$\partial_v([V_S]) = 0 \iff \partial_v([V_C]) = 0.$$

\medskip Let $\mathcal U = \mathcal U(L,t)$ be the set of $v \in \mathcal V'$ such that $\partial_v([V_S]) \not = 0$.
Let $v \in \mathcal U$; this implies that $v = v_q$ with $q \in \Pi_{S,\Phi}$. 
Since $S$ is unramified, 
the common, irreducible, symmetric factors of $S$ and $\Phi$  $({\rm mod} \ q)$ are
of even degree; hence $\partial_v([V_S])$ is contained in a product of Witt groups of the
type $W(\kappa,\sigma)$, where $\kappa$ is an even degree extension of ${\bf F}_q$  and 
$\sigma$ is a non-trivial involution of $\kappa$. Since these Witt groups are of order $2$, the (non-trivial) element
$\partial_v([V_S])$ is of order $2$; moreover, $\partial_v([V_S]) = \partial_v([V_C]) $.

\medskip

Note that $I = \{S,\Phi\}$, and recall that the map $(S,\Phi) : C(I) \to C(I)$ was defined in 
\S \ref{combinatorial section}, cf. Notation \ref{f,g}.  Set $E_S = {\bf Q}[X]/(S)$, $E_{\Phi} = {\bf Q}[X]/(\Phi)$ and $E = E_S \times E_{\Phi}$.
Let $(\lambda^v)$ be
the local data associated to $V$, and for all $v \in \mathcal V$, let $a^v = a(\lambda^v) \in \mathcal C^v$. 
Since $V_S$ and $V_C$ are global forms, we have $\underset {v \in \mathcal V} \sum a^v(f) = 0$ for
all $f \in I$.

\medskip
Let $N$ be the cardinality of $\mathcal U$. Assume first that $N \geqslant 2$; we claim that
in this case there exists an even unimodular lattice  $L'$ and 
a semi-simple isometry $t : L' \to L'$  with signature map $\tau$ such that the cardinal of  $\mathcal U(L',t')$ is $N-2$.

\medskip Let $v$ and $v'$ be two distinct elements of $\mathcal U$.
Set $b^v = (S,\Phi)(a^v)$ and $b^{v'} = (S,\Phi)(a^{v'})$; we have $b^v \in \mathcal C^v$ and $b^{v'} \in \mathcal C^{v'}$; note that $ b^v(f) + b^{v'} (f) = a^v(f) + a^{v'} (f)$
for all $f \in I$. 
Let $\mu^v, \mu^{v'}  \in T(E^v,\sigma)$  be such that $b^v = a(\mu^v)$ and $b^{v'} = a(\mu^{v'})$; let
us chose $\mu^v = (\mu^v_S,\mu^v_{\Phi})$ and  $\mu^{v'} = (\mu^{v'}_S,\mu^{v'}_{\Phi})$ such 
that $\partial_v(\mu^v_S) = 0$ and $\partial_{v'}(\mu^{v'}_S) = 0$ (hence also 
$\partial_v(\mu^v_{\Phi}) = 0$ and $\partial_{v'}(\mu^{v'}_{\Phi}) = 0$). 
For all $w \in \mathcal V$ with $w \not = v, v'$, set $\mu^{w} = \lambda^{w}$, and let 
$b^{w} = 
a(\mu^{w}) \in \mathcal C^w$; we have $\underset {w \in \mathcal V} \sum b^w(f) = \underset {w \in \mathcal V} \sum a^w(f)$
for all $f \in I$. This implies that $\underset {w \in \mathcal V} \sum b^w(f) = 0$ for all $f \in I$.
By Theorem \ref{realization}  there exists  an even unimodular lattice  $L'$ and 
a semi-simple isometry $t : L' \to L'$  with signature map $\tau$ and
local data $(\mu^w)$. Set $V' = L' \otimes_{\bf Z}{\bf Q}$. We have 
$\partial_v([V'_S]) = \partial_v(\mu^v_S)$, $\partial_{v'}([V'_S]) = \partial_{v'}(\mu^{v'}_S)$,
hence $\partial_v([V'_S]) = \partial_{v'}([V'_S]) = 0$. Since $\partial_w([V'_S]) =
\partial_w([V_S])$ if $w \not = v,v'$, this implies that
the cardinal of  $\mathcal U(L',t')$ is $N-2$, as claimed. 

\medskip Therefore we can assume that either $N = 0$ or $N = 1$.

\medskip
Note that the signature of $V_S$ is $(3,d-3)$, and that $V_C$ is negative definite. Hence
$V_S$ and $V_C$ do not contain even unimodular lattices. This implies that there exists at least one
$v \in \mathcal V'$ such that $\partial_v([V_S]) \not = 0$ and $\partial_v([V_C]) \not = 0$, and this
implies that $N \not = 0$.

\medskip Assume that $N = 1$, and let $\mathcal U = \{v\}$; if $v = v_p$, we are finished. Suppose
that $v \not = v_p$. 
Set $b^v = (S,\Phi)(a^v)$ and $b^{v_p} = (S,\Phi)(a^{v_p})$; we have $b^v \in \mathcal C^v$ and $b^{v_p} \in \mathcal C^{v_p}$; note that $ b^v(f) + b^{v_p} (f) = a^v(f) + a^{v_p} (f)$
for all $f \in I$. 
Let $\mu^v, \mu^{v_p}  \in T(E^v,\sigma)$  be such that $b^v = a(\mu^v)$ and $b^{v_p} = a(\mu^{v_p})$;
choose $\mu^v,$ so that $\partial_v(\mu^v_S) = 0$. 
For all $v' \in \mathcal V$ with $v' \not = v, v_p$, set $\mu^{v'} = \lambda^{v'}$, and let 
$b^{v'} = 
a(\mu^{v'}) \in \mathcal C^v$; we have $\underset {w \in \mathcal V} \sum b^w(f) = \underset {w \in \mathcal V} \sum a^w(f)$
for all $f \in I$. This implies that $\underset {w \in \mathcal V} \sum b^w(f) = 0$ for all $f \in I$.
By Theorem \ref{realization}  there exists an even unimodular lattice  $L'$ and 
a semi-simple isometry $t : L' \to L'$  with signature map $\tau$ and local data $(\mu^w)$. Set $V' = L' \otimes_{\bf Z}{\bf Q}$. We have 
$\partial_v([V'_S]) = \partial_v(\mu^v_S)$, 
hence $\partial_v([V'_S]) = \partial_{v'}([V'_S]) = 0$. Since $\partial_w([V'_S]) =
\partial_w([V_S])$ if $w \not = v,v_p$, this implies that $\partial_w([V'_S]) = 0$ for all
$w \in \mathcal V'$ with $w \not = v_p$. This completes the proof of the proposition.

\section{Unramified Salem polynomials of degre $10$ or $18$}\label{18 section}

Let $S$ be an unramified Salem polynomial of degree $d$, with $d = 10$ or $18$. 

\begin{theo}\label{10 and 18} Let $\tau$ be a Salem signature map with trivial obstruction of maximum $(3,d+1)$ and Salem factor $S$. Then
the lattice $\Lambda_{3,19}$ has a realizable isometry with signature map equivalent to $\tau$. 
Moreover, if $d = 10$, then this signature map can be chosen so that its polynomial is divisible by $(X-1)^8$. 

\end{theo}

\noindent {\bf Proof.} If $m \geqslant 1$ is an integer, we denote by $\Phi_m$ the $m$-th cyclotomic polynomial. Let $F$ be the polynomial associated to $\tau$. We have
$G_F = 0$ by definition; since $|S(1)S(-1)| = 1$, this implies that the polynomial $F$ is divisible by a cyclotomic polynomial
$\Phi_m$ with $m = 3,4,6$ or $12$ such that $\Pi_{S,\Phi_m} \not = \varnothing$. Let $\Phi_m$ be such
a polynomial, and let $p$ be a prime number such that $p \in \Pi_{S,\Phi_m}$. Let $f \in {\bf F}_p[X]$ be
an irreducible, symmetric common factor of $S$ and $\Phi_m$ $({\rm mod} \ p)$; the polynomial $f$ has
even degree. Let $M = {\bf F}_p[X]/(f)$, considered as a simple ${\bf F}_p[\Gamma]$-module, and let $[M]$ denote the only non-trivial element of $W_{\Gamma}({\bf F}_p,M)$. 

\medskip Set $C = \Phi_m^2$ if $m = 3,4$ or $6$ and $C = \Phi_m$ if $m = 12$, and let $\tau'$ be
the Salem signature map of maximum $(3,d+1)$ and associated polynomial $SC$ such that $\tau'(g) =
\tau(g)$ for all factors $g \in {\bf R}[X]$ of $S$. 
Let $t : L \to L$ be a semi-simple isometry with signature map $\tau'$ of an even unimodular lattice $L$
such that that $\partial_v([V_S]) = \partial_v([V_C]) = 0$ for all $v \in \mathcal V'$ with $v \not = v_p$; these
exist by Proposition \ref{unique p}. Note that $\tau'$ is equivalent to $\tau$.

\medskip Suppose first that $m = 4$, and note that with the notation above, the polynomial $f$ is the
image of $\Phi_4$ in ${\bf F}_p[X]$, and that $p  \equiv 3 \ {\rm (mod \ 4)}$. We now construct a negative definite lattice $N$ and a semi-simple isometry
$t_N : N \to N$ of type (P) such that 

$$\partial_v([N \otimes_{\bf Z} {\bf Q}]]) = \partial_v([V_C])$$ for all $v \in \mathcal V'$.

\medskip
Let $\mathcal O$ be the ring of integers of the quadratic field ${\bf Q}(\sqrt {-p})$, and let
$q_{\mathcal O} : \mathcal O \times \mathcal O \to {\bf Z}$ be the associated binary
quadratic form; $q_N$ is even, positive definite, and has determinant $p$. Let
$N = \mathcal O \oplus \mathcal O$ and let $q_N = (-q_N) \oplus (-q_N)$. For
all $x \in \mathcal O$, we denote by $x \mapsto \overline x$ the complex conjugation. Let
$t_N : N \to N$ be the isometry of $q_N$ given by $t_N(x,y) = (y,\overline x)$ for
all $x, y \in \mathcal O$. The characteristic polynomial of $t_N$ is $X^4 - 1 = (X^2 -1)(X^2+1)$,
and $t_N$ is of type (P). We have 
$\partial_v([N \otimes_{\bf Z} {\bf Q}]]) = \partial_v([V_C])$ for all $v \in \mathcal V'$, as
claimed; hence
 $\partial_v([N \otimes_{\bf Z} {\bf Q}]) = \partial_v([V_S])$ for all $v \in \mathcal V'$.
 
 \medskip
Set $V' = V_S \oplus (N \otimes_{\bf Z} {\bf Q})$, and let $t' : V' \to V'$ be the isometry
such that the restriction of $t'$ to $V_S$ is $t$, and the restriction to $N$ is $t_N$. The lattice
$L_S \oplus N$ is stabilized by $t'$. 
Note that $(L_S \oplus N) \otimes_{\bf Z} {\bf Z}_q$ is unimodular if $q \not = p$, and it is even
for $q = 2$. 

\medskip
Since $\partial_{v_p}([N \otimes_{\bf Z} {\bf Q}]]) = \partial_{v_p}[V_S])$, we have $\partial_{v_p}[V']) = 0$,
therefore $V'\otimes_{\bf Q}{\bf Q}_p$ contains a unimodular $\bf Z_p$-lattice stabilized by $t'$; let $L'_p$ be such
a lattice. If $q$ is a prime number $\not = p$, set $L'_q = (L_S \oplus N) \otimes_{\bf Z}{\bf Z_q}$, and let 
$L' = \underset{q} \cap L'_q$, where $q$ runs over all the prime numbers. The lattice $L'$ is even,
unimodular, and 
 the signature map of $t'$ is equivalent to $\tau'$. 

\medskip
If $d = 18$, then $L'$ is of rank $22$, isomorphic to $\Lambda_{3,19}$, hence the theorem is proved in this case.
Assume that $d = 10$; in this case, $L'$ is of rank $14$. Set $L'' = L' \oplus (-E_8)$, and
let $t'' : L'' \to L''$ be the isometry such that the restriction of $t''$ to $L'$ is equal to $t'$, and
that $t''$ is the identity on $-E_8$. 
The lattice $L''$ is isomorphic to $\Lambda_{3,19}$,  the isometry $t''$ is realizable, and its signature map is equivalent to $\tau'$, hence the theorem is
proved in this case as well.

\medskip

Suppose that $m = 3$. 
Since $S$ is unramified,  the image of the polynomial 
$\Phi_3$ in ${\bf F}_p[X]$ is the common irreducible factor of $S$ and $\Phi_3$ 
(${\rm mod}$ $p$). With the notation above, we have 
$\partial_{v_p}([V_S]) = \partial_{v_p}([V_C]) = [M]$, and 
$\partial_{v}([V_S]) = \partial_{v}([V_C]) = 0$ if $v \not = v_p$. 
If $L_C$ has no roots, then the restriction of $t$ to $L_C$ is of type (P); suppose that $L_C$ has roots.

\medskip Assume first that $p = 2$; in this case,  we have $L_C \simeq -D_4$, where $D_4$ is the
lattice associated to the root system $D_4$ (see for instance \cite {E}). Let $t' : L_C \to L_C$
be the trialitarian isometry of order $3$ of $D_4$, and note that $t'$
is of type (P). Let $V$ denote the
${\bf Q}[\Gamma]$-bilinear form on $L_C \otimes_{\bf Z}{\bf Q}$ obtained by letting a generator $\gamma$
of $\Gamma$ to act via $t'$; we have  $\partial_{v_2}([V]) = [M] =  \partial_{v_2}([V_C])$,
and hence $\partial_{v_2}([V]) = \partial_{v_2}([V_S])$. Let $t' : V_S \oplus V \to V_S \oplus V$
be defined by $t'(x,y) = (t(x),t'(y))$ for all $x \in V_S$ and $y \in V$. Since
$\partial_{v_2}([V]) = \partial_{v_2}([V_S])$, we have $t'(L) = L$. The isometry $t' : L \to L$ is
realizable, and its signature map is similar to $\tau'$.

\medskip
Suppose now that $p \not = 2$, and note that $p \not = 3$. Still assuming that $L_C$ has roots, this implies that  $L_C$ contains a sublattice of rank $2$ isomorphic to the
root lattice $A_2$. Let $L_0 \subset L_C$ be this sublattice, and let $L_1 \subset L_C$ be its
orthogonal complement in $L_C$; we have $t(L_1) = L_1$. Set $V_0 = L_0 \otimes_{\bf Z} {\bf Q}$ and
$V_1= L_1 \otimes_{\bf Z} {\bf Q}$.
We have $\partial_{v_p}[V_0] = 0$ and $\partial_{v_3}[V_0]= [({\bf F}_3,-1,id)]$.
This implies that $\partial_v[V_1] = 0$ for $v \not = v_3, v_p$, that
$\partial_{v_3}[V_1] =  [({\bf F}_3,1,id)]$, and $\partial_{v_p}[V_1] = \partial_{v_p}[V_C]$. 
Let $t' : L_C \to L_C$ be determined by $t'|L_0 = id$ and $t'|L_1 = t$. The isometry $t'$ is of type (P),
and determines an isometry $t' : L \to L$ of signature map $\tau'$.

\medskip The argument is the same for $m = 6$. In both cases, if $d = 18$, then $L$ is isomorphic
to $\Lambda_{3,19}$, hence the theorem is proved. If $d = 10$, then the same argument as in the
case $m = 4$ provides a realizable isometry of $\Lambda_{3,19}$ with signature map equivalent to $\tau'$.

\bigskip  Assume now that $m = 12$. If $L_C$ has roots, then either $p = 2$ and $L_C$ is isomorphic to  $- D_4$, or $p = 3$, and $L_C$ is isomorphic to $- A_2 \oplus - A_2$. In the first case, let
$t' : L_C \to L_C$ be obtained from the trialitarian isometry of order 3 of $D_4$; this isometry is
of type (P). We have $\partial_v[V_C,t'] = \partial_v[V_C,t]$
for all $v \in \mathcal V'$, therefore the lattice $L$ has a realizable isometry $t' : L \to L$ with
signature map equivalent to $\tau'$. If $p = 3$
and $L_C \simeq -A_2 \oplus - A_2$, then we have $\partial_{v_3}[V_C] = [{(\bf F}_3^2,q,\overline t)]$,
the characteristic polynomial of $\overline t$ being the image of $\Phi_4$ in ${\bf F}_3[X]$.
We identify the lattice $A_2$ to ${\bf Z}[\zeta_3] \times {\bf Z}[\zeta_3] \to {\bf Z}$, sending $(x,y)$
to  ${\rm Tr}_{{\bf Q}(\zeta_3)}(x \overline y)$, where $x \mapsto \overline x$
is the complex conjugation of the cyclotomic field ${\bf Q}(\zeta_3)$. Let
$t' : L_C \to L_C$ be  given  by $t'(x,y) = (y,\overline x)$ for $x,y \in A_2$; we have $\partial_{v_3}[V_C,t'] = \partial_{v_3}[V_C,t]$. The lattice $L$
has a realizable isometry with signature map equivalent to $\tau'$. If $d = 18$, then $L$ is isomorphic
to $\Lambda_{3,19}$, hence the theorem is proved. If $d = 10$, then the same argument as in the
case $m = 4$ gives us a realizable isometry of $\Lambda_{3,19}$ with signature map equivalent to $\tau'$. 

\section{Isometries of root lattices}\label{root section}

In the following section, we need more information on isometries of root lattices of characteristic polynomial
a power of a cyclotomic polynomial. This question was studied in several papers, in particular 
 \cite{BM}, th\'eor\`eme A.2. and \cite{B 99}, \S 3 (see also \cite{BB}, \cite{B 84}. 

\begin{prop} Let $\ell$ be an odd prime number. There exists a root lattice $\Lambda$  having an isometry
with minimal polynomial $\Phi_{\ell}$ or $\Phi_{2 \ell}$ if and only if the irreducible components of $\Lambda$
are 

\medskip
{\rm (a)} of type $A_{\ell -1}$.

\medskip
{\rm (b)} of type $A_{4}$ or $E_8$ if $\ell = 5$.

\medskip
{\rm (c)} $A_{2}$, $D_4$, $E_6$ or $E_8$ if $\ell = 3$.

\end{prop}

\noindent {\bf Proof.} See \cite{BM}, th\'eor\`eme A.2.

\begin{prop}\label{root lattice lemma}  Let $\ell$ be an odd prime, and let $\Lambda$ be a root lattice of type $A_{\ell-1}$ having
an isometry with characteristic polynomial $\Phi_{\ell}$ {\rm (}resp. $\Phi_{2 \ell}${\rm )}. Set $V = \Lambda \otimes_{\bf Z} {\bf Q}$ and let $[V] \in W_{\Gamma}({\bf Q})$ be the  Witt class of the ${\bf Q}[\Gamma]$-form $V$. Then $\partial_p[V] = 0$ if $p$ is a prime number
with $p \not = \ell$, and

$$\partial_{\ell}[V] = [({\bf F}_{\ell}, \langle 1 \rangle, id)],$$ resp.

$$\partial_{\ell}[V] = [({\bf F}_{\ell}, \langle 1 \rangle, - id)].$$ 

\end{prop}

\medskip
\noindent
{\bf Proof.} This is an easy computation. 

\medskip The following observation will be used several times :

\begin{prop}\label{equal residue} Let $L$ be a negative definite lattice with an isometry $t : L \to L$,
$V = L \otimes_{\bf Z} {\bf Q}$. Let $\ell$ be an odd prime number. Let $\Lambda \subset
L$ be a lattice stabilized by $t$, assume that $\Lambda$ is isometric to an orthogonal sum of
 lattices $-A_{\ell - 1}$, and that the
minimal polynomial of the restriction of $t$ to $\Lambda$ is  $\Phi_{\ell}$ or $\Phi_{2 \ell}$. Let $V'$ be the 
orthogonal of $\Lambda \otimes_{\bf Z} {\bf Q}$ in $V_L$. Then

$$\partial_p[V] = \partial_p[V']$$ for all prime numbers $p \not = \ell$.

\end{prop}

\noindent{\bf Proof.} Set $V_{\Lambda} = \Lambda \otimes_{\bf Z} {\bf Q}$. Let $p$ be a prime
number with $p \not = \ell$. We have
$$\partial_p[V] = \partial_p[V_{\Lambda}] + \partial_p [V'].$$ Since $p \not = \ell$, by Proposition
\ref{root lattice lemma} we have $ \partial_p[V_{\Lambda}] = 0$, hence $\delta_p[V] = \delta_p[V']$,
as claimed. 

\medskip We need a similar result for the lattice $E_6$.

\begin{prop}\label{E6} {\rm (a)} Let $\Lambda$ be a root lattice of type $E_6$ having
an isometry with characteristic polynomial $\Phi_9$ {\rm (}resp. $\Phi_{18}${\rm )}. Set $V = \Lambda \otimes_{\bf Z} {\bf Q}$ and let $[V] \in W_{\Gamma}({\bf Q})$ be the  Witt class of the ${\bf Q}[\Gamma]$-form $V$. Then $\partial_p[V] = 0$ if $p$ is a prime number
with $p \not = 3$.

\medskip
{\rm (b)} Let $L$ be a negative definite lattice with an isometry $t : L \to L$,
$V = L \otimes_{\bf Z} {\bf Q}$. Let $\Lambda \subset
L$ be a lattice stabilized by $t$, assume that $\Lambda$ is isometric to an orthogonal sum of
 lattices $-E_6$, and that the
minimal polynomial of the restriction of $t$ to $\Lambda$ is  $\Phi_{9}$ or $\Phi_{18}$. Let $V'$ be the 
orthogonal of $\Lambda \otimes_{\bf Z} {\bf Q}$ in $V_L$. Then
$$\partial_p[V] = \partial_p[V']$$ for all prime numbers $p \not = 3$.

\end{prop}

\noindent
{\bf Proof.} Part (a) is an immediate consequence of the definition of the root lattice $E_6$, and the
proof of (b) is the same argument as the one of the proof of Proposition \ref{equal residue}.

\section{Unramified Salem polynomials of degree $10$}\label{continued section}

Let $S$ be an unramified Salem polynomial of degree $10$. The case of Salem signature
maps of Salem factor $S$ and maximum $(3,11)$ was treated in \S \ref{18 section}; we now consider those
of maximum $(3,19)$. 

\begin{theo}\label{10 part 1} Let $\tau$ be a Salem signature map with trivial obstruction of Salem factor $S$ and maximum $(3,19)$.
Then the lattice $\Lambda_{3,19}$ has a realizable isometry with signature
map equivalent to $\tau$.

\end{theo}

\noindent
{\bf Proof.} 
Let $F$ be the polynomial associated to $\tau$; we have
$G_F = 0$ by definition. Since $S$ is unramified, this implies that the polynomial $F$ is divisible by a cyclotomic polynomial
$\Phi_m$ with $m \geqslant 3$  such that $\Pi_{S,\Phi_m} \not = \varnothing$. Let $\Phi_m$ be such
a polynomial, and let $p$ be a prime number such that $p \in \Pi_{S,\Phi_m}$. 
Let $f \in {\bf F}_p[X]$ be
an irreducible, symmetric common factor of $S$ and $\Phi_m$ $({\rm mod} \ p)$, and let $M = {\bf F}_p[X]/(f)$, considered as a simple ${\bf F}_p[\Gamma]$-module. Note that
the polynomial $f$ has
even degree, and let $[M]$ denote the only non-trivial element of $W_{\Gamma}({\bf F}_p,M)$. 

\medskip
The case where $m = 3,4,6$ or $12$ was already treated in Theorem \ref{10 and 18}. In the remaining cases,
we construct a polynomial $C$ of degree $12$, product of cyclotomic polynomials and divisible
by $\Phi_m$, such that $SC$ satisfies conditions (C 0) and (C 1) and that $G_{SC} = 0$. We then
consider the Salem signature map $\tau_C$ of maximum $(3,19)$ and associated polynomial $SC$ such that $\tau_{C}(g) =
\tau(g)$ for all factors $g \in {\bf R}[X]$ of $S$. 
The signature map $\tau_C$ is equivalent to $\tau$. 
Let $t : L \to L$ be a semi-simple isometry with signature map $\tau_C$ of an even unimodular lattice $L$
such that that $\partial_v([V_S]) = \partial_v([V_C]) = 0$ for all $v \in \mathcal V'$ with $v \not = v_p$; these
exist by Proposition \ref{unique p}. Note that $\tau_C$ is equivalent to $\tau$. If $L_C$ has no roots, then
$t$ is realizable, and we are done. Assume that $L_C$ contains roots; we then construct a realizable
isometry with signature map equivalent to $\tau_C$. 
 This strategy is the same in all cases, but the polynomial $C$ and
the construction of the realizable isometry vary slightly from case to case.

\medskip
{\it Case} 1 : Suppose that $m = 5,7,8, 9,11, 14, 18$ or $22$. Set $C = \Phi_m^2$ if $m = 7,9,14$ or $18$, 
$C(X)  = \Phi_m^2(X) (X-1)^2$ if $m = 5$ or $8$,  $C(X) = \Phi_m(X)(X-1)^2$ if $m = 11$ and $C(X) = \Phi_m(X)(X+1)^2$ if $m = 22$.

\medskip
Suppose that $L_C$ contains a root lattice $L_0$ with minimal polynomial $\Phi_m$,
and let $L_1$ be the orthogonal of $L_0$ in $L_C$. Set $V_1 = L_1 \otimes_{\bf Z} {\bf Q}$. Then by Propositions \ref{equal residue} and \ref{E6} we have $\partial_{v_p}[V_C] = \partial_{v_p}[V_1]$; note
that if $m = 8$, then $\partial_v[V_C] = \partial_v[V_1]$ for all $v \in \mathcal V'$. If 
$m = 11$ or $m = 22$, then this leads to a contradiction, hence $t$ is realizable. Otherwise, let
$t' : L_C \to L_C$ be such that $t'|L_0 = id$ and $t'|L_1 = t|L_1$.  The isometry $t'$ is of type (P),
and determines a realizable isometry $t : L \to L$ with signature map equivalent to $\tau_C$.

\medskip
{\it Case} 2 : Suppose that $m = 10$, and set $C(X) = \Phi_{10}(X)(X-1)^6(X+1)^2$; let $L'$ be the rank 2
sublattice of $L_C$ given by $L' = {\rm Ker}(t+1)$, and note that the restriction of $t$ to $L'$ is $-id$.
By Lemma \ref{-id} the lattice $L'$ has an isometry $t'$ of type (P) such that $\partial_v[L_1 \otimes_{\bf Z} {\bf Q},t_1] = \partial_v[L_1 \otimes_{\bf Z} {\bf Q},-id]$. We now continue as in Case 1 : the isometry
$t'$ obtained with this method is of type (P). 

\medskip
{\it Case} 3 : Suppose that $m = 15, 16, 20, 24$ or $30$. Set $C = \Phi_m\Phi_3^2$ if $m = 15$ or $24$,
$C(X) = \Phi_m(X)(X-1)^2(X+1)^2$ if $m = 16$,
$C = \Phi_m \Phi_4^2$ if $m = 20$ and $C = \Phi_m \Phi_6^2$ if $m = 30$. Suppose that $L_C$ contains a root lattice $L_0$ isomorphic to $-E_8$ with minimal polynomial $\Phi_m$,
and let $L_1$ be the orthogonal of $L_0$ in $L_C$. Set $V_1 = L_1 \otimes_{\bf Z} {\bf Q}$. Since $E_8$
is unimodular, we have $\partial_{v}[V_C] = \partial_{v}[V_1]$ for all $v \in \mathcal V'$. This implies
that $\Pi_{S,\Phi_n} \not =  \varnothing $ with $n = 3, 4$ or $6$, hence we obtain a realizable isometry
by Theorem \ref{10 and 18}. If $m = 16$ and $L_C$ contains a root lattice isomorphic to $-D_8$, then 
$\partial_{v}[V_C] = \partial_{v}[V_1]$ for all $v \in \mathcal V'$, and this leads to a contradiction. 
Otherwise, $L_C$ contains a root lattice isomorphic to $A_{\ell-1}$ for $\ell = 3$ or $5$, and we conclude
as in Case~1.

\medskip
{\it Case} 4 : Assume that $m = 21, 28, 36$ or $42$, and set $C = \Phi_m$. Suppose that $L_C$ contains a root lattice; this 
lattice is an orthogonal sum of root lattices of type $A_{\ell -1}$ for $\ell = 3, 5$ or $7$, and this
implies that $\partial_{v_{\ell}}[V_C]$ is non-trivial. Hence $\partial_{v_{\ell}}[V_C]$ is also non-trivial,
and therefore $\Pi_{S,\Phi} \not = \varnothing$ for some cyclotomic polynomial $\Phi$ of degree
at most $6$; we conclude by applying Theorem \ref{10 and 18}, or Cases 1 or 2. 

\bigskip

\centerline {\bf Part III : Automorphisms of $K3$ surfaces}

\bigskip This part combines the arithmetic results of the first two parts to obtain applications to
automorphisms of $K3$ surfaces. 
We start by recalling some definitions and results
(see \S \ref{K3}); in particular, the  results of McMullen in \cite{Mc2}  needed in 
the following sections, and some consequences of fundamental results on $K3$ surfaces, such as the strong Torelli theorem and the
surjectivity of the period map. 

\medskip 
 The following sections concern {\it realization of Salem numbers}; we prove the
 results announced in the introduction, and \S \ref{finiteness} gives some finiteness results
 
 \medskip The last two sections contain some remarks on automorphisms
 of projective $K3$ surfaces : examples of Salem numbers that cannot be realized by automorphisms of projective $K3$ surfaces in \S \ref{non-realizable section}, and an
 alternative approach to some results of Kondo \cite {Ko} and Vorontsov \cite{V} in \S \ref{Kondo section}.

\section{$K3$ surfaces}\label{K3}

The aim of this section is to give some basic facts on $K3$ surfaces and their automorphisms
in a more general framework than needed in this paper; we refer to \cite {H}, \cite {K}, \cite {Ca 14}, \cite{Ca 99}, \cite{Ca 01}, \cite{O 7} and \cite{O 8}
for details. 

\medskip A $K3$ surface $\mathcal X$ is a simply-connected compact complex surface with trivial
canonical bundle. The dimension of the complex vector space $H^2(\mathcal X,{\bf C})$ is 22, and we
have the Hodge decomposition $$H^2(\mathcal X,{\bf C}) = H^{2,0}(\mathcal X) \oplus H^{1,1}(\mathcal X) \oplus H^{0,2}(\mathcal X)$$ with ${\rm dim}H^{2,0} = {\rm dim}H^{0,2} = 1$. The {\it Picard group} of 
$\mathcal X$ is by definition $${\rm Pic}(\mathcal X) = H^2(\mathcal X,{\bf Z}) \cap H^{1,1}(\mathcal X).$$
The intersection form $H^2(\mathcal X,{\bf Z}) \times H^2(\mathcal X,{\bf Z}) \to {\bf Z}$ of $\mathcal X$
is an even unimodular lattice of signature $(3,19)$; the signature of its restriction to ${\rm Pic}(\mathcal X)$ 
depends on the geometry of the surface. Let $a(\mathcal X)$ be the algebraic dimension of $\mathcal X$ (see for instance \cite {K}, 3.2). The following result is well-known  (see \cite{H}, Chapter 17, Proposition 1.3, or  \cite{K}, Proposition 4.11; see also \cite{N}, 3.2).

\begin{prop}\label{algebraic dimension} Let $\rho({\mathcal X})$ be the rank of ${\rm Pic}(\mathcal X)$. 

\medskip
{\rm (a)} If $a(\mathcal X) = 2$, then  ${\rm Pic}(\mathcal X)$ is non-degenerate and its signature is
$$(1,\rho({\mathcal X}) - 1).$$

\medskip
{\rm (b)} If $a(\mathcal X) = 1$, then  ${\rm Pic}(\mathcal X)$  has a one-dimensional kernel, and the
quotient by the kernel is negative definite.

\medskip
{\rm (c)} If $a(\mathcal X) = 0$, then ${\rm Pic}(\mathcal X)$ is negative definite.

\end{prop}

The global Torelli theorem implies that  $K3$ surfaces are classified in terms of
their Hodge structures and intersection pairing (see for instance \cite{H}, Theorem 2.4) :

\begin{theo}\label{torelli} Two $K3$ surfaces $\mathcal X$ and $\mathcal X'$ are isomorphic if and only if there
exists an isomorphism $H^2(\mathcal X,{\bf Z}) \to H^2(\mathcal X',{\bf Z})$
respecting the Hodge decomposition and the intersection form.

\end{theo}

An {\it automorphism} of a $K3$ surface $\mathcal X$ is a biholomorphic map $T : \mathcal X \to \mathcal X$;
such a map induces an isomorphism $T^* : H^2(\mathcal X,{\bf Z}) \to H^2(\mathcal X,{\bf Z})$ preserving
the intersection form, the Hodge decomposition, as well as the {\it K\"ahler cone} $C_{\mathcal X}$, that
is, the subset of $H^{1,1}(\mathcal X)$ consisting of the K\"ahler classes of $\mathcal X$ (see for
instance \cite{H}, Chapter 8, \S 5).  The Torelli theorem implies that the map sending $T$ to $T^*$ induces an isomorphism 
$${\rm Aut}({\mathcal X}) \to {\rm Aut}(H^2(\mathcal X),H^{2,0}(\mathcal X),C_{\mathcal X})$$ (see for
instance \cite{H}, Chapter 15, Theorem 2.3). 

\medskip Recall from the introduction that we denote by $\lambda(T)$ the dynamical degree of $T$, and
by $\omega(T)$ the complex number such that $T^* : H^{2,0}(\mathcal X) \to H^{2,0}(\mathcal X)$ is
multiplication by $\omega(T)$.

\medskip Following McMullen (cf. \cite {Mc2}, \S 6), we introduce the notion of {\it $K3$ structure of a lattice};
if an isometry preserves  this structure, then it is induced by an automorphism of a $K3$ surface.

\begin{defn} Let $(L,q)$ be an even unimodular lattice of signature $(3,19)$. A $K3$ structure on $L$ consists
of the following data :

\medskip
(1) A {\it Hodge decomposition} $$L \otimes_{\bf Z} {\bf C} = L^{2,0} \oplus L^{1,1} \oplus L^{0,2}$$ of complex
vector spaces such that $L^{i,j} = \overline {L^{j,i}}$, and that the hermitian spaces $L^{1,1}$ and
$L^{2,0} \oplus L^{0,2}$ have signatures $(1,19)$ and $(2,0)$, respectively.

\medskip
(2) A {\it positive cone} $\mathcal C \subset L^{1,1} \cap (L \otimes_{\bf Z} {\bf R})$, consisting of
one of the two components of the locus $q(x,x)  > 0$.

\medskip
(3) A set of {\it positive roots} $\Psi^+ \subset \Psi = \{x \in L \cap L^{1,1} \ | \ q(x,x)= -2 \}$ such that
$\Psi = \Psi^+ \cup (-\Psi^+)$, and that the set (called the {\it K\"ahler cone}) 

\medskip
\centerline {$C(L) = \{x \in \mathcal C \ | \ q(x,y) > 0 $ for all $y \in \Psi^+ \}$ }

\medskip
\noindent is not empty. 

\bigskip

We say that a $K3$ structure  is {\it realized} by a $K3$ surface $\mathcal X$ if there exists
an isomorphism $\iota : L \to H^2(\mathcal X,{\bf Z})$ sending $L^{i,j}$ to $H^{i,j}(\mathcal X)$,
and sending $C(L)$ to the K\"ahler cone $C_{\mathcal X}$; an isometry $t : L \to L$
is {\it realized} by an automorphism $T : \mathcal X \to \mathcal X$ if $\iota$ can be chosen so
that $\iota \circ T^* = t \circ \iota$. 

\end{defn}

Note that if  a $K3$ structure on $L$ is realized by $\mathcal X$, then ${\rm Pic}(\mathcal X) = \iota(L \cap L^{1,1})$. 
The strong Torelli theorem and the surjectivity of the period map imply the following (see
\cite {Mc2}, Theorem 6.1) :

\begin{theo} Let $L$ be an even unimodular lattice of signature $(3,19)$. Then every $K3$ structure on $L$ is realized by a unique $K3$ surface $\mathcal X$,
and every isometry of $L$ preserving the $K3$ structure is realized by a unique
automorphism of  $\mathcal X$. 

\end{theo}

\medskip
Recall that the {\it transcendental lattice} of a $K3$ surface $\mathcal X$ is the primitive sublattice
${\rm Trans}(\mathcal X)$ of minimal rank of $H^2(\mathcal X,{\bf Z})$ such that the vector space
${\rm Trans}(\mathcal X)\otimes_{\bf Z}{\bf C}$ contains $H^{2,0}(\mathcal X) \oplus H^{0,2}(\mathcal X)$. 
The sublattices ${\rm Pic}(\mathcal X)$ and ${\rm Trans}(\mathcal X)$ of $H^2(\mathcal X,{\bf Z})$ are
orthogonal to each other; moreover, if
$a(\mathcal X) = 0$ or $2$, then ${\rm Pic}(\mathcal X) \oplus {\rm Trans}(\mathcal X)$ is a sublattice
 finite index in $H^2(\mathcal X,{\bf Z})$. 

\medskip
The following result is due to Oguiso, \cite{O 8}, Theorem 2.4 (1) :
 
 \begin{lemma}\label{transcendental} Let $T : \mathcal X \to \mathcal X$ be an automorphism of
 a $K3$ surface. The minimal polynomial of the restriction of $T^*$ to ${\rm Trans}(\mathcal X)$ is irreducible. 
 
 \end{lemma}
 
 \noindent
 {\bf Proof.} This is \cite{O 8}, Theorem 2.4 (1); we give a proof for the convenience of the reader.  Let $f \in {\bf Z}[X]$ be an irreducible factor of the minimal polynomial of the restriction of $T^*$ to ${\rm Trans}(\mathcal X)$, and suppose that $\omega(T)$ is
a root of $f$. Let $N = {\rm Ker}(f(T^*))$; then $N \otimes_{\bf Z} {\bf C}$ contains $H^{2,0}(X) \oplus H^{0,2}(X)$. By the minimality property of ${\rm Trans}(\mathcal X)$,
this implies that $N = {\rm Trans}(\mathcal X)$.

\medskip

The existence of an automorphism of dynamical degree $>1$ imposes some restrictions on the $K3$ surface; in
particular, its algebraic dimension is $0$ or $2$ (see Corollary \ref {0 or 2}).

\begin{prop}\label{LS} Let $T : \mathcal X \to \mathcal X$ be an automorphism with $\lambda(\mathcal X) > 1$, and
let $F$ be the characteristic polynomial of $T$. Then we have

\medskip
{\rm (a)} $F = SC$, where $S$ is a Salem polynomial, and $C$ is a product 
of cyclotomic polynomials. 



\medskip
{\rm (b)} Suppose that $\omega(T)$ is a root of $S$. Then the signature of 
${\rm Pic}(\mathcal X)$ is $(0,\rho({\mathcal X}))$, and ${\rm Pic}(\mathcal X) = {\rm Ker}(C(T^*))$; the signature of 
${\rm Trans}(\mathcal X)$ is $(3,19-\rho({\mathcal X}))$ and ${\rm Trans}(\mathcal X) =  {\rm Ker}(S(T^*))$.

\medskip
{\rm (c)} If $\omega(T)$ is a root of unity, then the signature of ${\rm Pic}(\mathcal X)$ is $(1,\rho({\mathcal X})-1)$ and
the signature of ${\rm Trans}(\mathcal X)$ is $(2,20-\rho({\mathcal X}))$.

\medskip
{\rm (d)} $T^*$ is semi-simple.

\end{prop}

\noindent
{\bf Proof.} (a) Recall that by  \cite{Mc1}, Corollary 3.3, the polynomial $F$ is a product of at most one Salem polynomial and of cyclotomic polynomials.
Since we are assuming that $\lambda(T) > 1$, this implies that $F$ is divisible by a Salem polynomial.

\medskip Let $L = H^2(\mathcal X,{\bf Z})$, $L_S = {\rm Ker}(S(T^*))$ and $L_C = {\rm Ker}(C(T^*))$. 
The
lattices $L_S$ and $L_C$ are orthogonal to each other, and $L_S \oplus L_C$ is of finite index in $L$. 
Set $d = {\rm deg}(S)$. 



\medskip (b) $L_S \otimes_{\bf Z} {\bf C}$ contains $H^{2,0} \oplus H^{0,2}$, since $\omega(T)$ is a root of $S$, therefore the signature of $L_S$ is $(3,d-3)$; 
we have ${\rm Trans}(\mathcal X)= L_S$ by the minimality of the lattice ${\rm Trans}(\mathcal X)$. This implies that  ${\rm Pic}(\mathcal X) = L_C$, and
and the signature of ${\rm Pic}(\mathcal X)$ is $(0,\rho({\mathcal X}))$.

\medskip (c) Since $\omega(T)$ is a root of unity, the lattice $L_S$ is orthogonal to $H^{2,0} \oplus H^{0,2}$, and
therefore $L_S$ is a sublattice of ${\rm Pic}(\mathcal X)$. This implies that the signature of ${\rm Pic}(\mathcal X)$ is $(1,\rho({\mathcal X}) - 1)$,
and  the signature of ${\rm Trans}(\mathcal X)$ is $(2,20-\rho({\mathcal X}))$.

\medskip (d) If $\omega(T)$ is a root of $S$, the lattice $L_C$ is negative definite, therefore the
restriction of $T^*$ to $L_C$ is semi-simple; the restriction of $T^*$ to $L_S$ is also semi-simple, because
$S$ is irreducible. Suppose that $\omega(T)$ is a root of unity. Then $S$ divides the minimal polynomial of $T^*|{\rm Pic}(\mathcal X)$, and the signature of ${\rm Pic}(\mathcal X)$ is $(1,\rho({\mathcal X}) - 1)$; this implies that the restriction of $T^*$ to ${\rm Pic}(\mathcal X)$ is semi-simple. The restriction of $T^*$ to ${\rm Trans}(\mathcal X)$ is semi-simple
by Lemma \ref{transcendental}. 

\begin{coro}\label{0 or 2} If $\mathcal X$ has an automorphism $T$ with $\lambda(T) > 1$, then $a(\mathcal X) = 0$ or
 $a(\mathcal X) = 2$. 

\end{coro}

\noindent {\bf Proof.} This follows from Theorem \ref{algebraic dimension}, combined with Proposition \ref{LS} (b) and (c).  

\medskip





\section{Realizable Salem numbers}\label{Mac}

We now return to the main topic of this paper; following McMullen
(\cite{Mc1}, \cite{Mc2}, \cite{Mc3}), we investigate the following question :

\begin{question}
{\it Which Salem numbers occur as dynamical degrees of automorphisms of complex $K3$ surfaces ?}
\end{question}

More precisely, {\it which pairs $(\lambda(T),\omega(T))$ occur} ? This leads to the definition~:

\begin{defn}\label{realizable definition} Let $\alpha$ be a Salem number, and let
$\delta \in {\bf C}^{\times}$ be such that $|\delta| = 1$. 
We say that the pair $(\alpha,\delta)$ is {\it realizable} if 
there exists an automorphism $T$ of a $K3$ surface with $\lambda(T) = \alpha$ and
$\omega(T) = \delta$.

\end{defn}

The following implies Theorem 3 of the introduction :

\begin{theo}\label{congr 4}  Let $S$ be Salem polynomial of degree $d$ with $4 \leqslant d \leqslant 18$, and let
$\delta \in {\bf C}^{\times}$ be such that $|\delta| = 1$. Assume that one of the following holds

\medskip

{\rm (a)}  $d  \equiv 0, 4 \ {\rm or} \  6 \ {\rm (mod \ 8)}$.  

\medskip

{\rm (b)} $d  \equiv 2 \ {\rm (mod \ 8)}$ and $|S(1)|$, $S(-1)$ are not both squares.

\medskip Then {$(\alpha,\delta)$ is realizable.}

\end{theo}

The proof uses the results of Parts I and II, as well as a criterion of McMullen; the 
following result is a reformulation of this criterion in the case we need here. See Definition
\ref{realizable isometry} for the notion of a realizable isometry.

\begin{theo}\label{Mc6.2}  Let $L$ be an even unimodular lattice of signature $(3,19)$, and let $t : L \to L$ be
a realizable isometry with signature map equivalent to $\tau_{S,\delta}$ and characteristic polynomial $S C$, where $C$ is
a product of a finite number of cyclotomic polynomials such that  ${\rm deg}(C) = 22-d$. Set $L_C = {\rm Ker}(C(t))$. 

\medskip
Then there exists a
$K3$ surface automorphism $T : \mathcal X \to \mathcal X$ with $t = T^*$. Moreover~:

\medskip {\rm (a)} We have $\lambda(T) = \alpha$ and $\omega(T) = \delta$. 

\medskip 
{\rm (b)} 
The Picard lattice of the $K3$ surface $\mathcal X$
is isomorphic to $L_C$. 

\end{theo}

\noindent
{\bf Proof.} Let us check that the hypotheses of \cite{Mc2}, Theorem 6.2 are fulfilled. We have 
$\alpha > 1$ by definition, and $\alpha$ is an eigenvalue of $t$, since $S$
divides the characteristic polynomial of $t$. Let $\nu = \delta+ \delta^{-1}$; then $E_{\nu} = {\rm Ker}(t + t^{-1} - \nu)$ has
signature $(2,0)$. Set $M = L \cap E_{\nu}^{\perp}$. Then $M$ is negative definite, and
$M = {\rm Ker}(C(t))$. Since $t$ is a realizable isometry, the restriction of $t$ to $M$ is of type (P). With
the terminology of \cite{Mc2}, this implies that the restriction of $t$ to  $M(-1)$
is positive. Therefore we can apply Theorem 6.2 of \cite{Mc2}, and conclude that
$t$ is realizable by a $K3$ surface automorphism $T : \mathcal X \to \mathcal X$. We
have $\lambda(T) = \alpha$ and $\omega(T) = \delta$ by construction, and ${\rm Pic}(\mathcal X)$
is isomorphic to $M = {\rm Ker}(C(t))$ by the above remark, hence (a) and (b) hold. This
completes the proof of the theorem. 

\medskip
With the terminology of \S \ref{Salem signature section}, this implies the following :

\begin{coro}\label{Mc cor} If the lattice $\Lambda_{3,19}$ has a realizable isometry with signature map
equivalent to $\tau_{S,\delta}$, then $(\alpha,\delta)$ is realizable. 

\end{coro}

\noindent
{\bf Proof of Theorem \ref{congr 4}.} Part (a) follows from Theorem \ref{realizable theorem} and Corollary \ref{Mc cor}, and
Part (b) from Proposition \ref{not both squares} combined with Corollary \ref{Mc cor}.

\begin{example}\label{16} Let $\alpha = \lambda_{16} = 1,2363179318...$ be the smallest 
degree $16$ Salem number (see for instance the table in Boyd \cite{Bo}). The corresponding Salem polynomial is
$$S(X) = X^{16} - X^{15} - X^8 - X + 1.$$ Theorem \ref{congr 4} (a)  implies that $(\alpha,\delta)$ is realizable
for every root $\delta$ of $S$ with $|\delta| = 1$. On the other hand, McMullen proved that $\alpha$ 
is {\it not} realizable by an automorphism of a projective $K3$ surface (see \cite{Mc3}, \S 9). 

\end{example}

\begin{example} Let $\alpha = 1,2527759374...$ be the Salem number with
polynomial $$S(X) = X^{18} - X^{12} - X^{11} - X^{10} - X^9 - X^8 - X^7 -X^6 + 1.$$ We have
$S(1) = -5$, hence Theorem \ref{congr 4} (b)  implies that $(\alpha,\delta)$ is
realizable for all roots $\delta$ of $S$ with $|\delta| = 1$. 

\end{example}

 We now prove Theorem 1 of the introduction :

\begin{theo}\label{theorem 1}
 Let $S$ be a Salem polynomial of degree $d$ with $4 \leqslant d \leqslant 22$, and let $\alpha$ be the corresponding Salem number. The following
are equivalent :

\medskip
{\rm (1)} There exists a root $\delta$ of $S$ with $|\delta| = 1$ such that $(\alpha,\delta)$ is realizable.

\medskip
{\rm (2)}  $(\alpha,\delta)$ is realizable for all roots $\delta$ of $S$ with $|\delta| = 1$.

\end{theo}

\noindent
{\bf Proof.} It is clear that (2) $\implies$ (1);  let us prove that (1) $\implies$ (2). 
Since $(\alpha,\delta)$ is realizable, there exists an even unimodular lattice $L$ of signature $(3,19)$ and a realizable isometry $t : L \to L$ with
signature map equivalent to 
 $\tau_{S,\delta}$; let us denote by $\tau$ this signature map. Let $F$ be the characteristic polynomial of $t$; we have $F = SC$, where $C$ is a product of cyclotomic polynomials. Set
$L_S = {\rm Ker}(S(t))$, $L_C = {\rm Ker}(C(t))$, and let $V_S = L_S \otimes_{\bf Z}{\bf Q}$, $V_C = L_C \otimes_{\bf Z}{\bf Q}$; both $V_S$ and $V_C$ are
stabilized by $t$. Let $V = L \otimes_{\bf Z}{\bf Q}$,
and let $q : V \times V \to {\bf Q}$ be the quadratic form given by the lattice $L$; we have $V = V_S \oplus V_C$, and this decomposition is orthogonal with respect
to the form $q$. Let $[V_S], [V_C] \in W_{\Gamma}({\bf Q})$ be the corresponding Witt classes. 
For all prime numbers $p$, let
$\delta^p_S = \partial_p[V_S]$, $\delta^p_C = \partial_p[V_C]$ be their images in $W_{\Gamma}({\bf F}_p)$; since $L$ is unimodular, we have $\delta^p_S \oplus \delta^p_C = 0$.

\medskip
For all prime numbers $p$, let  $a^p(S) \in \mathcal C_{\delta^p_S }$ and $a^{\infty}_{\tau_{S,\delta}}(S)$ be as in \S \ref{local data}, Notation \ref{a 1} and  \S \ref{Hasse section}. 
Let $\delta'$ be a root of $S$ with  $|\delta| = 1$. By Proposition \ref{dependance on tau}, we have
$a^{\infty}_{\tau_{S,\delta}}(S) = a^{\infty}_{\tau_{S,\delta'}}(S)$.  Theorem \ref{realization} implies that there exists a semi-simple isometry $t' : V_S \to V_S$ with signature map
$\tau_{S,\delta'}$ and local data $(a^p_S)$; let $[V'_S] \in W_{\Gamma}({\bf Q})$ be the corresponding Witt class.  
For all prime numbers $p$, let $N^p$ be a maximal ${\bf Z}_p$-lattice in $V_S \otimes_{\bf Q} {\bf Q}_p$ such that the class of $((N^p)^{\sharp}/N^p,t')$
in $W_{\Gamma}({\bf F}_p)$ is equal to $\partial_p[V'_S]$; note that for almost all $p$, we have $N^p = L_S \otimes_{\bf Z}{\bf Z}_p$. Let $N$ be the intersection
of the lattices $N^p$ in $V_S$; then $N$ is stabilized by the isometry $t'$. 
Let $t' : V \to V$ be such that the restriction of $t'$ to $V_S$ is $t'$, and the restriction
of $t'$ to $V_C$ is $t$. Since $\partial_p[V'_S] = \delta^p_S$ for all prime numbers $p$, the isometry $t' : V \to V$ stabilizes an even unimodular lattice $L'$ such
that ${\rm Ker}(C(t')) = L_C$. 
The isometry $t'$ is realizable, since the restriction of $t'$ to $V_C$ is equal to $t$. The signature
map of $t'$ is equivalent to $\tau_{S,\delta'}$, hence by Corollary \ref{Mc cor} $(\alpha,\delta')$ is realizable. 

\medskip This concludes the proof of the theorem.

\bigskip

The following is Theorem 2 of the introduction :

\begin{theo}\label{22} Suppose that $d = 22$ and that $|S(1)|$ and $S(-1)$ are squares. Then $(\alpha,\delta)$ is realizable.

\end{theo} 

\noindent
{\bf Proof.} This follows from Theorem \ref{22 theo} and Corollary \ref{Mc cor}.

\medskip

In the case of unramified Salem numbers of degree $10$ or $18$, we obtain
a necessary and sufficient criterion in terms of Salem signature maps.

\begin{theo}\label{equivalent conditions} Assume that $d = 10$ or $18$ and that $S$ is unramified.
 Then $(\alpha,\delta)$ is
realizable if and only if there exists a Salem signature map of maximum $(3,19)$ with trivial obstruction equivalent to $\tau_{S,\delta}$.

\end{theo}

\noindent
{\bf Proof.} Let $\tau$ be a Salem signature map of maximum $(3,19)$ with trivial obstruction equivalent to $\tau_{S,\delta}$; then by Theorem 
\ref{10 and 18} the lattice $\Lambda_{3,19}$ has a realizable isometry with signature
map equivalent to $\tau$. By \ref{Mc cor}, we conclude that $(\alpha,\delta)$ is realizable. 

\medskip
Conversely, suppose that $(\alpha,\delta)$ is realizable, and let $t : L \to L$ be a
semi-simple isometry $t$ of an even unimodular lattice $L$ with signature map $\tau$ equivalent to $\tau_{S,\delta}$.  Let $F$ be the
polynomial associated to $\tau$, i.e. the characteristic polynomial of $t$; we have $F = S C$, where
$C$ is a product of cyclotomic polynomials. The polynomial $F$ satisfies condition (C 1). Set $L_S = {\rm Ker}(S(t))$, $L_C = {\rm Ker}(C(t))$ and
let $V_S = L_S \otimes_{\bf Z} {\bf Q}$ and $V_C = L_C \otimes_{\bf Z} {\bf Q}$. We are assuming that $S$ is unramified, 
hence  $S$ satisfies condition (C 1); this implies that  $C$ also satisfies condition (C 1), and therefore we have
$\partial_v[V_S] = \partial_v[V_C] = 0$ for all $v \in \mathcal V'$. On the other hand, since $\delta$
is a root of $S$, the signature of $L_S$ is $(3,d-3)$ and $L_C$ is negative definite; hence $L_S$ and
$L_C$ cannot be unimodular. This implies that $C$ has an irreducible factor $\Phi$ such that
$\Pi_{S,\Phi} \not = \varnothing$. The hypothesis $S$ unramified  implies that $\Phi(X) \not = X-1$ or
$X+1$, hence $\Phi$ is a cyclotomic polynomial of even degree. There exists a product of cyclotomic polynomials $C'$ of degree $22-d$  divisible by $\Phi$ satisfying conditions (C 0) and (C 1)
such that $G_{SC'} = 0$. Indeed, if 
${\rm deg}(\Phi) = 22-d$, then by
Condition (C 1) we have $\Phi(1) = \Phi(-1) = 1$; in this case, set $C' = \Phi$. Otherwise, ${\rm deg}(\Phi) \leqslant 20-d$.
Using the properties of cyclotomic polynomials and the fact that $22-d$ is divisible by $4$, 
we check that such a polynomial $C'$ exists. 
Let $\tau$ be the signature map of maximum $(3,19)$ with polynomial $SC'$
equivalent to $\tau_{S,\delta}$; this is a Salem signature map, hence the theorem is proved.

\medskip
Recall that if $S$ is unramified and if $m$ is an integer $\geqslant 3$, then $\Pi_{S,\Phi_m}$ is the set of prime numbers $p$ such that 
 {$S  \ {\rm (mod \ {\it p})}$ and $\Phi_m  \ {\rm (mod \ {\it p})}$}
  have a common irreducible 
symmetric factor of even degree in ${\bf F}_p[X]$.

\begin{theo}\label {18 theorem} Assume that $d = 10$ or $18$ and that $S$ is unramified. If we have $\Pi_{S,\Phi_m} \not = \varnothing$ for some
$m = 3,4, 6$ or $12$, then 

\medskip {\rm (a)}  $(\alpha,\delta)$ is realizable. 

\medskip {\rm (b)} If moreover $d = 10$, then 
the lattice $\Lambda_{3,19}$ has a realizable isometry  with signature map equivalent to $\tau_{S,\delta}$ and polynomial divisible by $(X-1)^8$.

\end{theo}

\noindent{\bf Proof.} (a)
Set $C = \Phi_m^2$ if $m = 3,4$ or $6$, and  $C = \Phi_m$ if $m = 12$. Then the polynomial $SC$ 
satisfies conditions (C 0) and (C 1), and $\Pi_{S,\Phi_m} \not = \varnothing$  implies
that $G_{SC}= 0$.  Let $\tau$ be the signature map of maximum $(3,d+4)$ and associated polynomial
$SC$ be such that 
$\tau(f) = \tau_{S,\delta}(f)$ for all factors $f \in {\bf R}[X]$ of $S$, and note that
$\tau$ is a Salem signature map; hence Theorem \ref{equivalent conditions} implies that  $(\alpha,\delta)$
is realizable. Part (b) follows from Theorem \ref{10 and 18}.

\section{Salem numbers of degree 18}

We keep the notation of the previous section. 
Suppose that $d = 18$ and $S$ is unramified (for instance, 
$|S(1)S(-1)| = 1$). In this case, the necessary and sufficient condition of Theorem \ref{equivalent conditions}
can be reformulated as follows :

\begin{theo}\label {18 proposition} Suppose that $S$ is unramified. Then $(\alpha,\delta)$ is realizable if and only if $\Pi_{S,\Phi_m} \not = \varnothing$ for
some $m = 3,4,6$ or $12$.

\end{theo}

\noindent{\bf Proof.} We already know that if $\Pi_{S,\Phi_m} \not = \varnothing$ for $m = 3,4,6$ or $12$,
then $(\alpha,\delta)$ is realizable (cf. Theorem \ref{18 theorem}). The converse is a consequence of
Theorem \ref{equivalent conditions}, combined with some properties of cyclotomic polynomials. 

\begin{example} Let  $\alpha= 1,21972085590$ be the second smallest degree $18$ Salem
number, with Salem polynomial $$S(X) = X^{18} - X^{17} - X^{10} + X^9 - X^8 - X + 1$$
(see the table in \cite{Bo}). We have $\Pi_{S,\Phi_3} = \{5 \}$, therefore 
by Theorem \ref{18 proposition} the pair $(\alpha,\delta)$ is realizable for every root $\delta$ of $S$ 
with $|\delta| = 1$.

\end{example}

\medskip

On the other hand, we now show that there exist nonrealizable pairs $(\alpha,\delta)$.
If $f \in {\bf Z}[X]$ is a monic polynomial, we denote by ${\rm Res}(S,f)$ the resultant of the polynomials
$S$ and $f$. 

\begin{prop}\label{resultant}  Assume that $|S(1)S(-1)| = 1$, and that ${\rm Res}(S,\Phi_m) = 1$ for $m = 3,4$ and $6$. If 
$(\alpha,\delta)$ is realizable, then 
$\Pi_{S,\Phi_{12}} \not = \varnothing$. 

\end{prop}

\noindent
{\bf Proof.} Suppose that $(\alpha,\delta)$ is realizable. This implies that $\Lambda_{3,19}$ has
a semi-simple isometry with signature map $\tau$ of maximum $(3,19)$ and polynomial $SC$,
where $C$ is a product of cyclotomic polynomials. 
Assume that all the factors of $C$ belong to the set $\{\Phi_1,\Phi_2,\Phi_3,\Phi_4, \Phi_6 \}$. Then $S$ and $C$
are relatively prime over ${\bf Z}$. If $\Lambda_{3,19}$ has an isometry with characteristic polynomial $F$, then
$\Lambda_{3,19} = L_1 \oplus L_2$, where $L_1$ and $L_2$ are even unimodular lattices such that
$L_1$ has an isometry with characteristic polynomial $S$ and signature map  $\tau_{S,\delta}$, and $L_2$ has
an isometry with characteristic polynomial $C$. This implies that the signature of $L_1$ is $(3,15)$ and that
the signature of $L_2$ is $(0,4)$, and this is impossible. 
Therefore the only possiblity is $C = \Phi_{12}$.

\medskip Let us show that if $\Pi_{S,\Phi_{12}} = \varnothing$, then $(\alpha,\delta)$ is not realizable.
Note that $\Pi_{S,\Phi_{12}} = \varnothing$ implies that $G_{S \Phi_{12}} \not = 0$. We apply Example
\ref{Salem example}, (c); since $S$ and $\Phi_{12}$ both satisfy condition (C 1), the group
$G_{S \Phi_{12}}$ is $\not = 0$,  and $\delta$ is a root of $S$, the signature map $\tau$ is
not realizable, and this implies that $(\alpha,\delta)$ is not realizable.

\begin{example}\label{second smallest}
Let $\lambda_{18} = 1.1883681475...$, the smallest degree $18$ Salem number (it is also the second smallest known Salem number). Let
$S$ be the corresponding Salem polynomial, and let $\delta$ be a root of $S$ with
$|\delta| = 1$; then $(\lambda_{18},\delta)$ is {\it not
realizable}. Indeed, the
polynomial $S$ satisfies the conditions of Proposition 
\ref{resultant} : we have ${\rm Res}(S,f) = 1$ for all $f \in \{ \Phi_3,\Phi_4, \Phi_6 \}$. Therefore
by Proposition \ref{resultant}, if $(\lambda_{18},\delta)$ is realizable, then 
$\Pi_{S,\Phi_{12}} \not = \varnothing$.

\medskip
We have ${\rm Res}(S,\Phi_{12})= 169$, and the common factors modulo $13$
of $S$ and $\Phi_{12}$ in ${\bf F}_{13}[X]$ are $X + 6, X + 11 \in {\bf F}_{13}[X]$. These polynomials are not symmetric.
Therefore $\Pi_{S,\Phi_{12}} = \varnothing$; hence 
Proposition \ref{resultant} implies that $(\lambda_{18},\delta)$ is not realizable. 

\medskip
On the other hand, McMullen
proved that $\lambda_{18}$ is realized by an automorphism of a projective $K3$ surface
(see \cite{Mc3}, Theorem 8.1).

\end{example}

\begin{example} With the notation of Example \ref{second smallest}, set $\alpha_2 = \lambda_{18}^2$ and $\alpha_3 = \lambda_{18}^3$; let $S_2$ and $S_3$ be
the associated Salem polynomials. The Salem numbers $\alpha_2$ and $\alpha_3$ are realizable : indeed, we have $\Pi_{S_2,\Phi_6} = \{13 \}$ and $\Pi_{S_3,\Phi_3} = \{17 \}$.

\end{example} 

\section{Salem numbers of degree 10}\label{10 section}

We keep the notation of the previous sections : $S$ is a Salem polynomial of degree $d$ with Salem
number $\alpha$, and $\delta$ is a root of $S$ such that $|\delta| = 1$. 
Suppose that $d = 10$, and that $S$ is {\it unramified};  for instance, this is the case when
 $|S(1)S(-1)| = 1$.
Recall that if $m \geqslant 3$ is an integer, then $\Pi_{S,\Phi_m}$ is the set of prime numbers $p$ such that 
 {$S  \ {\rm (mod \ {\it p})}$ and $\Phi_m  \ {\rm (mod \ {\it p})}$}
  have a common irreducible 
symmetric factor of even degree in ${\bf F}_p[X]$.

\begin{theo}\label {10.2 theorem} Assume that there exists an integer $m \geqslant 3$ with
$\varphi(m) \leqslant 12$ and $m \not = 13$ or $26$ such that $\Pi_{S,\Phi_m} \not = \varnothing$. Then
$(\alpha,\delta)$ is realizable. 

\end{theo}

\noindent
{\bf Proof.} The hypothesis implies that there exists a product of
cyclotomic polynomials $C$ of degree $12$ such that the polynomial $F = S C$ satisfies conditions (C 0) and
(C~1), and $G_F = 0$. Let $\tau$ be the signature map of maximum $(3,19)$ and associated polynomial
$SC$ be such that 
$\tau(f) = \tau_{S,\delta}(f)$ for all factors $f \in {\bf R}[X]$ of $S$, and note that
$\tau$ is a Salem signature map; hence Theorem \ref{equivalent conditions} implies that  $(\alpha,\delta)$
is realizable. 

\medskip
Examples \ref{GM 10} and \ref{b} give infinite families of realizable degree $10$ Salem numbers~:
 
 \begin{example}\label{GM 10} Let $a \geqslant 1$ be an integer, let
 
 \medskip
 $R(X) = (X+1)^2(X^2-4)(X-a) - 1,$ and $S(X) = X^5 R(X + X^{-1}).$

 \medskip
 Then  $S$ is a Salem polynomial (see \cite{GM}, \S 7); let $\alpha_a$ be the corresponding Salem number. We have $S(i) = R(0) = 4a-1$. Since $a \geqslant 1$, this is congruent to $3$ modulo $4$, hence
 it is divisible by a prime number $p  \equiv 3 \ {\rm (mod \ 4)}$.
 Note that the resultant of $S$ and $\Phi_4$ is ${\rm N}_{{\bf Q}(i)/{\bf Q}}(S(i))$,  and that $\Phi_4$ is irreducible mod $p$. Therefore we have
$p \in \Pi_{S,\Phi_4}$, and  Theorem \ref{10.2 theorem} implies that $(\alpha,\delta)$ 
 is realizable for all roots $\delta$ of $S$ with $|\delta| = 1$. 
 
 \end{example}

\begin{example}\label{b}  Let $a,b, c \in {\bf Z}$, set $$R(X) = (X^2-4)(X^3 +aX^2 +(b-1)X +c) - 1,$$ and $$S(X) = X^5 R(X + X^{-1}).$$ 

\medskip
Suppose that $c \geqslant 0$ and that $a + c < -|b|$. Then we have $R(0) < 0$ and
$R(-1), R(1) > 0$. Since $R(-2), R(2) < 0$, the polynomial $R$ has $4$ roots between $-2$ and $2$. Moreover,
$R$ is irreducible (otherwise, one of the factors of $R$ would be a ``cyclotomic trace polynomial"). Therefore $S$ is a Salem polynomial. 

\medskip
We have $S(\zeta_3) = R(-1) = -3(a - b + c) - 1$; since  $a + c < -|b|$ by hypothesis, this number is
divisible by a prime number $p$ congruent to $2$ modulo $3$, and this implies that 
$p \in \Pi_{S,\Phi_3}$. Let $\alpha$ be the Salem number corresponding to $S$. By 
Theorem \ref{10.2 theorem}, the pair $(\alpha,\delta)$ 
 is realizable for all roots $\delta$ of $S$ with $|\delta| = 1$. 

\end{example}

\begin{example}\label{smallest Salem} Let $\alpha = \lambda_{10}$ be the smallest known Salem number (also called ``Lehmer number"), and let $S$ be
the corresponding Salem polynomial. The polynomial $S$ belongs to the family of Example \ref{a}, for $a = 1$; we have $3 \in \Pi_{S,\Phi_4}$, 
hence $(\alpha,\delta)$ is realizable for all roots $\delta$ of $S$ with $|\delta| = 1$. An explicit construction (for $\delta + \overline {\delta} = - 1.886...$) is given by McMullen in \cite{Mc2}, Theorem 7.1. 
 
 \medskip
 We also have $\Pi_{S,\Phi_{12}} = \{3\}$, $\Pi_{S,\Phi_{14}} = \{13\}$, $\Pi_{S,\Phi_{15}} = \{29\}$ and $\Pi_{S,\Phi_{36}} = \{3\}$, also leading to
 realizations of $(\alpha,\delta)$ for all roots $\delta$ of $S$ with $|\delta| = 1$. 

\medskip
The Lehmer number is also the dynamical degree of automorphisms of projective $K3$ surfaces, as shown by McMullen (see \cite{Mc3}, Theorems 7.1 and 7.2.)

\end{example}

These results and examples suggest that all degree $10$ Salem numbers might be realizable, suggesting the question :

\begin{question}\label{degree 10 question} Is every Salem number of degree $10$ realizable ?

\end{question}

In particular, I do not know the answer to the
following question :

\begin{question} Let $S$ be a Salem polynomial of degree $10$ with $|S(1)S(-1)| = 1$. Does there
exist an integer $m \geqslant 3$ with
$\varphi(m) \leqslant 12$ and $m \not = 13$ or $26$ such that $\Pi_{S,\Phi_m} \not = \varnothing$ ?

\end{question}

\section{Salem numbers of degree 20}\label{20 section}

The aim of this section is to give some partial results and examples concerning Salem numbers of degree $20$. More recently,
Takada proved that {\it every Salem number of degree $20$ is realizable} (see \cite{T}, Theorem 1.3). 

\medskip
Let $\alpha$ be a Salem number of degree $20$,
let $S$ be the corresponding Salem polynomial, and let $\delta$
be a root of $S$ with $|\delta| = 1$. Let $F(X) = S(X)(X^2 - 1)$, and let $\tau_{\delta}$, $\tau_{X^2-1}$ be
the signature maps of maximum $(3,19)$ and polynomial $F$ such 
that $\tau_{\delta}$ is equivalent to $\tau_{S,\delta}$, and that $\tau_{X^2-1}(X-1) = \tau_{X^2-1}(X+1) = 1$.

\begin{theo}\label{relatively prime} Suppose that $S(1)$ and $S(-1)$ are odd, square-free,  and relatively prime. Then 

\medskip
{\rm (a)} The lattice $\Lambda_{3,19}$ has a realizable isometry with signature map $\tau_{\delta}$.

\medskip
{\rm (b)} The lattice $\Lambda_{3,19}$ has an isometry with signature map $\tau_{X^2-1}$. 

\end{theo}

\noindent
{\bf Proof.} We prove (a) and (b) simultaneously. 
Set $D =S(1)S(-1)$, and let $K = {\bf Q}(\sqrt D)$; since $S(-1) < 0$ and $S(1) > 0$, the
field $K$ is an imaginary quadratic field. For all prime numbers $p$ dividing $D$, we denote by $P$
the ramified ideal of $K$ above $p$; let $I$ be the product of the ideals $P$ for $p$ dividing $S(1)$.
Let $q_I : K \times K \to {\bf Q}$ be the quadratic form defined by $q_I(x,y) =    {1 \over N(I)} {\rm Tr}_{K/{\bf Q}}(x c(y))$, where $x \mapsto c(x)$ denotes complex conjugation; the form $q_I$
is positive definite. We have $q_I(x,y) \in {\bf Z}$ for all $x,y \in I$, and $q_I(x,x) \in 2 {\bf Z}$ for all $x \in I$. 
The complex conjugation $c: I \to I$ is an isometry of $q_I$ with characteristic polynomial $X^2-1$.

\medskip  In case (a), set 
$\epsilon_p = \partial_p[K,-q_I, c]$ in $W({\bf F}_p,id)$ if $p$ divides $S(1)$, and in 
$W({\bf F}_p,-id)$ if $p$ divides $S(-1)$.

\medskip
In case (b), set 
$\epsilon_p = \partial_p[K,q_I, c]$ in $W({\bf F}_p,id)$ if $p$ divides $S(1)$, and in 
$W({\bf F}_p,-id)$ if $p$ divides $S(-1)$.

\medskip
Let $E = {\bf Q}[X]/(S)$, let $\alpha$ be the image of $X$ in $E$, and let $\sigma : E \to E$ be the involution induced by $X \mapsto X^{-1}$.
Let $E_0$ be the fixed field of $\sigma$ in $E$, and let $d \in (E_0)^{\times}$ be such that
$E = E_0(\sqrt d)$. 
If $p$ is a prime number, set $E_p = E \otimes_{\bf Q} {\bf Q_p}$.
For all prime numbers $p$ dividing $S(1)S(-1)$, there exists $\lambda_p \in T(E_p,\sigma)$ such
that $\partial_p[E_p,b_{\lambda_p},\alpha] = - \epsilon_p$; if $p$ does not divide 
$S(1)S(-1)$, we take $\lambda_p \in T(E_p,\sigma)$ such that $(E_p,b_{\lambda_p},\alpha)$ contains
an even unimodular lattice; this is possible by \cite{BT}, \S 6. 
Set $V = E \oplus K$ and $V_p = V \otimes_{\bf Q} {\bf Q}_p$.

\medskip In case (a), set 
$q_p = b_{\lambda_p} \oplus - q_I$; and in case (b), set $q_p = b_{\lambda_p} \oplus q_I$.
In both cases, we have $\partial_p[V_p,q_p,\alpha \oplus c] = 0$
for all prime numbers $p$. 

\medskip
 In case (a), let us choose $\lambda_{\infty}$ such that the quadratic form $b_{\lambda_{\infty}}$
has signature $(3,17)$ and that the signature map of $(b_{\lambda_{\infty}},\alpha)$ is 
equal to $\tau_{S,\delta}$. 
In case (b), we choose $\lambda_{\infty}$ such that the quadratic form $b_{\lambda_{\infty}}$
has signature $(1,19)$.

\medskip
Let us determine the invariants of the quadratic forms $(V,q_p)$ and $(E_p,b_{\lambda_p})$. We have ${\rm det}(q_p) = -S(1)^2S(-1)^2
= -1$ in ${\bf Q}_p^{\times}/{\bf Q}_p^{\times 2}$. The Hasse-Witt invariant $w(q_p)$ at $p$ is 0 if
$p \not = 2$, and it is $1$ for $p = 2$. 
We have ${\rm det}(b_{\lambda}) = - {\rm det}(q_I) = - {\rm det}(-q_I)$, therefore we have
$w(b_{\lambda_p} \oplus -q_I) = w(b_{\lambda_p}) + w(-q_I)$ and $w(b_{\lambda_p} \oplus q_I) = w(b_{\lambda_p}) + w(q_I)$;
this implies that the Hasse-Witt invariants of $(E_p,b_{\lambda_p})$ and of $q_I$ are equal at
$p \not = 2$, and different at $2$ and at infinity; since $q_I$ is a global form, the sum of its Hasse-Witt
invariants is 0, and therefore the same thing is true for the Hasse-Witt invariants of $(E_p,b_{\lambda_p})$.

\medskip
We have $w(b_{\lambda_p}) = w(b_1) + {\rm cor}_{E_p/{\bf Q}_p}
(\lambda_p,d)$ for all $p$ (see Proposition \ref{dim det w f}); since $b_1$ is a global form, the above argument shows that
$\underset {v \in \mathcal V} \sum  {\rm cor}_{E_v/{\bf Q}_v}(\lambda_p,d) = 0$.
By Theorem \ref{realization} this implies that there exists $\lambda \in E_0^{\times}$ such
that 
$(E,b_{\lambda}) \otimes_{\bf Q} {\bf Q}_p \simeq (E_p,b_{\lambda_p})$ for all $p$, and that
in case (a),
$(E,b_{\lambda},\alpha)$ has signature map $\tau_{S,\delta}$;
in case (b), 
$(E,b_{\lambda})$ has signature $(1,19)$.
Set
$(V,q,t) = (E,b_{\lambda},\alpha) \oplus (K,-q_I,c)$ in case (a), and 
$(V,q,t) = (E,b_{\lambda},\alpha) \oplus (K,q_I,c)$ in case (b). 

\medskip The quadratic form $(V,q)$ contains an
even unimodular lattice stabilized by $t$ everywhere locally, and the intersection of these lattices
is an even unimodular lattice stabilized by the isometry $t$. 
In case (a), the isometry $t$ is realizable by construction.



\begin{coro}\label{coro relatively prime} Suppose that $S(1)$ and $S(-1)$ are odd, square-free,  and relatively prime. Then $(\alpha,\delta)$ is realizable. 

\end{coro}

\noindent
{\bf Proof.} This follows from Theorem \ref{relatively prime} and Corollary \ref{Mc cor}.

\begin{example}\label{20 Mc}  Let
$$S(X) = X^{20} - X^{19} - X^{15} + X^{14} - X^{11} + X^{10} - X^9 + X^6 - X^5  - X +1.$$

The corresponding Salem number is $\alpha = \lambda_{20} = 1,2326135486...$, the smallest  degree $20$ Salem
number (cf. the table in \cite{Bo}). We have $S(1) = -1$ and $S(-1) = 11$, hence Corollary \ref{coro relatively prime} implies that 
 $(\alpha,\delta)$ is realizable for all roots $\delta$ of $S$ with $|\delta| = 1$; McMullen proved that $\lambda_{20}$ is not realized by an automorphism of a
 projective $K3$ surface (see \cite{Mc3}, \S 9).

\end{example}

\begin{example}\label{5,7} Let
$$S(X) = X^{20} - X^{17} - X^{16} - X^{15} - X^{11} - X^{10} - X^9 - X^5 - X^4 -X^3 +1.$$

The corresponding Salem number is $\alpha = 1,37892866...$.

\medskip We have $S(1) = -7$ and $S(-1) = 5$; by Corollary \ref{coro relatively prime}, the pair $(\alpha,\delta)$
is realizable for all roots $\delta$ of $S$ with $|\delta| = 1$.

\end{example}

\begin{remark} The method of Theorem \ref{relatively prime} and Corollary \ref{coro relatively prime} can be used to prove
that if $S(1)$ and $S(-1)$ are odd, then $(\alpha,\delta)$ is realizable.

\end{remark}

\section{Finiteness results}\label{finiteness}

McMullen proved that up to isomorphism, there are only countably many pairs $(\mathcal X,T)$ where $T : \mathcal X \to \mathcal X$ is an
automorphism of a $K3$ surface and $\omega(T)$ is not a root of unity (see \cite{Mc1}, Theorem 3.6); we record here some related finiteness results. 


\medskip
Let $\alpha$ be a Salem number with Salem polynomial $S$, let $d = {\rm deg}(S)$, and let $\delta$ be a root
of $S$ such that $|\delta| = 1$. 

\begin{prop}\label{finite}  There exist only finitely many isomorphism classes of $K3$ surfaces
having an automorphism $T$ with $\lambda(T) = \alpha$ and $\omega(T) = \delta$.
 
 \end{prop}

We start with a definition :

\begin{defn} Let $(L,q)$ be a lattice, and let $t : L \to L $ be an isometry of $L = (L,q)$; we
say that $(L,t) = (L,q,t)$ is an {\it isometric structure}. Two isometric structures $(L,q,t)$ and $(L',q',t')$
are {\it isomorphic} if there exists an isomorphism $f : L \to L$ such that $q'(f(x),f(y)) = q(x,y)$ and
$f \circ t = t' \circ f$. An isometric structure $(L,q,t)$ is said to be {\it realizable} if $t : L \to L$
is a realizable isometry of $(L,q)$. 

\end{defn}

\begin{example} Let $T : \mathcal X \to \mathcal X$ be an automorphism of a $K3$ surface; then
$(H^2(\mathcal X,{\bf Z}),T^*)$ is a realizable isometric structure.

\end{example}

\begin{prop}\label{EF} There exist only finitely many isomorphism classes of isometric structures
$(H^2(\mathcal X,{\bf Z}),T^*)$, where $T : \mathcal X \to \mathcal X$ is an automorphism
of the $K3$ surface $\mathcal X$ with $\lambda(T) = \alpha$ and $\omega(T) = \delta$.
 
 \end{prop}

 \noindent
 {\bf Proof.}  If $T : \mathcal X \to \mathcal X$ is an automorphism of a $K3$ surface with
 $\lambda(T) = \alpha$ and $\omega(T) = \delta$, then the characteristic polynomial of
 $T^*$ is equal to $S C$, where $C$ is a product of cyclotomic polynomials with ${\rm deg}(C) = 22-d$.
 There are 
  only finitely many cyclotomic polynomials of given degree, therefore it suffices to show the finiteness
  of the number of isometric structures 
  $(H^2(\mathcal X,{\bf Z}),T^*)$ with fixed characteristic polynomial. 
  Since $\omega(T)$ is a root of $S$, the isometry $T^*$ is semi-simple; therefore
  Proposition 4 of \cite{EF} implies 
  the desired finiteness result.

  \begin{prop}\label{classification}  Let $T: \mathcal X \to \mathcal X$ and $T' : \mathcal X' \to \mathcal X'$ be
  two automorphisms of $K3$ surfaces with 
  $\omega(T) = \omega(T') = \delta$. If the isometric structures ($H^2(\mathcal X,{\bf Z}),T^*)$
  and $(H^2(\mathcal X,{\bf Z}),(T')^*)$
  are isomorphic, then the $K3$ surfaces $\mathcal X$
  and $\mathcal X'$ are isomorphic. 
  
  \end{prop}

  \noindent {\bf Proof.} Let $L =  H^2(\mathcal X,{\bf Z})$, $L' =  H^2(\mathcal X',{\bf Z})$,
  $t = T^*$ and $t' = (T')^*$, and let $f : L \to L'$ be an isomorphism of the isometric structures
  $(L,t)$ and $(L',t')$. 
  
  \medskip Note that $H^{2,0}(\mathcal X)$ is the eigenspace of $\delta$ in $L \otimes_{\bf Z}{\bf C}$; this follows
  from the fact that $\delta$ is a simple root of $S$. 
  The isomorphism $f$ respects the eigenspaces, hence 
  $f(H^{2,0}(\mathcal X)) = H^{2,0}(\mathcal X')$; this implies that $f$ respects the Hodge
  decompositions. Therefore by Theorem \ref{torelli} the $K3$ surfaces $\mathcal X$ and
  $\mathcal X'$ are isomorphic. 
  
  \medskip
  \noindent
  {\bf Proof of Proposition \ref{finite}.} This follows from 
  Proposition \ref{EF} and Proposition \ref{classification}.

\begin{coro}\label{coro finite} Let $\alpha$ be a Salem number. Then there exist at most finitely many isomorphism classes of $K3$ surfaces
 of algebraic dimension $0$ having an automorphism of dynamical degree $\alpha$.
 
 \end{coro}
 
 \noindent
 {\bf Proof of Corollary \ref{coro finite}.} This follows from Proposition \ref{finite}, since $\lambda(T)$ and
 $\omega(T)$ are roots of the same Salem polynomial if $T$ is an automorphism of a $K3$ surface of
 algebraic dimension $0$.

\section {Non-realizable Salem numbers}\label{non-realizable section}

This section and the next one contain some remarks on automorphisms of projective $K3$ surfaces. 

\medskip If $T : \mathcal X \to \mathcal X$  is an automorphism
of a $K3$ surface, recall that we denote by $\lambda(T)$ the dynamical degree of $T$, and 
by $\omega(T)$ the complex number such that $T^* : H^{2,0}(\mathcal X) \to H^{2,0}(\mathcal X)$ is
multiplication by $\omega(T)$;
a pair $(\alpha,\delta)$ is  said to be realizable if 
there exists an automorphism $T$ of a $K3$ surface with $\lambda(T) = \alpha$ and
$\omega(T) = \delta$ (see Definition \ref{realizable definition}); we say that $\alpha$ is realizable if there exists a complex number $\delta$ with $|\delta| = 1$ such that $(\alpha,\delta)$ is
realizable. 
The previous sections were concerned with realizability of pairs $(\alpha,\delta)$, such that $\alpha$ and $\delta$ are roots of the same
Salem polynomial; this implies that the $K3$ surface has algebraic dimension $0$.

\medskip 
This section takes up this problem for projective $K3$ surfaces. Recall that if $T : \mathcal X \to \mathcal X$  is an automorphism
of a projective $K3$ surface, then $\omega(T)$ is a root of unity (see \cite{Mc1}, Theorem 3.5). Moreover, the characteristic
polynomial of the restriction of $T^*$ to ${\rm Trans}(\mathcal X)$ is a power of a cyclotomic polynomial (see for instance Lemma \ref{transcendental}). 


\medskip
{\it Degree $18$ Salem numbers}

\medskip The first result concerns degree $18$ Salem numbers : we give an example of a Salem number $\alpha$ of degree $18$ 
such that $\alpha$ is {\it not realizable}, in other words, there does not exist any automorphism of a $K3$ surface, projective or not, with
dynamical degree $\alpha$.
We start with a lemma.

\begin{lemma}\label{splitting} Let $S \in {\bf Z}[X]$ be a monic polynomial of degree $18$ with $|S(1)S(-1)| = 1$, and  such that ${\rm Res}(S,\Phi_m) =1$ for $m = 3,4,6$ and $12$.
Let $C$ be a product of cyclotomic polynomials with ${\rm deg}(C) = 4$,  let
$L$ be an even, unimodular lattice and let $t : L \to L$ be an isometry with characteristic polynomial $SC$.
Then $L \simeq L_1 \oplus L_2$, where $L_1$ and $L_2$ are even, unimodular lattices stables by $t$, and
the characteristic polynomial of the restriction of $t$ to $L_1$ is $S$, the the characteristic polynomial of the restriction of $t$ to $L_2$ is $C$.

\end{lemma}

\noindent
{\bf Proof.} Set $F = SC$; the polynomial $F$ satisfies condition (C 1), hence
$|F(1)|$ and $F(-1)$ are squares. Since  $|S(1)S(-1)| = 1$, this implies that $C(1)$ and $C(-1)$ are squares; hence $C$ cannot be equal to $\Phi_8$. Therefore
all the factors of $C$ are equal to $\Phi_m$, for $m = 1,2,3,4,6$ or $12$. All these polynomials are relatively prime to $S$ over $\bf Z$, therefore $L = L_1 \oplus L_2$,
where $L_1$ and $L_2$ are even, unimodular lattices stables by $t$, and the restriction of $t$ to $L_1$ has characteristic polynomial $S$, the characteristic polynomial of the
restriction of $t$ to $L_2$ is equal to $C$. 

\begin{prop}\label{no 18} Let $S$ be a Salem polynomial of degree $18$, and let $\alpha$ be the corresponding Salem number.  Assume that $|S(1)S(-1)| = 1$, and that ${\rm Res}(S,\Phi_m) =1$ for $m = 3,4,6$ and $12$. Then $\alpha$ is not realizable. 

\end{prop}

\noindent{\bf Proof.} By Proposition \ref{resultant} we already know that $(\alpha,\delta)$ is not realizable for any root $\delta$ of $S$. Suppose that $(\alpha,\zeta)$ is
realizable for some root of unity $\zeta$. Let $T : \mathcal X \to \mathcal X$  be an automorphism with $\lambda(T) = \alpha$ and $\omega(T) = \zeta$; note that
$\mathcal X$ is then a projective $K3$ surface. Set
$L = H^2(\mathcal X,{\bf Z})$, and $t = T^*$. Let $F$ be the characteristic polynomial of $t$; we have 
 $F = SC$, where $C$ is a
power of a cyclotomic polynomial with $C(\zeta) = 0$. By Lemma \ref{splitting}, we have $L \simeq L_1 \oplus L_2$, where $L_1$ and $L_2$ are even, unimodular lattices stables by $t$, and
the characteristic polynomial of the restriction of $t$ to $L_1$ is $S$, the the characteristic polynomial of the restriction of $t$ to $L_2$ is $C$.
Since $\omega(T) = \zeta$, the lattice $L_2$ is isomorphic to ${\rm Trans}(\mathcal X)$, and its signature is $(2,2)$; hence the lattice
$L_1$ is isomorphic to ${\rm Pic}(\mathcal X)$, and its signature is $(1,17)$. 
Therefore ${\rm Pic}(\mathcal X)$ is isomorphic to $\Lambda_{1,17}$, hence ${\rm Aut}(\mathcal X)$ is finite
(see \cite{PS}, \S 7, Example 2). This implies that all the automorphisms of $\mathcal X$ have dynamical degree 1, contradicting the assumption that
$\lambda(T) > 1$; hence $(\alpha,\zeta)$ is not realizable, as claimed. 

\begin{example}\label{18 example} I thank Chris Smyth for this example. Let $a$ be an integer, let $R_a(X) = X^2(X^2-4)(X^2-3)(X^2-1)(X-a) - 1$, and set $S_a(X) = X^9 R_a(X + X^{-1})$.
For $a = 3$, the polynomial $S = S_a$ is a Salem polynomial, and the corresponding Salem number is $2,618575...$
We have  $|S(1)S(-1)| = 1$ and  ${\rm Res}(S,\Phi_m) =1$ for $m = 3,4,6$ and $12$, hence by Proposition \ref{no 18}, the Salem number $\alpha = 2,618575...$ is not realizable.

\end{example}

 {\it Degree $20$ Salem numbers}

\medskip We now give examples of degree $20$ Salem numbers that are not realized by automorphisms of {\it projective} $K3$ surfaces.

\begin{prop}\label{projective 20} Let $S$ be a Salem polynomial of degree $20$, and let $\alpha$ be the corresponding Salem number.
Let $T : \mathcal X \to \mathcal X$  is an automorphism
of a projective $K3$ surface such that $\lambda(T) = \alpha$. Then one of the following holds :

\medskip
{\rm (i)} $S(-1)$ is a square.

\medskip
{\rm (ii)} $|S(1)|$ is a square.

\medskip
{\rm (i)} There exist integers $m$ and $n$ such that $|S(1)| = 2m^2$ and $S(-1) = 2n^2$.

\end{prop}

\noindent{\bf Proof.} Let $F$ be the characteristic polynomial of $T^*$; then $F = S C$, where $C$ is a power of a cyclotomic polynomial. Therefore
we have either $C(X) = (X-1)^2$, $C(X) = (X+1)^2$, or $C = \Phi_m$, for $m = 3,4$ or $6$. The polynomial $F$ satisfies condition (C 1), hence
(i) holds if $C(X) = (X-1)^2$ or $C = \Phi_3$, (ii) holds if $C(X) = (X+1)^2$ or $C = \Phi_6$, and (iii) holds if $C = \Phi_4$. 

\begin{coro}\label{no 20} Let $p$ and $q$ be two distinct prime numbers, let $S$ be a Salem polynomial of degree $20$ such that $S(1) = -p$ and $S(-1) = q$, and let
$\alpha$ be the corresponding Salem number. Then $\alpha$ is not realizable by any automorphism of a projective $K3$ surface.

\end{coro} 

\noindent
{\bf Proof.} This follows from Proposition \ref{projective 20}, since $S$ does not satisfy any of the conditions (i), (ii) or (iii). 

\begin{example} It is easy to find Salem numbers satisfying the hypothesis of Corollary \ref{no 20} : for instance, $\alpha = 1,378928...$ (with $S(1) = -7$, $S(-1) = 5$),
 $\alpha = 1,394338...$ (with $S(1) = -5$, $S(-1) = 3$),  $\alpha = 1,464501...$ (with $S(1) = -13$, $S(-1) = 3$),  $\alpha = 1,464843...$ (with $S(1) = -7$, $S(-1) = 29$),...
 By Corollary \ref{no 20}, these Salem numbers are not realized by any automorphism of a projective $K3$ surface; note that they 
 are realized by automorphisms of non-projective $K3$ surfaces (see Corollary \ref{coro relatively prime}, or Takada \cite{T}, Theorem 1.3).

\end{example}

\begin{remark} Let $\alpha = \lambda_{20} = 1,2326135486...$ be the smallest degree $20$ Salem number, and let $S$ be the corresponding Salem polynomial.
We have $S(1) = -1$ and $S(-1) = 11$, hence $S$ satisfies condition (ii) of Proposition \ref{projective 20}. Nevertheless, $\alpha$ is not realized by any automorphism
of a projective $K3$ surface, as shown by McMullen in \cite{Mc3}, \S 9. 

\end{remark}

\section {Automorphisms which act trivially on Picard groups}\label{Kondo section}

In this section, we consider automorphisms of projective $K3$ surfaces that induce the identity on their Picard groups - in particular, they have dynamical degree one. 
The aim is to give an alternative approach to some results of Vorontsov \cite {V}, Kondo \cite{Ko}, Machida-Oguiso \cite{MO} and Oguiso-Zhang \cite{O 0}; see also Brandhorst \cite{Bra} and the references therein for more
recent results. 


\medskip

Let $C$ be a power of a cyclotomic polynomial of even degree, i.e $C = \Phi_m^n$, where $m \geqslant 3$ and $n \geqslant 1$ are integers, and assume that
$2 \leqslant {\rm deg}(C) \leqslant 20$. 

\medskip

\begin{lemma}\label{power of cyclotomic} Let $T : \mathcal X \to \mathcal X$ be an automorphism of a projective $K3$ surface such that $T^*|{\rm Pic}(\mathcal X)$ is the identity, and that the
characteristic polynomial of $T^*|{\rm Trans}(\mathcal X)$ is equal to $C$. Then the following conditions hold

\medskip
{\rm (i)} $C(-1)$ is a square.

\medskip
{\rm (ii)} If $C(1) = 1$, then  ${\rm deg}(C)  \equiv \ 4  \ {\rm (mod \ 8)}$.

\medskip
Moreover, if  $C(1) = 1$, then the lattice ${\rm Trans}(\mathcal X)$ is unimodular.


\end{lemma} 

\noindent
{\bf Proof.} Set $F(X) = C(X) (X-1)^{22-{\rm deg}(C)}$; the polynomial  $F$ is the characteristic 
polynomial of $T^*$.
Since $T^*$ is an isometry of a unimodular lattice, the polynomial $F$ satisfies condition (C 1) (see \S \ref{local}), and
hence $F(-1)$ is a square; this implies that $C(-1)$ is a square; therefore
condition (i) holds. 

\medskip Set $L = H^2(\mathcal X,{\bf Z})$, set $L_1 = {\rm Ker}(C(T^*))$, and let $L_2$ be the sublattice of fixed points of $T^*$.
If $C(1) = 1$, then the polynomials $X-1$ and $C$ are relatively prime over $\bf Z$; this implies that 
$L = L_1 \oplus L_2$, hence $L_1$ and $L_2$ are both unimodular. We have $L_2 = {\rm Pic}(\mathcal X)$ and $L_1 = {\rm Trans}(\mathcal X)$; this shows that
the lattice ${\rm Trans}(\mathcal X)$ is unimodular. 
The
signature of the lattice $L_1$ is $(2,{\rm deg}(C)-2)$; this lattice is unimodular, therefore ${\rm deg}(C) - 4  \equiv \ 0  \ {\rm (mod \ 8)}$, and this implies that
${\rm deg}(C)  \equiv \ 4  \ {\rm (mod \ 8)}$, hence (ii) holds.

\medskip

\begin{prop}\label{irreducible} Suppose that $C$ is a cyclotomic polynomial. There exists an automorphism 
$T : \mathcal X \to \mathcal X$ of a projective $K3$ surface such that  $T^*|{\rm Pic}(\mathcal X)$ is the identity and that the
characteristic polynomial of $T^*|{\rm Trans}(\mathcal X)$ is equal to $C$ if and only if the following conditions hold :

\medskip
{\rm (i)} $C(-1) = 1$.

\medskip
{\rm (ii)} If  $C(1) = 1$, then  ${\rm deg}(C)  \equiv \ 4  \ {\rm (mod \ 8)}$.

\medskip
The $K3$ surface $\mathcal X$ with this property is unique up to isomorphism. 
Moreover, the lattice ${\rm Trans}(\mathcal X)$ is unimodular if and only if  $C(1)  = 1$.

\end{prop}

\noindent
{\bf Proof.} Let $T : \mathcal X \to \mathcal X$ be an automorphism of a projective $K3$ surface such that $T^*|{\rm Pic}(\mathcal X)$ is the identity, and that the
characteristic polynomial of $T^*|{\rm Trans}(\mathcal X)$ is equal to $C$. Condition (ii) follows from part (ii) of Lemma \ref {power of cyclotomic}.
Since $C$ is a cyclotomic polynomial, $C(-1)$ is a square if and only if $C(-1) = 1$, therefore part (i) of Lemma  \ref {power of cyclotomic} implies (i).

\medskip
Conversely, assume that conditions (i) and (ii) hold, and let us prove the existence of $T : \mathcal X \to \mathcal X$ with the required properties. 
Set $F(X) = C(X) (X-1)^{22-{\rm deg}(C)}$;
by (i), the polynomial $F$ satisfies condition (C 1). 

\medskip
If $C(1) > 1$ and ${\rm deg}(C) < 20$, then $G_F = 0$; if ${\rm deg}(C) = 20$ and $C(1)$ is not a square,
then $G_F (D_+)  = 0$ (see \S \ref{combinatorial section} and \S \ref{sufficient section}). By Corollary \ref{G=0} there exists an even, unimodular lattice $L$ of signature $(3,19)$ and an isometry $t : L \to L$ with characteristic polynomial $F$ and signature
map $\tau$ satisfying $\tau(C) = (2,{\rm deg}(C) - 2)$.  Let $L_1 = {\rm Ker}(C(t))$ and let $L_2$ be the sublattice of $L$ of fixed points by $t$. 

\medskip Suppose now that $C(1)$ is a square; since $C$ is a cyclotomic polynomial, this implies that $C(1) = 1$; condition (ii) implies that ${\rm deg}(C)  \equiv \ 4  \ {\rm (mod \ 8)}$.
 Hence the polynomial $C$ satisfies Condition (C 1), and by   \cite{BT}, Theorem A (or Corollary \ref{G=0})  there exists an even, unimodular lattice $L_1$ of signature $(2,{\rm deg}(C)-2)$ having
an isometry $t_1 : L_1 \to L_1$ of characteristic polynomial $C$. Let $L_2$ be an even, unimodular lattice of signature $(1,22-{\rm deg}(C))$. Then $L = L_1 \oplus L_2$
is isomorphic to $\Lambda_{3,19}$. Let $t : L \to L$ be the isometry that is equal to $t_1$ on $L_1$, and is the identity on $L_2$.

\medskip
The  isometry $t : L \to L$ is realizable (see McMullen, \cite{Mc3}, Theorem 6.1), and 
gives rise to an automorphism $T : \mathcal X \to \mathcal X$ of a projective $K3$ surface such that ${\rm Pic}(\mathcal X)$ is isomorphic to $L_2$, and
${\rm Trans}(\mathcal X)$ is isomorphic to $L_1$. 

\medskip
If $C(1) = 1$, then the lattice ${\rm Trans}(\mathcal X)$ is unimodular by Lemma \ref{power of cyclotomic}.
Conversely, suppose that ${\rm Trans}(\mathcal X)$ is unimodular. Its signature is $(2,{\rm deg}(C) - 2)$, and this implies that
${\rm deg}(C)  \equiv \ 4  \ {\rm (mod \ 8)}$. Moreover, $C$ satisfies condition (C 1), hence $C(1)$ and $C(-1)$ are squares; since $C$ is cyclotomic, this
implies that $C(1) = C(-1) = 1$.

\medskip
The uniqueness of the $K3$ surface up to isomorphism follows from a result of Brandhorst, \cite{Bra}, Theorem 1.2.

\bigskip We now show that Proposition \ref{irreducible} implies some results of Vorontsov, Kondo, Machida-Oguiso and Oguiso-Zhang. If $\mathcal X$ is a $K3$ surface, we denote by
$H_{\mathcal X}$ the kernel of the map ${\rm Aut}({\mathcal X}) \to {\rm O}({\rm Pic}{\mathcal X})$; it is a cyclic group of finite order (see \cite{N}, Corollary 3.3). Let $m_{\mathcal X}$ be the order of $H_{\mathcal X}$. Set 
$\Sigma = \{12, \ 28, \ 36, \ 42, \ 44, \ 66 \}$ and $\Omega = \{ 3, \ 9, \ 27, \ 5, \ 25, \ 7, \ 11, \ 13, \ 17, \ 19 \}$. The following corollaries are a combination of results of 
\cite{V}, \cite{Ko}, \cite{MO} and \cite{O 0} :

\begin{coro}\label{unimodular}
Let $\mathcal X$ be a projective $K3$ surface, and suppose that the rank of ${\rm Trans}({\mathcal X})$ is equal to $\varphi(m_{\mathcal X})$. Assume that
the lattice ${\rm Trans}({\mathcal X})$ is unimodular. Then $m_{\mathcal X} \in \Sigma$. 
Conversely, for each $m \in \Sigma$, there exists a unique {\rm (up to isomorphism)} projective
$K3$ surface with $m_{\mathcal X} = m$  
and ${\rm rank}({\rm Trans}({\mathcal X})) = \varphi(m)$.
 
 \end{coro}
 
 \begin{coro}\label{nonunimodular}
Let $\mathcal X$ be a projective $K3$ surface, and suppose that the rank of ${\rm Trans}({\mathcal X})$ is equal to $\varphi(m_{\mathcal X})$. 
Assume that
the lattice ${\rm Trans}({\mathcal X})$ is not unimodular. Then $m_{\mathcal X} \in \Omega$. Conversely, for each $m \in \Omega$, there exists a unique
{\rm (up to isomorphism)} projective
$K3$ surface with $m_{\mathcal X} = m$ and ${\rm rank}({\rm Trans}({\mathcal X})) = \varphi(m)$. 
 
 \end{coro}
 
 \noindent{\bf Proof of Corollaries \ref{unimodular} and \ref{nonunimodular}.} Let $T$ be a generator of the cyclic group $H_{\mathcal X}$, and note that $m_{\mathcal X} = m$ if and only if the characteristic polynomial of the restriction of $T^*$ to ${\rm Trans}({\mathcal X})$ is equal to $\Phi_m$. 
 
 \medskip Suppose that the the characteristic polynomial of the restriction of $T^*$ to ${\rm Trans}({\mathcal X})$ is equal to $\Phi_m$. Since the rank of ${\rm Trans}({\mathcal X})$ is at most $20$,
 we have ${\rm deg}(\Phi_m) \leqslant 20$. By Proposition \ref{irreducible}, $C = \Phi_m$ satisfies conditions (i) and (ii), and this implies that $m \in \Sigma$ if $\Phi_m(1)  = 1$, and 
 $m \in \Omega$ otherwise; therefore $m_{\mathcal X} \in \Sigma$ if ${\rm Trans}({\mathcal X})$ 
is unimodular, and $m_{\mathcal X} \in \Omega$ otherwise. 

\medskip Conversely, if  $m \in \Sigma$, then $\Phi_m(1) = \Phi_m(-1) = 1$ and 
${\rm deg}(\Phi_m)   \equiv \ 4  \ {\rm (mod \ 8)}$; if $m \in \Omega$, then $\Phi_m(-1) = 1$ and $\Phi_m(1) > 1$, hence Proposition \ref{irreducible} implies
the desired result.

\bigskip
\bigskip
Eva Bayer--Fluckiger 

EPFL-FSB-MATH

Station 8

1015 Lausanne, Switzerland

\medskip

eva.bayer@epfl.ch


\bigskip






\begin{thebibliography}{99}


\bibitem[BB 92] {BB} C.  Bachoc, C. Batut, \textit{\'Etude algorithmique de réseaux construits avec la forme trace}, Experiment. Math. {\bf 1} (1992), 183-190.

 \bibitem[BHPV 04]{BHPV} W. Barth, K. Hulek, C. Peters, A. Van de Ven, \textit{Compact complex surfaces},
 Second edition. Ergebnisse der Mathematik und ihrer Grenzgebiete , 3. Folge, {\bf 4} (2004).

\bibitem[B 84]{B 84} E. Bayer-Fluckiger, \textit{Definite unimodular lattices having an automorphism of
given characteristic polynomial}, Comment. Math. Helv. {\bf 59} (1984), 509-538.



 
 
 \bibitem[B 99]{B 99} E. Bayer-Fluckiger,
\textit{Lattices and number fields}, Algebraic geometry : Hirzebruch 70 (Warsaw, 1998), Contemp. Math.
{\bf 241}, Amer. Math. Soc., Providence, RI (1999), 69-84.



\bibitem[B 15]{B 15} E. Bayer-Fluckiger,
\textit{Isometries of quadratic spaces}, J. Eur. Math. Soc.  \textbf{17} (2015), 1629-1656.

\bibitem[B 21]{B 21} E. Bayer-Fluckiger,
\textit{Isometries of lattices and Hasse principles}, J. Eur. Math. Soc. 
online first (2023), DOI 10.4171/jems/1334.

\bibitem[BM 94]{BM} E.  Bayer-Fluckiger, J. Martinet, \textit{Formes quadratiques liées aux algèbres semi-simples} J. reine angew. Math. {\bf 451} (1994), 51-69.

 \bibitem[BM 79]{EF} E. Bayer, F. Michel, \textit{Finitude du nombre de classes d'isomorphisme des
 structures isom\'etriques enti\`eres}, Comment. Math. Helv. {\bf 54} (1979), 378-396. 



\bibitem[BT 20]{BT} E. Bayer-Fluckiger, L. Taelman,
\textit{Automorphisms of even unimodular lattices and equivariant Witt groups}, J. Eur. Math. Soc.
\textbf{22} (2020), 3467-3490.

\bibitem[BG 05]{bg} M. Belolipetsky, W. T. Gan, \textit{The mass of unimodular lattices},
J. Number Theory {\bf 114} (2005), 221-237.

 

\bibitem[Bo 77]{Bo} D. W. Boyd, \textit{Small Salem numbers}, Duke Math. J. {\bf 44} (1977), 315-328.

 
\bibitem[Br 19]{Bra} S. Brandhorst, \textit{The classification of purely non-symplectic automorphisms of high order on $K3$ surfaces},
J. Algebra {\bf 533} (2019), 229-265.

\bibitem[Br 20]{Br} S. Brandhorst, \textit{On the stable dynamical spectrum of complex surfaces},
Math. Ann. {\bf 377} (2020), 421-434.

\bibitem[BGA 16]{BG} S. Brandhorst, V. Gonzalez-Alonso, \textit{Automorphisms of minimal entropy
on supersingular $K3$ surfaces}, J. London Math. Soc. {\bf 97} (2016), 282-305. 

\bibitem[Ca 99]{Ca 99} S. Cantat, \textit{Dynamique des automorphismes des surfaces projectives complexes},
C. R. Acad. Sci. Paris {\bf 328} (1999), 901-906. 

\bibitem[Ca 01]{Ca 01} S. Cantat, \textit {Sur la dynamique du groupe d'automorphismes des
surfaces $K3$}, Transform. Groups. {\bf 6} (2001), 201-214.

\bibitem[Ca 14]{Ca 14} S. Cantat, \textit {Dynamics of automorphisms of compact complex surfaces}, Frontiers in Complex
Dynamics : In celebration of John Milnor's 80th birthday (463-514), Princeton Math. Ser.  {\bf 51}, Princeton Univ. Press,
Princeton, NJ, 2014.


\bibitem[Co 89]{Co} D. Cox, \textit {Primes of the form $x^2 + n y^2$, Fermat, class field theory and complex multiplication},
A Wiley-Interscience Publication, John Wiley \& Sons, Inc., New York, 1989.

\bibitem[E 13] {E} W. Ebeling, \textit {Lattices and codes}, Third edition, Advanced Lectures in Mathematics, Springer Spectrum,
Wiesbaden, 2013.

\bibitem[GH 01]{GH} E. Ghate, E. Hironaka, \textit{The arithmetic and geometry of Salem numbers},
Bull. Amer. Math. Soc. {\bf 38} (2001), 293-314.




\bibitem[GM 02]{GM} B. Gross, C. McMullen,
\textit{Automorphisms of even, unimodular lattices and unramified Salem numbers}, J. Algebra {\bf 257} (2002), 265--290. 

\bibitem[H 16]{H} D. Huybrechts, \textit{Lectures on $K3$ surfaces}, Cambridge Studies in Advanced Mathematics {158}, Cambridge
University Press, Cambridge, 2016.



\bibitem[IT 21]{IT} K. Iwasaki, Y. Takada, \textit{$K3$ surfaces, Picard numbers and Siegel disks}, J. Pure Appl. Algebra,
{\bf 227} (2023), Paper No. 107215, 31 pages. 

\bibitem[K 92]{Ko} S. Kondo, \textit{Automorphisms of algebraic $K3$ surfaces which act trivially on Picard groups}, J. Math. Soc. Japan {\bf 44} (1992), 
75-98.

\bibitem[K 20]{K} S. Kondo, \textit{$K3$ surfaces}, Tracts in Mathematics {\bf 32}, European Mathematical
Society (2020). 

\bibitem[MO 98]{MO} N. Machida; Oguiso, K, \textit{On $K3$ surfaces admitting finite non-symplectic group actions}, J. Math. Sci. Univ. Tokio {\bf 5} (1998), 273-297.


 
\bibitem [McM 02]{Mc1} C. McMullen, Dynamics on K3 surfaces: Salem numbers and Siegel disks. J. Reine Angew. Math. {\bf 545} (2002), 201-233.

\bibitem[McM 11]{Mc2} C. McMullen, \textit K3 surfaces, entropy and glue. J. Reine Angew. Math. {\bf 658} (2011), 1-25. 

\bibitem[McM 16]{Mc3} C. McMullen, \textit
Automorphisms of projective K3 surfaces with minimum entropy. Invent. Math. {\bf 203} (2016), 179-215.

\bibitem[M 69]{M} J. Milnor, \textit {Isometries of inner product spaces}, Invent. Math. {\bf  8}  (1969), 83--97. 

\bibitem[N 79]{N} V.V. Nikulin, \textit {Finite groups of automorphisms of K\"ahlerian $K3$ surfaces}, Trudy Moskov Mat. Obshch.
{\bf 38}  (1979), 75-137.



\bibitem[O 07]{O 7} K. Oguiso, \textit {Automorphisms of hyperkähler manifolds in the view of topological entropy}, Algebraic geometry, 173-185, Contemp. Math {\bf 422}, Amer. Math. Soc., Providence, RI (2007).

\bibitem[O 08]{O 8} K. Oguiso, \textit{ Bimeromorphic automorphism groups of non-projective hyperk\"ahler manifolds—a note inspired by C. T. McMullen},  J. Differential Geom. {\bf 78}  (2008), 163-191.


\bibitem[O 10]{O} K. Oguiso, \textit {The third smallest Salem number in automorphisms of $K3$ surfaces}, Algebraic
geometry in East Asia - Seoul 2008, 331-360, Adv. Stud. Pure Math. {\bf 60}, Math. Soc. Japan, Tokyo, 2010.

\bibitem[OZ 00]{O 0} K. Oguiso, D-Q. Zhang, \textit {On Vorontsov's theorem on $K3$ surfaces with non-symplectic group actions}, Proc. Amer. Math. Soc. {\bf 128} (2000) 1571-1580.

\bibitem[PS 71]{PS} I.I. Piatetski-Shapiro, I.R. Shafarevich, \textit {A Torelli theorem for algebraic surfaces of type $K3$}, Math. USSR Izvestija {\bf 5} (1971), 547-588.

\bibitem[R 12]{R} P. Reschke, \textit{Salem numbers and automorphisms of complex surfaces}, Math. Res. Lett. {\bf 19}
(2012), 475-482.


\bibitem[Sa 45]{Sa} R. Salem, \textit{Power series with integral coefficients}, Duke Math. J. {\bf 12} (1945),
153-172.

\bibitem[Sch 85]{Sch}  W. Scharlau, \textit {Quadratic and hermitian forms}, Grundlehren der Mathematischen Wissenschaften {\bf 270}, Springer-Verlag, Berlin, 1985.

\bibitem[S 77]{S} J-P. Serre, \textit {Cours d'arithm\'etique}, Presses Universitaires de France, 1977.

\bibitem [Sm 15]{Sm} C. Smyth, \textit {Seventy years of Salem numbers}, Bull. Lond. Math. Soc.
{\bf 47} (2015), 379-395. 

\bibitem[T 22]{T} Y. Takada, \textit{Lattice isometries and $K3$ surface automorphisms : Salem numbers of degree $20$}, J. Number Th. {\bf 252} (2023), 195-242,
arXiv:2210.12946v1




\bibitem[V 83]{V} S.P. Vorontsov, \textit{Automorphisms of even lattices that arise in connection 
with automorphisms of algebraic $K3$ surfaces}, Vestnik Mosk. Univ. Mathematika {\bf 38} (1983), 19-21. 

\bibitem [Z 19]{Z} S. Zhao, \textit {Automorphismes loxodromiques de surfaces ab\'eliennes r\'eelles}, Ann. Fac. Sci. Toulouse Math
{\bf 28}, 109-127.




\end{thebibliography}
\end{document}